\documentclass[authoryear,preprint,12pt]{elsarticle}
\pdfoutput=1

\def\today{2nd March 2012}

\usepackage{natbib}

\usepackage{amssymb}

\usepackage{alltt,times,pifont,amssymb,url,graphicx} 
\usepackage{latexsym} 

\def\Section#1{\section{#1}}
\def\Subsection#1{\subsection{#1}}
\def\Subsubsection#1{\subsubsection{#1}}
\def\Hide#1{\relax}

\def\Lab#1{{\it(#1)}}
\DeclareSymbolFont{AMSb}{U}{msb}{m}{n}
\DeclareSymbolFontAlphabet{\mathbb}{AMSb}

\def\Func#1{{\mathsf{#1}}}

\def\fin{\mathbb{F}}
\def\zahl{\mathbb{Z}}
\def\NatSort{{\cal N}}
\def\VecSort{{\cal V}}
\def\ScaSort{{\cal R}}
\def\SetSort{{\cal P}}

\def\IP#1{\Func{IP}^{#1}}
\def\HS#1{\Func{HS}^{#1}}
\def\NS#1{\Func{NS}^{#1}}
\def\BS#1{\Func{BS}^{#1}}
\def\VS#1{\Func{VS}^{#1}}
\def\MS{\Func{MS}}
\def\NSA#1{\NS{#1}_{{+}}}
\def\BSA#1{\BS{#1}_{{+}}}
\def\MSA{\Func{MS}_{{+}}}
\def\Tfin{T_{\mbox{\scriptsize fin}}}

\def\LL{{\cal L}}
\def\A{\LL_{A}}
\def\AR{\LL_{AR}}
\def\AV{\LL_{AV}}
\def\AN{\LL_{AN}}
\def\LA{\LL^{+}}
\def\LB{\LL_{B}}
\def\LI{\LL_{I}}
\def\LM{\LL_{M}}
\def\LN{\LL_{N}}
\def\LNA{\LL^{+}_{N}}
\def\LV{\LL_{V}}

\def\sG{\mathbb{G}}
\def\sH{\mathbb{H}}
\def\sK{\mathbb{K}}
\def\sL{\mathbb{L}}
\def\sLO{\sL_0}
\def\sM{\mathbb{M}}
\def\sS{\mathbb{S}}
\def\sY{\mathbb{Y}}
\def\sX{\mathbb{X}}

\def\Abs#1{|#1|}

\def\Dim#1{\Func{dim}(#1)}

\def\Card#1{|#1|}
\def\C{\Func{c}}
\def\Dimc{\mathbb{V}^{\C}}
\def\Exp#1{\Func{exp}(#1)}

\def\Log#1{\Func{log}(#1)}
\def\Max{\Func{max}}
\def\NRM{\Func{NRM}}
\def\Seq#1{#1^{\nat}}
\def\Size#1{\#(#1)}
\def\Tuple#1{\overline{#1}}
\def\Cos#1{\Func{cos}(#1)}
\def\Sin#1{\Func{sin}(#1)}

\def\Sup{\Func{sup}}
\def\pA#1{\Func{A}(#1)}
\def\pACC#1{\Func{ACC}(#1)}
\def\pADJ#1{\Func{ADJ}(#1)}
\def\pB#1{\Func{B}(#1)}
\def\pD{\Func{D}}
\def\pDef{\Func{Def}}
\def\pE#1{\Func{E}(#1)}
\def\pEP#1{\Func{EP}(#1)}
\def\pEPX#1{\Func{EPX}(#1)}
\def\pESD#1{\Func{ESD}(#1)}
\def\pG#1{\Func{G}(#1)}
\def\pGamma#1{{\sf\Gamma}\left(#1\right)}
\def\pH#1{\Func{H}(#1)}
\def\pI#1{\Func{I}(#1)}
\def\pInf{\Func{Inf}}
\def\pMult{\Func{Mult}}
\def\pN#1{\Func{N}(#1)}
\def\pNRM#1{\NRM(#1)}
\def\pNTIMES#1{\Func{NTIMES}(#1)}
\def\pNat#1{\Func{Nat}(#1)}
\def\pO#1{\Func{O}(#1)}

\def\pP#1{\Func{P}(#1)}
\def\pPeriodic{\Func{Periodic}}
\def\pPi#1{{\sf\Pi}\left(#1\right)}
\def\pPeano{\Func{Peano}}
\def\pQTIMES#1{\Func{QTIMES}(#1)}
\def\pQ{\Func{Q}}
\def\pR#1{\Func{R}(#1)}

\def\pRTIMES#1{\Func{RTIMES}(#1)}
\def\pS{\Func{S}}
\def\pSIN#1{\Func{SIN}(#1)}
\def\pT#1{\Func{T}(#1)}
\def\pU#1{\Func{U}(#1)}
\def\pUG{\Tuple{\Func{UG}}}
\def\pU#1{\Func{U}(#1)}

\def\pW#1{\Func{W}(#1)}
\def\pX#1{\Func{X}(#1)}
\def\pXAX#1{\Func{XAX}(#1)}
\def\pYAX#1{\Func{YAX}(#1)}
\def\pZ#1{\Func{Z}(#1)}
\def\pZTIMES#1{\Func{ZTIMES}(#1)}

\def\Conv#1{\Func{conv}(#1)}
\def\Sconv#1{\Func{sconv}(#1)}

\def\fAE{\Al\Ex}
\def\fEA{\Ex\Al}
\def\fAIA{{\Al}{\Imp}{\Al}}

\def\IsDef{\mathrel{{:}{=}}}
\def\Al{\forall}
\def\Ex{\exists}
\def\Imp{\Rightarrow}
\def\Pmi{\Leftarrow}
\def\Iff{\Leftrightarrow}
\def\And{\land}
\def\Or{\lor}
\def\Not{\lnot}

\def\Res#1#2{#1{|}_{#2}}
\def\ST{\mathrel{|}}

\def\PreSpecial#1{{#1}^{\mbox{{\scriptsize\sf R}}}}
\def\QElim#1{{#1}^{\mbox{{\scriptsize\sf QE}}}}
\def\Special#1{{#1}^{\mbox{{\scriptsize\sf S}}}}
\def\Star#1{{#1}^{*}}

\def\Va{\mathbf{a}}
\def\Vb{\mathbf{b}}
\def\Vc{\mathbf{c}}
\def\Vd{\mathbf{d}}
\def\Ve{\mathbf{e}}
\def\Vf{\mathbf{f}}
\def\Vg{\mathbf{g}}

\def\Vl{\mathbf{l}}

\def\Vp{\mathbf{p}}
\def\Vq{\mathbf{q}}
\def\Vr{\mathbf{r}}
\def\Vs{\mathbf{s}}

\def\Vu{\mathbf{u}}
\def\Vv{\mathbf{v}}
\def\Vw{\mathbf{w}}
\def\Vx{\mathbf{x}}
\def\Vy{\mathbf{y}}
\def\Vz{\mathbf{z}}
\def\VO{\mathbf{0}}
\def\VQ{\mathbf{Q}}

\def\Norm#1{{\|}#1{\|}}

\def\Diff{\mathop{\backslash}}

\newtheorem{Theorem}{Theorem}
\newtheorem{Lemma}[Theorem]{Lemma}
\newtheorem{Corollary}[Theorem]{Corollary}

\def\Proof{\par \noindent{\bf Proof: }}
\def\Done{\hfill\rule{0.5em}{0.5em}}

\newcommand{\nat}{\mbox{$\protect\mathbb N$}}
\newcommand{\num}{\mbox{$\protect\mathbb Z$}}
\newcommand{\pow}{\mbox{$\protect\mathbb P$}}
\def\Pow#1{\pow(#1)}
\newcommand{\rat}{\mbox{$\protect\mathbb Q$}}
\newcommand{\real}{\mbox{$\protect\mathbb R$}}

\newcommand{\spot}{{\cdot}}

\newcommand{\all}[1]{\forall #1 \spot\:}
\newcommand{\ex}[1]{\exists #1 \spot\:}
\newcommand{\exu}[1]{\exists! #1 \spot\:}

\newcommand{\Ands}{\bigwedge}

\newcommand{\BEQ}{\mbox{\raise4pt\hbox{$\ulcorner$}}}

\newcommand{\EEQ}{\mbox{\raise4pt\hbox{$\urcorner$}}}

\let\subset=\subseteq

\newcommand{\Union}{\cup}
\newcommand{\Inter}{\cap}

\newcommand{\BA}{\begin{array}[t]{l}}
\newcommand{\EA}{\end{array}}

\newcommand{\floor}[1]{\mbox{$\left\lfloor #1 \right\rfloor$}}

\makeatletter
\def\imod#1{\allowbreak\mkern10mu({\operator@font mod}\,\,#1)}
\makeatother

\newcommand\HOLSpacing{12pt}

\newlength{\hsbw}
\newcommand{\ty}{\!:\!}

\newcommand{\inner}[1]{\mbox{$\left\langle #1 \right\rangle$}}

\journal{Annals of Pure and Applied Logic}

\begin{document}


\begin{frontmatter}



\title{Some new results on decidability for elementary algebra and geometry}

\author{Robert M. Solovay}
\ead{solovay@gmail.com}
\address{PO Box 5949, Eugene, OR 97405, USA}

\author{R. D. Arthan\corref{fn-rda}}
\ead{rda@lemma-one.com}
\address{
Lemma 1. Ltd., 27 Brook St., Twyford, Berkshire, RG1 2HX, UK \& \\
School of Electronic Engineering and Computer Science, \\
Queen Mary, University of London, London, E1 4NS, UK}

\author{John Harrison}
\ead{johnh@ichips.intel.com}
\address{Intel Corporation, RA2-451, 2501 NW 229th Avenue, Hillsboro, OR 97124, USA}

\cortext[fn-rda]{Corresponding author}

\begin{abstract}
We carry out a systematic study of decidability for theories {\em(a)} of real
vector spaces, inner product spaces, and Hilbert spaces and {\em(b)} of normed
spaces, Banach spaces and metric spaces, all formalised using
a 2-sorted first-order language.  The theories for list {\em(a)} turn out to be
decidable while the theories for list {\em(b)} are not even arithmetical:
the theory of 2-dimensional Banach spaces, for example, has the same many-one
degree as the set of truths of second-order arithmetic.

We find that the purely universal and purely existential fragments of the
theory of normed spaces are decidable, as is the $\fAE$ fragment of the theory
of metric spaces.  These results are sharp of their type: reductions of Hilbert's $10^{\mbox{\scriptsize th}}$ problem show that the $\fEA$ fragments for metric and normed
spaces and the $\fAE$ fragment for normed spaces are all undecidable.

\end{abstract}

\begin{keyword}
Decidability \sep undecidability \sep vector spaces \sep inner product spaces \sep normed spaces \sep Hilbert spaces \sep Banach spaces \sep metric spaces

\MSC 03D35 \sep 03C10 \sep 03C65

\end{keyword}

\end{frontmatter}

{\small \tableofcontents}
\Section{Introduction}\label{sec:theory}

It is natural to formulate the theory of real vector spaces using a
2-sorted first-order language with a sort for the scalars and a sort for
the vectors.
Introduction of coordinates reduces
the theory $\VS{n}$ of a vector space of a given finite dimension $n$§
to the first-order theory of the real numbers, known to be decidable since the
pioneering work of \citet{tarski-decision} to which our title alludes.
The purpose of this paper is to investigate decidability for more
general classes of real vector spaces.

We will consider real vector spaces equipped with
an inner product or a norm, possibly required to be complete (i.e., to
be Hilbert spaces or Banach spaces) and under various restrictions on the
dimension and often with multiplication disallowed.  So, for example, we
will find that the theories $\IP{\infty}$ and $\HS{\infty}$ of
infinite-dimensional inner product spaces and Hilbert spaces
respectively are both decidable, and in fact by the same decision
procedure, so that the two theories coincide. By contrast, we will see
that the analogous theories $\NS{\infty}$ and $\BS{\infty}$ of
infinite-dimensional normed spaces and Banach spaces differ, and both
are undecidable, as is the purely additive fragment $\BSA{d}$ of the
theory of $d$-dimensional Banach spaces for $d \ge 2$.

In fact, all the theories of inner product spaces we consider are
decidable, while for normed spaces, only the most trivial example,
namely the theory of a 1-dimensional space, is decidable. The
undecidable normed space theories are not recursively axiomatizable or
even arithmetical, as we will see by constructing primitive recursive
reductions of the set of truths of second-order arithmetic to these
theories. In fact, if we restrict to normed spaces of finite dimension,
the normed space theories have the same degree of unsolvability
(many-one degree) as second-order arithmetic, while for arbitrary
dimensions, the normed space and Banach space theories are many-one
equivalent to the set of true $\Pi^2_1$ sentences in third-order
arithmetic.

Normed spaces and inner product spaces are vector spaces with a metric that
relates nicely to the algebraic structure. We therefore consider metric
spaces as a source of motivating examples and for their own interest.  The
theory $\MS$ of metric spaces is known to be undecidable
\citep{bondi:dec-ms,kurz-metric}. We give an alternative proof which
shows that the theory is not arithmetical.

We obtain positive decidability results for normed spaces by
restricting the use of quantifiers: rather trivially, the set of valid
purely existential sentences is decidable, but much more interestingly,
so is the set of valid purely universal sentences.  The decision
procedure for the purely universal case is via a computational process
which (at least in principle) produces a concrete counter-example in the
shape of an explicit norm on $\real^n$ for some $n$ which fails to
satisfy a given invalid sentence.  This algorithm has been implemented
in the special case where multiplication is not allowed.  For metric
spaces, we do even better: the set of valid $\fAE$ sentences is
decidable, as we see using an analogue of the Bernays-Sch\"{o}nfinkel
decision procedure for valid $\fAE$ sentences in a first-order language
with no function symbols.  However, by reducing satisfiability for
quantifier-free formulas of arithmetic to the dual satisfiability
problem, we will find that validity for the $\fEA$ fragment is
undecidable for both metric spaces and normed spaces, as is validity for
the $\fAE$ fragment for normed spaces.
Finding other useful decidable fragments is an interesting challenge.

The structure of the sequel is as follows:

Section~\ref{sec:initobs} introduces notation and terminology and
then gives some preliminary observations and results.
We assume that the reader is acquainted with the concept of a many-sorted
first-order language as described, for example, in the book by
\citet{manzano-efol}.
However, as an {\it aide-memoire} to make the material more
easily accessible to readers without a professional
background in pure mathematics, we review many of the ideas
from vector algebra and affine geometry that we will use.
We then make
some initial observations on the possibilities for decision procedures in the
theories of interest. This leads on to a number of examples showing the
expressive strength of the language of normed spaces compared with the language
of inner product spaces.  For example, while we will later prove that a
first-order property of inner product spaces that holds in all finite
dimensional spaces holds in any inner product space, there are very simple
first-order properties of normed spaces that only have infinite-dimensional models.  The
section concludes with a proof that there are first-order properties that hold
in all Banach spaces but not in all normed spaces.

Section~\ref{sec:prelim} introduces our basic method for proving
undecidability in a language equipped with a sort whose intended interpretation
is the real numbers. The method is to exhibit a structure $\cal M$ for the
language and a formula $\nu(x)$ with the indicated free variable of the real sort
which holds in the structure iff $x$ is interpreted as a natural number.
If such a structure $\cal M$ exists, the method provides a reduction of the
set of truths of second-order arithmetic to the set of sentences that hold
in any class of structures containing $\cal M$. Thus the method shows
that a theory for which such a structure $\cal M$ exists is not even
arithmetical.

Section~\ref{sec:metric-spaces} applies the method of the previous section to
the case of metric spaces, which gives a new proof that the first-order theory
of a metric space is undecidable.  Here we also give a decision procedure for
(a superset of) the $\fAE$ fragment of the theory of metric spaces and show
that this is the best possible result of its type by reducing the
satisfiability problem for Diophantine equations to $\fAE$ satisfiability for
metric spaces.

Section~\ref{sec:undec} gives the main undecidability results for normed spaces
and Banach spaces: it turns out that in every dimension $d \ge 2$ we can apply
the methods of Section~\ref{sec:prelim} and prove undecidability of the
corresponding theories of normed spaces and Banach spaces, even for the purely
additive fragments where multiplication is disallowed. This section concludes
with a more detailed investigation into the degrees of unsolvability of these
theories: the theories for spaces of finite dimension $d \ge 2$ turn out to
have the same many-one degree as the set of truths in second-order arithmetic,
while if we allow infinite-dimensional spaces, the theories have the
same many-one degree as the set of true $\Pi^2_1$ sentences in third-order
arithmetic.

In Section~\ref{sec:qelim} we turn to inner product spaces and find that they
are quite tractable: the key result implies that a sentence holds in every
space of dimension $d \ge k$ iff it holds in $\real^k$ where $k$ is the
number of distinct vector variables in the sentence.  From this we find that
the theories of inner product spaces and Hilbert spaces with various
dimensional constraints can all be decided via a simple reduction to the first
order theory of the real numbers.

Section~\ref{sec:decidable-fragments} complements our investigation with some
results on decidable fragments of the normed space theory analogous to the
decidability results for metric spaces in Section~\ref{sec:metric-spaces}.  The
purely existential fragment admits a very simple reduction to the first-order
language of the real numbers.  The purely universal fragment is also decidable via a
more sophisticated method.

Again these results are the best possible of their type: in
Section~\ref{sec:ea-ae-fragments} we give reductions of satisfiability for
Diophantine equations to both the $\fAE$ and the $\fEA$ satisfiability problems
for normed spaces; in fact both these reductions are subsumed by our final
result which gives the undecidability of the set of $\fAIA$ sentences valid
for normed spaces.

Some of the results presented here have been foreshadowed by several authors
and some have been strengthened since the present paper was first written.
We conclude in Section~\ref{sec:concluding-remarks} with a brief
survey of related work.

The genesis of this paper was a question about decision
procedures for vector spaces asked several years ago of
Solovay by Harrison, and quickly answered with the first
proofs of decidability and quantifier elimination for inner
product spaces. Some time later, Harrison became interested
in corresponding questions for the theory of normed spaces
and implemented a decision procedure for the universal
additive theory. Arthan conjectured, however, that the full
theory of normed spaces is undecidable. On hearing this
conjecture, Solovay rapidly proved it and precisely
characterized the theory as many-one equivalent to the
fragment of third-order arithmetic discussed below. Arthan
refined these results to cover finite-dimensional spaces,
purely additive theories and formulas with limited
quantifier alternations, while Harrison extended the
decidability to the full universal theory and has done
further practical work on implementations. All hands have
contributed to the numerous improvements leading to the
present account.

\Section{The languages and their interpretation}\label{sec:initobs}

We will study sublanguages of a 2-sorted first-order language $\LL$.
$\LL$ itself
provides a full repertoire of first-order features for work in a Hilbert space. It
includes the operations of a vector space equipped with an inner product
together with the induced norm and metric. After introducing $\LL$,
we define sublanguages for other kinds of structure: $\LV$ for vector spaces,
$\LM$ for metric spaces, $\LN$ for normed spaces and
$\LI$ for inner product spaces.

\Subsection{Sorts}

The two sorts in $\LL$ are as follows:

\begin{enumerate}

\item $\ScaSort$ --- scalars

\item $\VecSort$ --- vectors or points

\end{enumerate}

The variables and constants in our many-sorted languages all carry a
label indicating their sort.
In $\LL$, we
adopt the familiar convention of bold font ($\Vx$, $\Vy$, $\VO$, etc.) for
vectors or points and regular font ($x$, $y$, $0$, etc.) for scalars.
If we need to write sort labels explicitly,
we will use superscripts, e.g., $x^{\cal N}$ will be a
variable of the natural number sort in second-order arithmetic.

\Subsection{Language}

We describe here the constant, function and predicate symbols
of $\LL$ and then define important sublanguages in later sections.
Following mathematical custom,
we overload many of the arithmetic operations like `+' for both scalars and
vectors, but this should not cause confusion given that our notation
distinguishes vector variables and constants from their scalar counterparts.
$\LL$ has the following constants and function symbols:

\begin{enumerate}

\item Scalar constants $n$ for all rational numbers $n$.

\item Addition ($x + y$), negation ($-x$) and multiplication ($x y$) of
scalars.

\item The zero vector or origin, $\VO$.

\item Addition ($\Vx + \Vy$) and negation ($-\Vx$) of vectors.

\item Multiplication of a scalar and a vector, with type $\ScaSort \times
\VecSort \to \VecSort$. We write the product of a scalar $c$ and vector $\Vx$
as $c \Vx$.

\item The inner (dot) product of vectors, with type $\VecSort \times \VecSort
\to \ScaSort$. We write the inner product of vectors $\Vx$ and $\Vy$ as
$\inner{\Vx,\Vy}$.

\item The norm operation on vectors, with type $\VecSort \to \ScaSort$. We
write $\|\Vx\|$ for the norm of a vector $\Vx$.

\item The distance function for metric spaces, with type $\VecSort \times
\VecSort \to \ScaSort$. We write $d(\Vx,\Vy)$ for the distance between $\Vx$
and $\Vy$.

\end{enumerate}

We will also use the usual shorthands such as $x-y$ (for $x + (-y)$),
$x^2$ (for $xx$) and $\Vv/2$ (for $(1/2)\Vv$).
Nothing of substance would change if we also added a
multiplicative inverse operation. However, it can always be eliminated if
necessary. In any case, if a multiplicative inverse is to be included, adapting
the results of this paper is much more straightforward and efficient if
$0^{-1}$ has a specific known value.

The predicate symbols are:

\begin{enumerate}

\item Equality $\Vv = \Vw$ of vectors.

\item All the usual equality and inequality comparisons for scalars:
$x = y$, $x < y$, $x \leq y$, $x > y$, $x \geq y$.

\end{enumerate}

We use $\Abs{x}$ as a shorthand for the absolute value of $x$ and
$\Max\{x, y\}$ as a shorthand for the maximum of $x$ and $y$:
$\phi(\Max\{x, y\})$ stands for $x \ge y \And \phi(x) \Or x < y \And \phi(y)$ and $\phi(\Abs{x})$ stands for $\phi(\Max\{x, -x\})$.
Recursively, we write
$\Max\{x_1, x_2, \ldots, x_k\}$ for
$\Max\{x_1, \Max\{x_2, \ldots, x_k\}\})$.

\Subsection{Interpretation}

All the languages considered here include a symbol for equality for every sort
and this is to be interpreted as actual equality in any structure.  For
sublanguages of $\LL$, unless otherwise stated, we require the sort $\ScaSort$
and the symbols for the field operations and the ordering to be interpreted as
the ordered field of real numbers. Thus all the first-order properties of
$\real$ form part of the theory while we may make free use of higher-order
properties such as completeness when we reason about it.

\Subsubsection{Vector spaces}

The language $\LV$ of vector spaces has the scalar constant, function and
predicate symbols together with the constant $\VO$ and addition, negation and scalar multiplication
for vectors.  A vector space is a structure for this language satisfying the
vector space axioms listed below.  These state that the vectors form an Abelian
group on which the field of scalars acts a ring of homomorphisms:

\begin{itemize}

\item $\all{\Vu\;\Vv\;\Vw} \Vu + (\Vv + \Vw) = (\Vu + \Vv) + \Vw$

\item $\all{\Vv\;\Vw} \Vv + \Vw = \Vw + \Vv$

\item $\all{\Vv} \VO + \Vv = \Vv$

\item $\all{\Vv} -\Vv + \Vv = \VO$

\item $\all{a\;\Vv\;\Vw} a (\Vv + \Vw) = a \Vv + a \Vw$

\item $\all{a\;b\;\Vv} (a + b) \Vv = a \Vv + b \Vv$

\item $\all{\Vv} 1 \Vv = \Vv$

\item $\all{a\;b\;\Vv} (a b) \Vv = a (b \Vv)$.

\end{itemize}


The simplest example of a vector space comprises the single element $\VO$ and
is called $0$. One can define a vector space structure component-wise on
the set $\real^n$ of $n$-tuples of real numbers, and $0$ can be considered the
degenerate case $n = 0$. The space $\real^n$ contains the {\em standard basis}
$\{\Ve_1,\ldots,\Ve_n\}$:
\begin{eqnarray*}
\Ve_1 & = & (1,0,0,\cdots,0)    \\
\Ve_2 & = & (0,1,0,\cdots,0)    \\
    & \cdots &                  \\
\Ve_n & = & (0,0,0,\cdots,1).
\end{eqnarray*}

A fundamental result is that every vector space $V$ has a {\em basis}, i.e., a
set of vectors $B$ such that {\em(i)} any vector $\Vx \in
V$
can be represented as a linear combination
$\Vx = x_1\Vb_1 + \ldots + x_m\Vb_m$ for some $m \in \nat$, $x_i \in \real$ and $\Vb_i \in B$
($B$ {\em spans} $V$) and {\em(ii)} this representation is unique ($B$ is {\em linearly independent}).  The standard basis
$\{\Ve_1,\ldots,\Ve_n\}$ is indeed a basis for $\real^{n}$. Any two bases of a
vector space $V$ have the same cardinality called the {\em dimension} of $V$
and we will write $\Dim{V} = n$ if $V$ has a finite basis with $n$ elements,
otherwise we write $\Dim{V} = \infty$.

A {\em subspace} of a vector space is a substructure that also
interprets $\ScaSort$ as the field of all real numbers.  A subspace is
automatically a vector space, since the vector space axioms are purely
universal.  An analogous definition applies to all our notions of a
``space'', a subspace being given by any subset of the vectors or points
that is closed under all relevant operations.  Two subspaces $U$ and $W$
of a vector space $V$ are said to be {\em complementary} if every $\Vv
\in V$ can be written uniquely as $\Vv = \Vu + \Vw$ with $\Vu \in U$ and
$\Vw \in W$, in which case the dimension of $W$ depends only on $U$ and $V$ and
is said to be the {\em codimension} of $U$ in $V$.  Any subspace of a
vector space has at least one complementary subspace.

If $A$ is any set, the set $A \to \real$ of all real-valued functions
on $A$ becomes a vector space if one defines $(\Vf+\Vg)(a) = \Vf(a) + \Vg(a)$ and $(x\Vf)(a) = x\Vf(a)$.
Taking $A = \nat$, the elements of $A \to \real$ are sequences of real numbers
and we define $\real^*$ to be the subspace
comprising sequences $(\Vx_0, \Vx_1, \ldots)$ whose support $\{ n \mid \Vx_n
\not= 0 \}$ is finite. This space is infinite-dimensional since the unit
vectors $(0,\ldots,0,1,0,\ldots)$ are linearly independent.
Identifying the $n$-tuple $(x_1, \ldots, x_n)$ with the sequence $(x_1, \ldots, x_n, 0,  \ldots)$,
$\real^*$ can be viewed as the union of the spaces $\real^n$ for $n \in \nat$.

Many useful geometric notions can be defined just in terms of the vector
space operations.  If $\Vv$ and
$\Vw$ are distinct vectors, the {\em affine line} passing through them comprises
the set of points that can be written as linear combinations $a\Vv + b\Vw$
where $a + b = 1$.  The points of this form with $a, b \ge 0$ comprise the {\em closed
line segment} $[\Vv, \Vw]$, while those with $a, b > 0$ form the {\em open line
segment} $(\Vv, \Vw)$.  We say the line segment $[\Vv, \Vw]$
is {\em parallel} to a subspace $W$ iff, for some $\Vu$, $[\Vu +
\Vv, \Vu + \Vw]$ is contained in $W$.

A set of vectors $A$ is said to be {\em convex} if it contains the line segment
connecting any two of its points. Following the convention that quantifiers
have lower precedence than propositional operators (so the scope of a quantifier
extends as far to the right as possible), we express this formally as follows:
$$ \all{\Vv\;\Vw} \Vv \in A \And \Vw \in A
   \Imp \all{a\;b} 0 \leq a \And 0 \leq b \And a + b = 1
                   \Imp a \Vv + b \Vw \in A.
$$
If $A$ is any set of vectors, its {\em convex hull}, $\Conv{A}$, is the
smallest convex set containing $A$ (this  is well-defined because the
intersection of any family of convex sets is convex). $\Conv{A}$ comprises
all the {\em convex combinations} of elements of $A$,
i.e., all finite sums $a_1\Vv_1 + \ldots + a_m\Vv_m$ where $a_i \ge 0$,
$\Vv_i \in A$, $i = 1, \ldots m$, $a_1 + \ldots + a_m = 1$ and $m \ge
1$.  If $A$ is a finite set with $n$ elements and if each element of $\Conv{A}$ 
has a unique representation as a convex combination of elements of $A$, then
$\Conv{A}$ is said to be an {\em $(n-1)$-simplex} and the points of $A$ are its
{\em vertices}.  So, for example, a 1-simplex is a closed line segment while a
2-simplex is a triangle.

\Subsubsection{Metric Spaces}

The language $\LM$ of metric spaces has all the scalar constant, function and
predicate symbols together with the metric $d(\_, \_)$ as the only function
symbol involving the point type $\VecSort$. A metric space is a structure for
this language satisfying the metric space axioms listed below: positive
definiteness, symmetry and the triangle inequality.

\begin{itemize}

\item $\all{\Vx\;\Vy} d(\Vx,\Vy) \geq 0 \And (d(\Vx,\Vy) = 0 \Iff \Vx = \Vy)$

\item $\all{\Vx\;\Vy} d(\Vx,\Vy) = d(\Vy,\Vx)$

\item $\all{\Vx\;\Vy\;\Vz} d(\Vx, \Vz) \leq d(\Vx,\Vy) + d(\Vy,\Vz)$.

\end{itemize}

%

A metric space is said to be {\em complete} if every Cauchy sequence converges.
Unsurprisingly, it turns out to be impossible to capture this notion by
first-order axioms in our language (see Theorem~\ref{thm:incompleteness}), but if we
allow quantification over infinite sequences of points, we can express it
as follows, where $\Vx$ ranges over such sequences:
$$ \all{\Vx}
\BA (\all{\epsilon} \epsilon > 0
\Imp \ex{N} \all{m\;n} m \geq N \And n \geq N
\Imp d(\Vx(n), \Vx(m))< \epsilon) \\
\Imp \ex{\Vl} \all{\epsilon} \epsilon > 0
\Imp \ex{N} \all{n} n \geq N \Imp d(\Vx(n),\Vl) < \epsilon.
\EA
$$

\Subsubsection{Normed spaces}\label{subsubsec:normed-spaces}

The language $\LN$ of normed spaces includes all the symbols of the language
$\LV$ of vector spaces together with a norm $\Norm{\_}$. A metric $d$ may also
be used  as a notational convenience (see below).  A normed space is a
structure for this language that satisfies the axioms for a vector space
together with the axioms for norms listed below: positive definiteness, scaling
and the triangle inequality:

\begin{itemize}

\item $\all{\Vv} \Norm{\Vv} \ge 0 \And (\|\Vv\| = 0 \Iff \Vv = \VO)$

\item $\all{a\;\Vv} \|a \Vv\| = |a| \|\Vv\|$

\item $\all{\Vv\;\Vw} \|\Vv + \Vw\| \leq \|\Vv\| + \|\Vw\|$.

\end{itemize}

As a function from the space to the real numbers the norm is continuous with
respect to a topology defined by the {\em induced metric:} $d(\Vv, \Vw) = \|\Vv -
\Vw\|$.  This will be a very useful fact in our later arguments. The continuity
of the norm at the point $\Vv$ can be expressed in our first-order language as
follows:
$$ \all{\epsilon}
   \epsilon > 0
   \Imp \ex{\delta}
       \delta > 0
       \And \all{\Vw} \| \Vw \| < \delta \Imp
             | (\|\Vv + \Vw\| - \|\Vv\|) | < \epsilon.
$$


A {\em Banach space} is a normed space that is also metrically complete, i.e.,
with respect to the induced metric, every Cauchy sequence converges.  As with
metric spaces, we shall prove later that it is impossible to capture this by
first-order axioms in our language (see Theorem~\ref{thm:incompleteness}).

The usual euclidean norm on $\real^n$ is defined by $\|\Vx\| =
\sqrt{\sum_{i=1}^n \Vx_i^2}$, but there are plenty of other possibilities
satisfying the axioms, such as the 1-norm (``Manhattan distance'') $\|\Vx\| =
\sum_{i=1}^n |\Vx_i|$ and the $\infty$-norm $\|\Vx\| = \Max \{ |\Vx_i| \mid 1
\leq i \leq n \}$. Similar norms can be defined on $\real^*$ by summing or
maximizing over the support rather than from $1$ to $n$. Other examples from
functional analysis include the norm $\|f\| = \Sup \{ |f(x)| \mid x \in [a,b]
\}$ on the Banach space of continuous functions $f \ty [a,b] \to \real$.

All norms on a finite-dimensional vector space, $V$, can be shown to be
equivalent in the sense that if $\Norm{\_}_1$ and $\Norm{\_}_2$ are norms on
$V$, then there are positive real numbers $s$ and $t$ such that for any $\Vv
\in V$, $\frac{1}{s}\Norm{\Vv}_1 \le \Norm{\Vv}_2 \le t\Norm{\Vv}_1$.  Although
this implies that many properties of interest, in particular topological ones,
are independent of the norm, we shall see that there are very great differences
in the general first-order properties satisfied by different norms on the same
finite-dimensional vector space.

Each norm defines a corresponding {\em unit circle}
$S = \{\Vx \mid \|\Vx\| = 1\}$ and a {\em unit disc} $D = \{\Vx \mid \|\Vx\| \leq
1\}$.
In spaces of higher dimension we also sometimes refer to $S$
and $D$ as the {\em unit sphere} and {\em unit ball} respectively.
For the usual euclidean norm $\|(x_1,x_2)\| = \sqrt{x_1^2 + x_2^2}$
on $\real^2$, $S$ and $D$ are indeed a circle and a disc respectively.
However many other shapes are possible,
e.g. for the $\infty$-norm on $\real^2$, $D$ is a square.
However, $D$ is always a convex set: if $\Vx$ and
$\Vy$ are in $D$ then $\|a \Vx + b \Vy\| \leq \|a \Vx\| + \|b \Vy\| =
|a| \|\Vx\| + |b| \|\Vy\| \leq |a| + |b|$, and if $a, b \geq 0$ with $a + b =
1$ we have $|a| + |b| = 1$.
$D$ is also always symmetric about the origin in the sense that $\Vv \in D$ iff $-\Vv \in D$. As $D$ is convex and $S \subseteq D$, any convex
combination of unit vectors (i.e., members of $S$) has norm at most 1.

Conversely, it is often convenient to define a norm by nominating a suitable
set as its unit disc $D$ and defining the norm by taking $\|\Vx\|$ to be the
smallest non-negative real number $\lambda$ such that for some
$\Vd \in D$ one has $\lambda \Vd = \Vx$. Provided the set $D$ is convex and
meets every line through the origin in a closed line segment $[-\Vv, \Vv]$
where $\Vv \not= \VO$, this is well-defined and satisfies the norm properties. For example, if $\Vx
= \|\Vx\|\Vd$ and $\Vy = \|\Vy\| \Ve$ for $\Vd, \Ve \in D$ then
$$ \Vx + \Vy = (\|\Vx\|\ + \|\Vy\|) ({\|\Vx\|\ \over \|\Vx\|\ + \|\Vy\|} \Vd +
                                {\|\Vy\| \over \|\Vx\|\ + \|\Vy\|} \Ve)
$$

\noindent and since by convexity $({\|\Vx\|\ \over \|\Vx\|\ + \|\Vy\|} \Vd +
{\|\Vy\| \over \|\Vx\|\ + \|\Vy\|} \Ve) \in D$, we have $\|\Vx + \Vy\| \leq
\|\Vx\|\ + \|\Vy\|$, i.e., the triangle inequality holds.

Under the euclidean norm, the unit circle meets any affine line in at most two
points.  However, there are many interesting norms for which this is not the
case: in the $\infty$-norm in $\real^2$, for example, the unit circle comprises
the union of four line segments. In working with such norms, it is useful to
note that if $L$ is an affine line, then $L \cap D$, the set of points on $L$
of norm at most 1, is a bounded convex subset of $L$ whose endpoints are
contained in $S$. So if $L \cap D$ is non-empty, either $L
\cap D = L \cap S = \{\Va\}$ for some $\Va$ with $\Norm{\Va} = 1$, or $L \cap D
= [\Va, \Vb]$ for some distinct $\Va$, $\Vb$ with $\Norm{\Va} = \Norm{\Vb} =
1$.  In the latter case, either $L \cap D \subseteq S$, i.e., $\Norm{\Vx} = 1$
for every $\Vx \in [\Va, \Vb]$, or $L \cap S = \{\Va, \Vb\}$ and $\Norm{\Vx} <
1$ for every $\Vu \in (\Va, \Vb)$.  In particular, the condition $\Norm{\Vv} =
\Norm{\Vw} = \Norm{(\Vv + \Vw)/2} = 1$ implies that $[\Vv, \Vw] \subseteq S$.

More generally, let $\Vv_1, \ldots, \Vv_n$ be the vertices of an
$(n-1)$-simplex, $\Delta$, and let $\Vu = a_1\Vv_1 + \ldots + a_n\Vv_n$ be any
proper convex combination of the $\Vv_i$, i.e., $a_i > 0$, $i = 1 \ldots n$,
and $a_1 + \ldots + a_n = 1$.  Then $\{\Vu, \Vv_1, \ldots, \Vv_n\} \subseteq S$
implies that  $\Delta \subseteq S$, i.e., $\Norm{\Vu} = \Norm{\Vv_1} = \ldots =
\Norm{\Vv_n} = 1$ implies that $\Norm{\Vx} = 1$ for every $\Vx \in \Delta$. To
see this, first note that $\Delta \subseteq D$ because a convex combination of
unit vectors has norm at most 1.  For $1 \le i \le n$, let $L_i$ be the affine
line passing through $\Vu$ and $\Vv_i$ and let $\Delta_i$ be the
$(n-2)$-simplex whose vertices are $\Vv_1, \ldots, \Vv_{i-1}, \Vv_{i+1}, \ldots
\Vv_n$.  $L_i$ meets $\Delta_i$ at a point $\Vu_i$, say, that must be a proper
convex combination of the vertices of $\Delta_i$.  Since $\Vv_i$ and $\Vu$ are
unit vectors and $\Vu$ lies on the open line segment $(\Vu_i, \Vv_i)$, we have
$\Norm{\Vu_i} \ge 1$.  Then as $\Vu_i \in \Delta \subseteq D$, $\Norm{\Vu_i} =
1$ and so, by induction, $\Delta_i \subseteq S$.  Let $\Vx$ be any point of
$\Delta$. As $\Delta \subseteq D$, to show that $\Norm{\Vx} = 1$, it suffices
to show that $\Norm{\Vx} \ge 1$.  If $\Vx = \Vu$, then $\Norm{\Vx} = 1$ by
assumption.  So assume $\Vx \not= \Vu$. For some $i$, the half-line starting at
$\Vx$ and passing through $\Vu$ meets the $(n-2)$-simplex $\Delta_i$ at a point
$\Vw$. As $\Vw$ and $\Vu$ are both unit vectors and $\Vu$ lies on the open line
segment $(\Vx, \Vw)$, we find $\Norm{\Vx} \ge 1$ completing the proof.
As a special case we have that if the unit vectors $\Vv_1, \ldots, \Vv_n$ are
the vertices of an $(n-1)$-simplex $\Delta$, then $\Delta \subseteq S$ iff the
barycentre $\frac{1}{n}(\Vv_1 + \ldots + \Vv_n)$ is a unit vector.
Note that this gives a considerable economy from a logical point of
view: we can assert $\Delta \subseteq S$ without using any
quantifiers or scalar variables.

\Subsubsection{Inner product spaces}

The language $\LI$ of inner product spaces includes all the
symbols of the language $\LV$ of vector spaces together with an inner product
$\inner{\_, \_}$. A norm may also be used as a notational convenience (see
below).  An inner product space satisfies the axioms for a vector space
together with the axioms asserting that inner product is a
positive definite symmetric bilinear form, which means:

\begin{itemize}

\item $\all{\Vv\;\Vw} \inner{\Vv,\Vw} = \inner{\Vw,\Vv}$

\item $\all{\Vu\;\Vv\;\Vw} \inner{\Vu + \Vv,\Vw} = \inner{\Vu,\Vw} +
\inner{\Vv,\Vw}$

\item $\all{a\;\Vv\;\Vw} \inner{a \Vv,\Vw} = a \inner{\Vv,\Vw}$

\item $\all{\Vv} \inner{\Vv,\Vv} \geq 0 \And (\inner{\Vv,\Vv} = 0 \Iff \Vv = \VO)$.

\end{itemize}

For example, $n$-dimensional euclidean space is $\real^n$ equipped with the
inner product $\inner{\Vx,\Vy} = \sum_{i=1}^n \Vx_i \Vy_i$. Note that
$\inner{\Vx,\Vx} = \|\Vx\|^2$ for the euclidean norm, and in general given
any inner product we define the {\em induced norm} by $\|\Vx\| =
\sqrt{\inner{\Vx,\Vx}}$.

A {\em Hilbert space} is an inner product space that is also complete for the
induced norm. Any finite-dimensional inner product space is a Hilbert
space.  The vector space of sequences $\Vx \ty \nat \to \real$ such that the
sum $\sum_{i=0}^{\infty}\Vx_i^2$ is convergent is an infinite-dimensional
Hilbert space under an inner product defined by $\inner{\Vx,\Vy} =
\sum_{i=0}^{\infty} \Vx_i \Vy_i$.  This Hilbert space, called $l_2$, is one of
many Hilbert spaces that occur naturally in functional analysis.
The vector space $\real^{*}$ of finitely-supported sequences viewed as a
subspace of $l_2$ gives an example of an incomplete inner product space.

If $\Vu$ and $\Vv$ are elements of an inner product space $V$, we say $\Vv$ is
orthogonal to $\Vu$ if $\inner{\Vu, \Vv} = 0$. If $\Vu$ is non-zero then the
set $W$ of all vectors orthogonal to $\Vu$ forms a subspace $W$ of $V$ called
the {\em orthogonal complement} of $\Vu$. Every element $\Vv$ of $V$ can be
written uniquely in the form $\Vv = a\Vu + \Vw$ where $\Vw$ is a member of $W$.

\Subsection{Additive sublanguages}

The so-called linear fragment of real arithmetic admits a very simple
quantifier elimination procedure \citep{hodes,ferrante-rackoff} and enjoys many other
pleasant properties. Here ``linear'' means that the multivariate
polynomials that are the terms of the language are restricted to have total
degree at most one.  To define an analogous notion for $\LL$ (or any of
its sublanguages or extensions thereof by the addition of extra vector
constants), we say a term or formula is {\em additive} if the left operand of
every subterm of the form $xy$, $x\Vv$ or $\inner{\Vv, \Vw}$ is a constant.  In
$\LL$ itself, which has only rational scalar constants and the vector constant
$\VO$, an additive formula is equivalent to one in which multiplication and
inner product do not occur.
E.g., one can write $\Vq + \Vq = \Vp + \Vr$ rather than
$\Vq = (\Vp + \Vr)/2$ to indicate that $\Vq$ is the midpoint of the line
between $\Vp$ and $\Vr$.
We write $\LA$ for the additive sublanguage of $\LL$.

Unless otherwise stated, in a structure for one of the additive sublanguages,
we will require the sort $\ScaSort$, the symbols for the additive group
operations and the ordering to be interpreted as the ordered additive group of
real numbers with the rational number constants interpreted accordingly.

\Subsection{Initial observations on decidability}

The principal results of this paper are connected with decidability or
undecidability for the various 2-sorted languages introduced above. We now make
some initial observations about the possibilities for decision procedures,
e.g., via quantifier elimination, and about the interrelations among the
decision problems for the languages.

\Subsubsection{Reductions among decision problems}

Recall that every vector space $V$ has a basis, i.e., a subset $B$ such that
any vector $\Vx \in V$ can be written uniquely as a sum $\sum_{\Vb \in B}
x_{\Vb} \Vb$ (where all but finitely many $x_{\Vb}$ are zero).  Given any basis
we can regard the scalar coefficient $x_{\Vb}$ as the ``$\Vb^{\it th}$
coordinate'' of $\Vx$, and for $\Vy = \sum_{\Vb \in B} y_{\Vb}\Vb$,
we can define $\inner{\Vx,\Vy} = \sum_{\Vb \in B} x_{\Vb} y_{\Vb}$
and show that this satisfies the inner product properties.  Thus every vector
space can be made into an inner product space; in logical parlance, this
implies that the theory of inner product spaces is a conservative extension of
the theory of vector spaces:

\begin{quote}
A formula using neither the inner product nor norm operation holds in all
vector spaces [optionally with constraints on the dimension] iff it holds in
all inner product spaces [with corresponding constraints].
\end{quote}

As noted already, in any inner product space we can define $\|\Vx\| =
\sqrt{\inner{\Vx,\Vx}}$ and this satisfies the norm properties, so any model of
the inner product space axioms immediately gives a
model of the normed space axioms. The converse is
not true, i.e. not every normed space is an inner product space.  (See the
remarks at the end of Section~\ref{sec:qelim} for a more quantitative statement
on this topic). However if a normed space {\em is} derived from an inner
product as above, the inner product can be recovered from the norm, e.g. by
$\inner{\Vx,\Vy} = \frac{1}{2}(\|\Vx + \Vy\|^2 - \|\Vx\|^2 - \|\Vy\|^2)$.
It is a classic result of \citet{Jordan-von-Neumann} that a norm is induced
by an inner product iff it satisfies the parallelogram identity:
$$
\all{\Vx\;\Vy} \Norm{\Vx + \Vy}^2 + \Norm{\Vx - \Vy}^2 = 2(\Norm{\Vx}^2 + \Norm{\Vy}^2)
$$
\noindent See Section~\ref{subsubsec:further-expressiveness-results}
below for more about characterizations of inner product spaces.

Let $\iota$ be a sentence in the language $\LI$ of inner product spaces
asserting that `$\inner{\_,\_}$' satisfies the inner product axioms.
Given any formula $\phi$ in $\LI$, let $\phi^*$ be the corresponding formula
in the language $\LN$ of normed spaces where each term $\inner{\Va,\Vb}$
is replaced by $(\|\Va + \Vb\|^2 - \|\Va\|^2
- \|\Vb\|^2)/2$. If $M$ is an inner product space in which $\phi$ holds,
  then $\iota \And \phi$ holds in $M$. In that case $\iota^* \And \phi^*$ holds in the
normed space $N$ derived from $M$ by defining $\|\Vx\|_N =
\sqrt{\inner{\Vx,\Vx}_M}$. Conversely, if $\iota^* \And \phi^*$ holds in a
normed space $N$, then setting $\inner{\Vx,\Vy}_M = (\|\Vx + \Vy\|_N^2 -
\|\Vx\|_N^2 - \|\Vy\|_N^2) / 2$ makes $N$ into an inner product space,
$M$ say, in which $\phi$ holds.  Both these constructions preserve
dimensions and completeness and so restating in terms of validity, we have:

\begin{quote}
A sentence $\phi$ in the language of inner product spaces holds in all inner
product spaces [with or without constraints on the dimension and with or
without the requirement for completeness] iff the sentence $\iota^* \Imp \phi^*$ (as
defined above) holds in all normed spaces [with or without corresponding
limitations].
\end{quote}

This establishes that the decision problem for normed spaces is at least as
general as the decision problem for inner product spaces, which in turn is at
least as general as the decision problem for vector spaces. It will emerge
in what follows that the decision problem for normed spaces is in fact
dramatically harder than the other two. Intuitively, one might see this as
expressing the fact that one has freedom to describe very ``exotic'' norms,
whereas the freedom to define inner products is more constrained.

\Subsubsection{Possibility of quantifier elimination}\label{subsubsec:poss-qe}

It is not hard to see that we cannot have quantifier elimination in the basic
language we are considering, for any of the vector space theories. For if so,
any closed formula would be equivalent to a ground formula. Now the
vector-valued subterms in a ground formula are formed from $\VO$
using addition and scalar multiplication and so evaluate to $\VO$
in any model. Thus the truth of a ground formula is independent of the
space in which it is interpreted.  So quantifier elimination would imply
that all models are elementarily equivalent. This is certainly
not the case however: we can write down formulas expressing non-trivial
properties of the dimension and/or the norm. For example, the dimension is
finite and $\leq n$ iff there is a spanning set of at most $n$ vectors:
$$\ex{\Vv_1 \ldots \Vv_n} \all{\Vw}
        \ex{a_1 \ldots a_n} a_1 \Vv_1 + \cdots + a_n \Vv_n = \Vw.
$$

We will see in Section~\ref{sec:qelim} that, if these sentences
$\pD_{{\le}n}$ are treated as atomic predicates, there {\em is} a full
quantifier elimination algorithm for vector spaces and for inner product
spaces. This allows us to decide validity in all vector spaces or all those
with a specific restriction on the dimension. Moreover, it implies the
existence, for any formula $\phi$ in this theory, of a bound $k$ such that $\phi$
holds in all vector (or inner product) spaces iff it holds in all those of
dimension at most $k$. In other words, if a formula $\phi$ in the language of inner
product spaces is satisfiable, it is satisfiable in an inner product space
with a specific finite upper bound on the dimension.

If we turn to normed spaces, however, the situation changes dramatically. We
will see in Section~\ref{sec:undec} that the theory is undecidable, so
no algorithmically useful quantifier elimination in an expanded
language exists. We will show below that there are satisfiable formulas that
are satisfiable only in infinite-dimensional normed spaces. Moreover,
quantifier elimination in the unexpanded language must even fail for purely
additive formulas (no scalar multiplication or inner products, and scalar-vector
multiplication only for integer constants), since we can for example express
the fact that the dimension is $\leq 1$ by:
$$ \ex{\Vx} \all{\Vy} \|\Vy\| = 1 \Imp \Vy = \Vx \Or \Vy = -\Vx $$

\noindent and distinguish the 1-norm and 2-norm by:
$$ \all{\Vx\;\Vy}
   \|\Vx\| = \|\Vy\| \And \|\Vx + \Vy\| = \|\Vx\| + \|\Vy\| \Imp \Vx = \Vy.
$$

(This holds for the euclidean norm in any number of dimensions, but fails in
$\real^2$ with the 1-norm $\|(x_1,x_2)\| = |x_1| + |x_2|$, as can be seen by
setting $\Vx = (1,0)$, $\Vy = (0,1)$.)

\Subsubsection{Further expressiveness results for normed spaces}
\label{subsubsec:further-expressiveness-results}

There are (purely additive) formulas in the language of normed spaces that are
satisfiable yet have only infinite-dimensional models. To see this,
define a 1-place predicate $\pE{\Vv}$ that holds iff $\Vv$ is a unit vector
that is not the midpoint of the line connecting two distinct vectors in the
unit disc, i.e., $\Vv$ is an extreme point of the unit disc:
$$ \pE{\Vv} \IsDef
     \|\Vv\| = 1 \And
     \all{\Vu\;\Vw} \|\Vu\| \leq 1 \And \|\Vw\| \leq 1 \And \Vv = (\Vu + \Vw) / 2
               \Imp \Vu = \Vw.
$$

Now consider the sentence $\pInf$ asserting that there exist non-zero vectors but
that the unit disc has no extreme points:
$$ \pInf \IsDef (\ex\Vv \Vv \not= \VO) \And (\all{\Vv} \Not \pE{\Vv}).
$$

\noindent In a finite-dimensional normed space, the Krein-Milman theorem
implies that the unit disc is the convex hull of its extreme
points, so $\pInf$ cannot hold in finite dimensions.  But when equipped
with the $\infty$-norm, the space
$\real^*$ considered above (sequences of real numbers with finite
support) has a unit disc with no extreme points:
given any unit vector $\Vv$, pick an $n$ so that $\Vv_n = 0$ and set
$\Vu_n = -1$, $\Vw_n = 1$ and $\Vu_i = \Vw_i = \Vv_i$ for $i \not= n$;
then $\Vv = (\Vu + \Vw) / 2$.  Hence $\pInf$ holds in $\real^*$, so $\pInf$ is
satisfiable but only has infinite-dimensional models.

\begin{figure}
\begin{center}
\includegraphics[angle=0,scale=0.9]{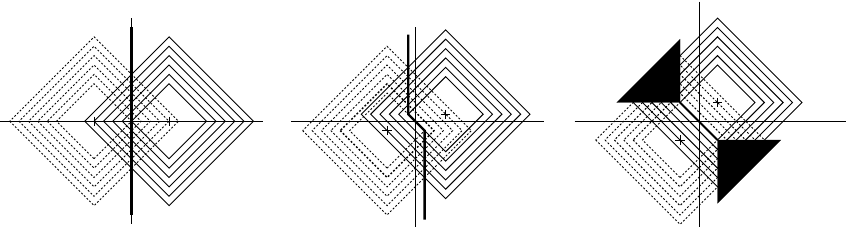}
\caption{Examples of $\{\Vx \ST \pO{\Vx, \Vy}\}$ in the 1-norm on $\real^2$}
\label{fig:pseudoorthog}
\end{center}
\end{figure}

It is also interesting to observe that using the norm, we can find
purely additive sentences that are satisfiable, but only in certain
models with a specific finite dimension. In fact without using
multiplication we can even characterize specific norms, e.g., the
1-norm and the euclidean norm on $\real^n$.

To give these characterizations, we will
write $\pO{\Vx,\Vy}$ for $\|\Vx + \Vy\| = \|\Vx - \Vy\|$,
i.e.  $\pO{\Vx, \Vy}$ holds iff $\Vx$ is equidistant from $\Vy$ and $-\Vy$.
Intuitively, this is intended as an approximation to the concept of
orthogonality in an inner product space. Indeed, for a norm derived in the
usual way from an inner product, this says exactly that $\inner{\Vx,\Vy} = 0$.
In a general normed space, $\pO{\Vx, \Vy}$ will not enjoy all the properties of
orthogonality, and, in particular, the ``orthogonal complement'', $C_{\Vy} \IsDef
\{\Vx \mid \pO{\Vx,\Vy}\}$, need not be a
subspace, as illustrated  for the 1-norm on $\real^2$
in Figure~\ref{fig:pseudoorthog}.

Assume $\Norm{\Vx + \Vy} < \Norm{\Vx - \Vy}$ for some $\Vx$ and $\Vy$ in
some normed space $V$, so that certainly $\Vx, \Vy \not = \VO$.  Let
$f(s) = \Norm{s\Vx + (2-s)\Vy} - \Norm{s(\Vx - \Vy)}$ so that $f(s) = 0$
iff $\pO{s\Vx + (1-s)\Vy, \Vy}$ holds. $f(s)$ is a continuous function
of $s$ with $f(0) = 2\Norm{\Vy} > 0$ and $f(1) = \Norm{\Vx + \Vy} -
\Norm{\Vx - \Vy} < 0$.  By the intermediate value theorem, $f(t)$ = 0
for some $t > 0$.  So $\Vx = a\Vy + b\Vz$, where $a = \frac{t-1}{t}$, $b
= \frac{1}{t}$ and $\Vz = t\Vx + (1 - t)\Vy$. As $f(t) = 0$, $\pO{\Vz,
\Vy}$ holds.  Similarly, if $\Norm{\Vx + \Vy} > \Norm{\Vx - \Vy}$, we
can also find $\Vz$ such that $\pO{\Vz, \Vy}$ holds and $\Vx = a\Vy + b\Vz$
for some $a$ and $b$.
Since $\Norm{\Vx + \Vy} = \Norm{\Vx - \Vy}$ implies $\pO{\Vx, \Vy}$,
any $\Vx \in V$ can be
written as a linear combination $\Vx = a\Vy + b\Vz$, where $\pO{\Vz,
\Vy}$, i.e., $\Vy$ and $C_{\Vy}$ span $V$.

Now assume that for some $\Vy \not= \VO$, the set $C_{\Vy}$
is a subspace. Then if $\Vz \in C_{\Vy}$, so also is $b\Vz$. Thus
any $\Vx \in V$ can be written as $a\Vy + \Vz$ where $\Vz \in C_{\Vy}$,
and, as $\Vy \not\in C_{\Vy}$, this representation is unique.
Thus $C_{\Vy}$ has codimension 1 and $\Vy$ spans a complementary
subspace.  If $\Vz \in C_{\Vy}$ and $a \not= 0$, then $\Vz/a \in
C_{\Vy}$ and we have
$
\Norm{a\Vy + \Vz} = |a|\cdot\Norm{\Vy + \Vz/a} = |a|\cdot\Norm{\Vy -
\Vz/a} = \Norm{{-a}\Vy + \Vz}.
$
Thus, if $C_{\Vy}$ is a subspace, there
is a (unique) linear isometry from $V$ to itself that fixes $C_{\Vy}$
and maps $\Vy$ to $-\Vy$.
For example, for $\Vy \not= \VO$ in
$\real^2$ under the 1-norm, $C_{\Vy}$ is a subspace iff $\Vy$ lies on
one of the coordinate axes, in which case $C_{\Vy}$ is the other
axis and reflection in it gives the linear isometry mapping $\Vy$ to $-\Vy$
(see Figure~\ref{fig:pseudoorthog}).

For any $n \in \nat$, there is a sentence $\phi_n$ of $\LNA$ which holds in a normed
space iff there are vectors $\Ve_1$, \ldots, $\Ve_n$ such that:

\begin{itemize}

\item $\Norm{\Ve_i} = 1$ for each $i$

\item $\pO{\Ve_i,\Ve_j}$ for each $i \not= j$

\item $\all{\Vv\;\Vw} \pO{\Vv,\Ve_i} \And \pO{\Vw,\Ve_i} \Imp \pO{\Vv + \Vw,\Ve_i}$
for each $i$

\item $\all{\Vv} \pO{\Vv,\Ve_i} \Imp \pO{{1 \over 2} \Vv,\Ve_i}$ for each $i$

\item $\all{\Vv} \pO{\Vv,\Ve_1} \And \cdots \And \pO{\Vv,\Ve_n} \Imp \Vv = \VO$.

\end{itemize}

I claim that in any model $V$ of $\phi_n$, the set
$W_i = C_{\Ve_i} = \{\Vx \mid \pO{\Vx,\Ve_i}\}$ is a
subspace, and hence, by the above remarks, a subspace of codimension 1.
To see that $W_i$ is indeed a subspace, note that, by induction, if $\pO{\Vv,\Ve_i}$ holds then so does $\pO{{m \over 2^k} \Vv,\Ve_i}$
for any integers $m$ and $k$, and so by continuity $\pO{a\Vv, \Ve_i}$ holds for all real $a$. Now, setting
$V_0 = V$ and $V_{i+1} = W_{i+1} \Inter V_i$, we see that each $V_{i+1}$ is a
subspace of $V_i$ of codimension 1. By the final hypothesis in our list, we
must have $V_{n+1} = 0$, and so $V$ must have dimension $n$.

Moreover, if we add the additional property $\|{1 \over n} (\Ve_1 + \cdots + \Ve_n)\|
= 1$, then the resulting sentence actually has a {\em unique} model up
to isomorphism, namely $\real^n$ with the 1-norm w.r.t. the usual basis
$\{\Ve_1,\ldots,\Ve_n\}$. For, by the remarks at the end of
Section~\ref{subsubsec:normed-spaces}, these revised hypotheses imply that
the $(n-1)$-simplex with vertex set $\{\Ve_1, \ldots, \Ve_n\}$ is
contained in the unit sphere. Also, there is a linear isometry
mapping $\Ve_i$ to $-\Ve_i$ and fixing the other basis elements,
which means that each of the $2^n$ $(n-1)$-simplices with vertex sets
$\{\pm\Ve_1, \ldots, \pm\Ve_n\}$ is contained in the unit sphere.
It follows that the unit sphere is the generalised
octahedron whose facets are these simplices and this
is the unit sphere of the 1-norm w.r.t the basis $\{\Ve_1, \ldots,
\Ve_n\}$.

Many characterizations of inner product spaces amongst normed spaces have been
discovered and rediscovered over the years, often based on abstractions of
orthogonality (our ``isosceles orthogonality'' $\pO{\Vv, \Vw}$ was proposed by
\citet{james-orthogonality}).  \citet{amir-characterizations} gives a
systematic presentation of some 350 characterizations involving
a wide range of ideas from geometry and analysis.
One characterization due to \citet{aronszajn-characterization} says that a normed
space is an inner product space if the norms of two sides and of one diagonal
of any parallelogram determine the norm of the other diagonal. 
Aronszajn's theorem implies that if we add the following additional
hypothesis to our original $\phi_n$, we obtain a purely additive characterization of euclidean
$n$-space: $$
\begin{array}{@{}l@{}l@{}}
\all{\Vv_1\;\Vw_1\;\Vv_2\;\Vw_2} &
        \Norm{\Vv_1} = \Norm{\Vv_2} \And
        \Norm{\Vw_1} = \Norm{\Vw_2} \And
        \Norm{\Vv_1 - \Vw_1} = \Norm{\Vv_2 - \Vw_2} \\
        & \quad\quad {} \Imp \Norm{\Vv_1 + \Vw_1} = \Norm{\Vv_2 + \Vw_2}.
\end{array}
$$
See \citet{amir-characterizations} or \citet{arthan-on-aronszajn} for a proof
of Aronszajn's characterization and see \citet{mok-metric-characterization} for
another interesting purely additive characterization.

\Subsubsection{Completeness in Metric Spaces and Normed Spaces}

In Section~\ref{sec:qelim} we will show that the theories of inner product
spaces and of Hilbert spaces coincide. In this section we investigate the
analogous question for metric spaces compared with complete metric spaces and
for normed spaces compared with Banach spaces and find, by contrast, that for
these theories the assumption of completeness does make a difference to the
first-order theory.

Consider the following properties of a relation $R$ between the real numbers
and the points of a metric space $X$.
\begin{itemize}
\item $R$ is a partial function whose domain comprises positive
numbers:
$$
  \all{x\; \Vp\; \Vq} R(x, \Vp) \And R(x, \Vq) \Imp x > 0 \And \Vp = \Vq.
$$
\item The domain of $R$ has no positive lower bound:
$$
  \all{\epsilon > 0} \ex{x\; \Vp} x < \epsilon \And R(x, \Vp).
$$
\item $R$ satisfies a form of the Cauchy criterion as its argument tends to 0:
$$
\all{\epsilon > 0} \ex{\delta > 0}\all{x\; y\; \Vp\; \Vq}
   x < \delta \And y < \delta \And R(x, \Vp) \And R(y, \Vq)
                     \Imp d(\Vp, \Vq) < \epsilon.
$$
\item $R$ has no limit as its argument tends to 0:
$$
\all{\Vq}\ex{\epsilon > 0} \all{\delta > 0} \ex{x\; \Vp} x < \delta \And R(x, \Vp) \And d(\Vp, \Vq) \ge \epsilon.
$$
\end{itemize}

Write $Q_R$ for the conjunction of the above properties and say
$R$ {\em represents} a sequence, $\Vs_n$, of points of $X$ iff there is a
strictly decreasing subsequence $x_n$ contained in the domain of $R$ such that
$x_n$ tends to $0$ as $n$ tends to $\infty$ and $R(x_n, \Vs_n)$ for all $n$.
Thus $Q_R$ implies that $R$ represents at least one Cauchy sequence but that no
Cauchy sequence represented by $R$ has a limit, so that $Q_R$ cannot hold in a
complete metric space.  Moreover if $Q_R$ holds and $R$ is definable in some
space $S$ by a formula $\pR{x, \Vp}$ of the language of metric spaces, then
the sentence asserting $\Not Q_{\sf R}$ belongs to the theory of complete metric spaces but not to the theory of metric spaces in general, since it does not hold in $S$.
A similar argument applies to normed spaces and Banach spaces,
the construction below being slightly complicated by the need for a parameter
in the formula $R(x, \Vv)$.

\begin{figure}
\begin{center}
\includegraphics[angle=0,scale=0.9]{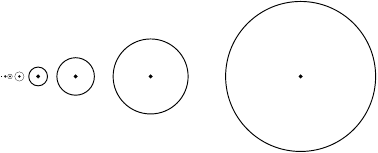}
\caption{The incomplete  metric space $\sM$}
\label{fig:I}
\end{center}
\end{figure}

For the metric space case, consider the subset $\sM$ of the real plane comprising
points $\Vp_n = (\frac{1}{2^n}, 0)$ and circles $C_n$ of radius $\frac{1}{2^{n+2}}$
with centre $\Vp_n$ for $n = 1, 2, \ldots$ (see Figure~\ref{fig:I}). Taking $\sM$
as a metric space under the euclidean metric, the sequence $\Vp_n$ is Cauchy
but has no limit in $\sM$, since its limit in the plane is the origin, which is
not in $\sM$. Define a predicate $\pP{x, \Vp}$ as follows:
$$
\pP{x, \Vp} \IsDef (\ex{\Vq} \Vq \not= \Vp \And d(\Vp, \Vq) = x) \And
    (\all{\Vq} \Vq \not= \Vp \Imp d(\Vp, \Vq) \ge x).
$$
I.e., $\pP{x, \Vp}$ holds iff $\Vp$ is an isolated point such that for some
$\Vq$, $x = d(\Vp, \Vq)$ is minimal for $\Vq \not= \Vp$, i.e., in $\sM$, iff $x =
\frac{1}{2^{n+2}}$ and $\Vp = \Vp_n$ for some $n$, $\Vq$ being any point of $C_n$.
Thus $\pP{x, \Vv}$ represents the divergent sequence $\Vp_n$ and $Q_{\Func{P}}$
holds in $\sM$ so that a first-order sentence asserting $\Not Q_{\Func{P}}$
holds in all complete metric spaces but not in $\sM$.

For the normed space case, we start with the vector space $\real^*$ of finitely
non-zero sequences of real numbers, which we think of as the union of the finite
dimensional spaces $\real^n$.  We will construct a normed space $\sY$ by making
modifications to the euclidean unit ball to make a certain divergent sequence
representable. So until further notice we work with
the euclidean metric on $\real^*$ which we write as $d(\Vx, \Vy) = \Norm{\Vx -
\Vy}$.  Also if $X$ is any non-empty subset of $\real^*$, we write $d(\Vv, X)$
for the distance between $\Vv$ and $X$, i.e., the infimum of the numbers
$d(\Vv, \Vx)$ as $\Vx$ ranges over $X$.

If $\Vv$ and $\Vw$ are distinct, non-antipodal unit vectors (i.e., $\Norm{\Vv}
= \Norm{\Vw} = 1$ and $\Vv \not= \pm \Vw$), the {\em great circle} through $\Vv$ and
$\Vw$ is defined to be the intersection of the unit sphere $S$ in $\real^*$ and
the plane through the origin spanned by $\Vv$ and $\Vw$.  Writing $\Ve_1,
\Ve_2, \ldots$ for the standard basis vectors, define a sequence of unit
vectors $\Vv_1, \Vv_2, \ldots$ as follows:
\begin{eqnarray*}
\Vv_1 &=& \Ve_1 \\
\Vv_{n+1} &=& \mbox{\parbox[t]{0.75\hsize}{
the unique
point on the great circle through $\Vv_n$ and $\Ve_{n+1}$ such that
$d(\Vv_{n+1}, \Vv_{n}) = \frac{1}{4^n}$ and $d(\Vv_{n+1}, \Ve_{n+1}) < d(\Vv_n,
\Ve_{n+1})$.
}}
\end{eqnarray*}
So each $\Vv_n$ lies in $\real^{n} \Diff \real^{n-1}$
and of the two points on the great circle at
distance $\frac{1}{4^n}$ from $\Vv_n$, $\Vv_{n+1}$ is the one on the same side of
$\real^{n}$ as $\Ve_{n+1}$ in $\real^{n+1}$.

It is a straightforward exercise in using the triangle inequality to prove
the following bounds on the distance between two members of the
sequence $\Vv_n$ (e.g., prove the upper bound first by induction on $k$
and then derive the lower bound using the upper bound for $d(\Vv_{n+1}, \Vv_{n+k})$).
\[
\frac{2}{3}\cdot\frac{1}{4^n} < d(\Vv_n, \Vv_{n+k}) < \frac{4}{3}\cdot\frac{1}{4^n} = \frac{1}{3}\cdot\frac{1}{4^{n-1}}
\]
These upper bounds show that the $\Vv_n$ form a Cauchy sequence.  Also if
$\alpha$ is the angle between $\Vv_n$ and $\Vv_{n+1}$, one has that
$d(\Vv_{n+1}, \real^n) = \Sin{\alpha} \ge \Sin{\frac{\alpha}{2}} =
\frac{1}{2}d(\Vv_{n}, \Vv_{n+1}) = \frac{1}{2}\frac{1}{4^n}$.  Whence using the
triangle inequality and the above bounds, we have
$d(\Vv_{n+k}, \real^n) \ge (\frac{1}{2}- \frac{1}{3})\frac{1}{4^n} =
\frac{1}{6}\frac{1}{4^n}$.  It follows that a limit of the $\Vv_n$ could not
belong to any $\real^n$, and so the sequence $\Vv_n$ has no limit in $\real^*$.

We have $d(\Ve_1, \Vv_n) = d(\Vv_1, \Vv_n) < \frac{1}{3}$ implying the
following bound for any $m, n \ge 1$.
\[
d(\Vv_m, -\Vv_n) > d(\Ve_1, -\Ve_1) - \frac{2}{3} = \frac{4}{3}.
\]

Let $O_n$ be the open disc with centre $\Vv_n$ and radius
$\frac{1}{2}\frac{1}{4^n}$.  Our estimates imply that the sets $O_1, -O_1, O_2,
-O_2, \ldots$ have pairwise disjoint closures. Let $E$ be the convex hull of
the set $A \cup \{\Vv_1, -\Vv_1, \Vv_2, -\Vv_2, \ldots \}$ where $A$ is the set
obtained from the (euclidean) unit disc $D$ in $\real^*$ by removing any points
that are within $\frac{1}{2}\frac{1}{4^n}$ of $\pm\Vv_n$, i.e., $A = D \Diff
\bigcup\{O_1, -O_1, O_2, -O_2, \ldots\}$.

$E$ satisfies the conditions for a unit disc in a normed space.  Let $T$ be the
unit sphere in this normed space, i.e., the boundary of $E$. Writing $S$ for
the unit sphere in $\real^{*}$, $T$ comprises $S \Diff \bigcup\{O_1, -O_1, O_2,
-O_2, \ldots\}$ together with a set of truncated cones made up of line segments
$[\pm\Vv_n, \Vw]$ joining each $\pm\Vv_n$ to each (euclidean) unit vector $\Vw$
such that $d(\pm\Vv_n, \Vw) = \frac{1}{2}\frac{1}{4^n}$.  Since the closures of
the sets $\pm O_n$ are pairwise disjoint, the $\Vv_n$ are the only isolated
extreme points of $T$ and the points on the open line segments $(\pm\Vv_n,
\Vw)$ are the only points of $T$ that are not extreme points.  (All these
claims are most easily seen by considering the possible ways in which $T$ can
intersect a plane through the origin).

Clearly, $\frac{1}{2}D \subseteq E \subseteq D$. Thus writing
$\Norm{\_}_X$ for the norm with unit disc $X$ (so $\Norm{\_}_D$ is the euclidean norm), we have that our two norms are equivalent in the sense that each is bounded by a constant multiple of the other:
\[
        2\Norm{\Vv}_D \ge \Norm{\Vv}_E \ge \Norm{\Vv}_D \ge \frac{1}{2}\Norm{\Vv}_E.
\]

As a consequence, under $\Norm{\_}_E$, just as under the euclidean norm, the
$\Vv_n$ form a Cauchy sequence that has no limit in the normed space $\sY$ whose
underlying vector space is $\real^{*}$ and whose unit disc is $E$.  Now let
$\pR{x, \Vv, \Ve}$ be a formula in the language of normed spaces expressing
the following properties:
\begin{description}
\item[$(i)$]  $\Norm{\Ve - \Vv} < \frac{2}{3}$;
\item[$(ii)$] $\Vv$ is an isolated point in the set of extreme points of the
unit disc;
\item[$(iii)$] there exists an extreme point $\Vw \not= \Vv$ of the unit disc
such that the line segment $[\Vv, \Vw]$ lies on the unit disc and $x =
\Norm{\Vw - \Vv}$.
\end{description}

In $\sY$, take $\Ve = \Ve_1 = \Vv_1$, and let $\Vv \in \sY$  and $x \in \real$ be
given. By the above estimates and remarks, conditions $(i)$ and $(ii)$
are satisfied iff $\Vv$ is one of the $\Vv_n$. If $\Vv = \Vv_n$,
then condition $(iii)$
is satisfied iff $x = \Norm{\Vv - \Vw}_E$ where $\Vw$ is a (euclidean) unit
vector with $\Norm{\Vv - \Vw}_D = \frac{1}{2}\frac{1}{4^n}$, and, for such a $\Vw$, we have:
\[
\frac{1}{4^n} = 2\Norm{\Vv - \Vw}_D \ge \Norm{\Vv - \Vw}_E = x \ge \Norm{\Vv - \Vw}_D = \frac{1}{2}\cdot\frac{1}{4^n}.
\]

We conclude that when the parameter $\Ve$ is interpreted by $\Ve_1$, the relation defined in $\sY$ by $\pR{x, \Vv, \Ve}$ represents the divergent sequence $\Vv_n$.
Thus $\ex{\Ve} Q_{\Func{R}}$ holds in $\sY$ and a sentence asserting
$\all{\Ve} \Not Q_{\Func{R}}$ holds in all Banach spaces but does not hold
in the normed space $\sY$.

In Section~\ref{sec:qelim}, we shall prove that for every set of sentences $A$
in the language of inner product spaces there is a subset $D$ of $\nat \cup
\{\infty\}$, such that an inner product space $V$ is a model of $A$ iff
$\Dim{V} \in D$ (see Corollary~\ref{:ip-axiomatizability}). So if $A$ is any set
of sentences in the language of metric spaces, then the class of metric space
models of $A$ cannot coincide with the class of complete metric spaces, since
if that were the case then the inner product space models of $A$ would comprise
precisely the class of Hilbert spaces, but this is impossible since, if $D$ is
the set of dimensions associated with $A$, either $\infty \not\in D$, so that
$A$ does not admit any infinite-dimensional Hilbert space as a model, or
$\infty \in D$, so that $A$ admits every infinite-dimensional inner product
space as a model and hence any incomplete
inner product space is a model of
$A$ (incomplete spaces being necessarily infinite-dimensional).  Essentially the same argument shows that no set of sentences in the
language of normed spaces can have the class of Banach spaces as its class of
models.

Collecting together the results of this section gives us the following theorem.

\begin{Theorem}\label{thm:incompleteness}
There are first-order sentences that hold in all complete metric spaces (resp.
Banach spaces) but not in all metric spaces (resp. normed spaces).  However, the
class of complete metric spaces (resp. Banach spaces) is not an axiomatizable
subclass of the class of metric spaces (resp. normed spaces).
\Done
\end{Theorem}


\Section{On undecidability in languages with a sort for the real numbers}\label{sec:prelim}

We will demonstrate the undecidability of various theories over
languages containing a sort for the real numbers by showing
how to interpret second-order arithmetic in them. In this section
we describe a general procedure for doing this.

\Subsection{Interpreting first-order arithmetic}

Consider a first-order language $L$ that includes
symbols for the field operations and the ordering relation
on a sort $\cal R$ whose intended interpretation is the ordered field $\real$,
e.g., our language $\LN$ for normed spaces. Let $\cal C$ be some
class of structures for $L$ in which the sort $\cal R$ has its intended
interpretation, e.g., the class of all Banach spaces is such a class for $\LN$.

Let a formula $\nu(x)$ of $L$ with one free variable of sort $\cal R$ be given.
The following sentence $\pPeano$ holds in a structure $\cal M$ in the class $\cal C$
iff in $\cal M$, $\nu(x)$ defines the set $\nat \subseteq \real$.
\begin{eqnarray*}
\pPeano &\IsDef&
 \BA \nu(0) \And {} \\
       (\all{x} \nu(x) \Imp x \ge 0 \And \nu(x + 1)) \And {} \\
       (\all{x\;y} \nu(x) \And \nu(y) \And x \not= y \Imp |x - y| \ge 1).
   \EA
\end{eqnarray*}

Now take any sentence $\phi$ in the language of first-order
arithmetic and reinterpret it as a sentence $\phi_{\nat}$ of $L$
by labelling all variables and constants in $\phi$
with sort $\cal R$ and relativizing all quantifiers using the formula $\nu(x)$,
i.e., replacing every subformula of the form $\ex{x}\psi$ by $\ex{x}\nu(x) \And \psi$
and every subformula of the form $\all{x}\psi$ by $\all{x}\nu(x) \Imp \psi$.

I claim that if $\pPeano$ is satisfiable in $\cal C$, then $\pPeano \Imp
\phi_{\nat}$ holds in $\cal C$ iff $\phi$ holds in $\nat$.  For, in any structure
with the intended interpretation of $\cal R$ in which $\pPeano$ holds,
$\phi_{\nat}$ holds iff $\phi$ holds in $\nat$.  So if $\pPeano$ holds in some
structure $\cal M \in {\cal C}$, then $\pPeano \Imp \phi_{\nat}$ holds in $\cal M$ iff $\phi$
is true, iff $\pPeano \Imp \phi_{\nat}$ holds in $\cal C$.  Thus, if we can find a
single model of the sentence $\pPeano$ in the class $\cal C$, then the theory
of $\cal C$ must be undecidable, since a decision procedure for it would lead
to a decision procedure for the set of
truths of first-order arithmetic, contradicting Tarski's
theorem on the undefinability of truth.

This method of relativization has often been used to show that extending
decidable theories such as Presburger arithmetic or the theory
of a real closed field with a new
uninterpreted unary function or predicate leads to undecidability
\citep{tarski-undecidable,downey-presmon}. Even though our $\nu(x)$ is not just an
uninterpreted unary predicate but rather a complex formula in a language with a
constrained interpretation, we have to exhibit just {\em
one} model of the characterizing sentence $\pPeano$ in order to get a
reduction of first-order arithmetic to the theory of the class $\cal C$.

\Subsection{Interpreting second-order arithmetic}

We will obtain still stronger undecidability results by observing that in a
first-order theory of the real numbers with a predicate for the natural numbers, one can
interpret not only first-order arithmetic as we did above but even {\em
second-order} arithmetic. This is ``well-known'' but since we know of no
reference for it in the literature we will give the proof. The setting
is as in the previous section with $L$ a language including a sort $\cal R$
for the real numbers with the usual operations and $\cal C$ a class of structures
for $L$ in which these things have their intended interpretation.

First we will briefly describe second-order arithmetic; see, e.g.,
\citet{simpson-subsystems} for more details. The language $\A^2$ of
second-order arithmetic is a 2-sorted language with a sort $\NatSort$ called
``type 0'' whose intended interpretation is the set of natural numbers
$\nat$ and a  sort $\SetSort$ called ``type 1'' whose intended interpretation
is the set $\Pow{\nat}$ of all sets of natural numbers.  The expressions are
those of first-order arithmetic which have sort $\NatSort$ together with
variables of sort $\SetSort$. Atomic formulas can be built from numeric terms
by the usual predicates of first-order arithmetic, and also if $t$ is a numeric
term and $A$ a set variable we can form the atomic formula $t \in A$.
Quantification is allowed over both numeric and set variables.

We have already seen how to interpret first-order arithmetic by relativizing
quantifiers using the natural number predicate $\nu(x)$. In order to interpret
type 1 variables in the first-order theory of the real numbers, we use the mapping taking a
set $A$ with characteristic function $\chi_A$:
$$ \chi_A(n) = \left\{ \begin{array}{ll}
                        1 & \mbox{if $n \in A$}     \\
                        0 & \mbox{otherwise}        \\
                     \end{array} \right.
$$
into the real number whose ternary expansion is determined by the values
$\chi_A(n)$:
$$ \sharp A = \sum_{n=0}^{\infty} \chi_A(n) / 3^n .
$$

Note that a binary version of the same method would not give an injective map
because of $1.000 \cdots = 0.111 \cdots$ etc., and so would require workarounds
like treating terminating expansions differently or encoding the function in
even digits of the binary expansion. Using ternary, we can straightforwardly
and unambiguously recover the set $A$ from the number $\sharp A$. Let $h_n(x)$
be the value of the first $n$ ternary digits of $x$, considered as an integer,
where the zeroth `digit' is simply $\floor{x}$:
$$ h_n(x) = \floor{3^n x} .
$$

Then defining
$$ d_n(x) = \left\{ \begin{array}{ll}
                        h_0(x) & \mbox{if $n = 0$}               \\
                        h_n(x) - 3 h_{n-1}(x) & \mbox{otherwise} \\
                     \end{array} \right.
$$
we have $d_n(\sharp A) = \chi_A(n)$. We will show below that the function
$d_n(x)$ is definable, or more precisely that we can find a formula $\pD(n,x)$ of
our language with two free variables whose interpretation corresponds to
$d_n(x) = 1$ in all standard models (i.e. those interpreting the real sort in
the usual way).

Assume that $\cal M$ is a structure for the language $L$ and that $\nu(x)$ is
a formula in $L$ with the indicated free variable of sort $\ScaSort$
which defines the natural numbers in $\cal M$, i.e., $\nu(x)$ holds in $\cal M$
iff $x$ is interpreted as a natural number.
Then the relation $\pD(n,x)$ can be defined in terms of $\nu(x)$
using the following
relational translations of the definitions given above, first for $h_n(x)$:
$$ h_n(x) = l \Iff
   \nu(n) \And \nu(l) \And \ex{k} \nu(k) \And
   3^n = k \And l \leq k \cdot x \And k \cdot x < l + 1
$$
\noindent then $d_n(x)$:
$$ d_n(x) = y \Iff
   \BA \nu(n) \And \nu(y) \And \\
       \BA((n = 0 \And h_0(x) = y) \Or \\
           \;(\ex{m\;l\;k} \BA \nu(m) \And \nu(l) \And \nu(k) \And \\
                              n = m + 1 \And
                              h_n(x) = l \And h_m(x) = k \And \\
                              l = y + 3 \cdot k)) \EA \EA
   \EA
$$
\noindent and finally:
$$
\pD(n,x) \Iff d_n(x) = 1.
$$

It remains to define the exponential relation $3^n = k$ in $L$, but this can be
done by taking any of the usual definitions in the language of first-order
arithmetic, e.g. the one given by \citet{smullyan-godel}, and translating into
$L$ using the numeric sort $\ScaSort$ and its operations and relativizing with
respect to the predicate $\nu(x)$.

If we define:
$$
\pS(x) \IsDef x \ge 0 \And \all{y} y \ge 0 \And (\all{n}\pD(n, x) \Iff \pD(n, y)) \Imp x \le
y
$$
\noindent then $\pS(x)$ holds iff $x = \sharp \{n \in \nat \mid \pD(n, x)\}$.
Thus, we can interpret second-order arithmetic in $L$
using $\pD(n, x)$ to represent sets of natural
numbers as real numbers and using $\pS(x)$ to pick a canonical representative:
given a formula  $\phi$ of second-order arithmetic, we take
$\phi^*$ to be the result of the following sequence of transformations:

\begin{enumerate}

\item Replace subformulas of the form
$\ex{x^{\cal N}}\psi$ by $\ex{x^{\cal R}}\nu(x^{\cal R}) \And \psi$
and subformulas of the form
$\all{x^{\cal N}}\psi$ by $\all{x^{\cal R}}\nu(x^{\cal R}) \Imp \psi$;

\item Replace subformulas of the form
$\ex{A^{\cal P}}\psi$ by $\ex{A^{\cal R}}\pS(A^{\cal R}) \And \psi$
and subformulas of the form
$\all{A^{\cal P}}\psi$ by $\all{A^{\cal R}}\pS(A^{\cal R}) \Imp \psi$;

\item Replace remaining occurrences of the sort labels $\cal N$
and $\cal P$ by $\cal R$;

\item Replace subformulas of the form $t \in A$ by $\pD(t, A)$.

\end{enumerate}

Here recall that each variable and constant in our many-sorted language
comprises a name labelled with a sort, which we write as a superscript, and note
that there are no constants of sort $\cal P$.
Now given a {\em sentence} $\phi$ of second-order arithmetic, we
may assume (up to a logical equivalence) that bound variables have been renamed
if necessary so that no variable name appears in $\phi$ with two different sorts and
distinct variables remain distinct even after a relabelling that identifies two
sorts.  Assuming that $\nu(x)$ does indeed define the natural numbers, we then
find by induction on the structure of a formula in which no variable name
appears with two different sorts that the sentence $\phi$ is true iff $\phi^*$ holds
in the structure $\cal M$.  The details of the induction are
straightforward: in the inductive step for the type 1 quantifiers, one notes
that by the discussion above, $\sharp$ defines a 1-1 correspondence between
$\Pow{\nat}$ and the set of real numbers $s$ such that $\pS(s)$ holds.

\begin{Theorem}\label{thm:mult-undec}
Let $L$ be a (many-sorted) first-order language including a sort $\cal R$,
constants $0 : {\cal R}$ and $1 : {\cal R}$ and function symbols $\_+\_,
\_\times\_ : {\cal R} \times {\cal R} \rightarrow {\cal R}$ whose intended
interpretations form the field of the real numbers.  Let $\cal C$ be some class
of structures for $L$ in which $\cal R$ and these symbols have their intended
interpretations and let $\cal T$ be the theory of $\cal C$, i.e., the set of
all sentences that hold in every member of $\cal C$.  If there is a formula $\nu(x)$
of $L$ with one free variable $x$ of sort $\cal R$ such that in some structure
$\cal M$
in the class $\cal C$, $\nu(x)$ defines the set of natural numbers, then there is
a primitive recursive reduction of second-order arithmetic to $\cal T$.
\end{Theorem}
\Proof
The reduction maps a sentence $\phi$ of second-order arithmetic to the sentence $\phi' \IsDef \pPeano \Imp \phi^{*}$
where $\pPeano$ is defined as above using the $\nu(x)$ that we are given
by hypothesis and $\phi^{*}$ is the above translation of $\phi$ into the language $L$.
By the discussion above, $\phi'$ then holds in every member of $\cal C$ iff $\phi$ is true.
\Done

\Subsection{Interpretation in an additive theory}

Since the {\em linear} theory of integer arithmetic is decidable
\citep{presburger} we need multiplication in our language in order to interpret
the full, undecidable theory, even though the characterizing formula $\pPeano$
itself does not involve multiplication. But we will later want to show the
undecidability of additive theories of metric and vector spaces where
multiplication is not available. In some interesting cases we can construct a
structure in which we can define not only the natural numbers but also the
graph of the multiplication function $(x, y) \mapsto xy$.  In order to
interpret first-order arithmetic we only need to be able to define and
characterize the multiplication of natural numbers. But to achieve the full
reduction of second-order arithmetic, we require multiplication of arbitrary
real numbers, since this is used in the formulas defining $h_n(x) = l$ and
$d_n(x) = l$ above.

To make this programme work, we need an analogue $\pMult$ of the sentence
$\pPeano$, asserting that a formula $\mu(x,y,z)$ with three free variables
defines the multiplication relation $x \cdot y = z$ on the real numbers.
Let us define $\pMult$ as follows:
\begin{eqnarray*}
\pMult &\IsDef&
\BA (\all{x\;y}\exu{z} \mu(x,y,z)) \And {}\\
      (\all{x\;y\;z} \mu(x,y,z) \Imp \mu(y,x,z)) \And {}\\
      (\all{y\;z} \mu(0,y,z) \Iff z = 0) \And {}\\
      (\all{y\;z} \mu(1,y,z) \Iff z = y) \And {}\\
      (\all{x_1\;x_2\;y\;z_1\;z_2} \mu(x_1,y,z_1) \And \mu(x_2,y,z_2)
                                   \Imp {}\\
        \quad \quad \mu(x_1+x_2,y,z_1+z_2)) \And {}\\
      (\all{x\;y\;z}\all{\epsilon > 0}
        \mu(x,y,z)
        \Imp \ex{\delta> 0} \all{x'\;z'} \\
        \quad\quad  |x - x'| < \delta \And \mu(x',y,z')
                                 \Imp |z - z'| < \epsilon).

  \EA
\end{eqnarray*}

The first conjunct asserts that $\mu(x, y, z)$ does indeed define a function $f(x,y)
= z$, and the second that $f(x,y) = f(y,x)$. The next three conjuncts ensure that
this function coincides with multiplication in the case where $x$ is a natural
number because they give $f(0,y) = 0$, $f(1,y) = y$ and $f(x+1,y) = f(x,y) +
y$. They also imply that this holds for $x \in \num$, because $f(-x,y) + f(x,y)
= 0$ and therefore $f(-x,y) = -f(x,y)$. Using the additivity property
repeatedly we also see that for any real number $x$ and natural number $q > 0$ we
have $f(x,y) = f(x/q + \cdots + x/q,y) = f(x/q,y) + \cdots + f(x/q,y) = q \cdot
f(x/q,y)$ and therefore $f(x/q,y) = f(x,y)/q$. Together these imply that
$f(x,y) = x \cdot y$ when $x \in \rat$. Now the final `continuity' conjunct
implies that for any $y \in \real$ the function $g(x) = f(x,y) - x \cdot y$ is
continuous.  Since $\{x \mid g(x) \not= 0 \}$ is the preimage of the open set
$\real - \{0\}$ under a continuous function, it is open. Since it contains no
rational numbers, it must be empty, so $g(x)$ is identically zero as
required.

Hence our characterizing formula $\pMult$ works as claimed and we can
summarize the import of all this in the following result:

\begin{Theorem}\label{thm:linear-undec}
Let $L$ be a (many-sorted) first-order language including a sort $\cal
R$, together with function symbol $\_+\_ : {\cal R} \times {\cal R}
\rightarrow {\cal R}$, a binary predicate symbol $\_<\_$ on the sort
$\cal R$ and a constant $1 : \cal R$ whose intended interpretations form
the ordered group of real numbers under addition with $1$ as a
distinguished positive element.  Let $\cal C$ be some class of
structures for $L$ in which $\cal R$ and these symbols have their
intended interpretations and let $\cal T$ be the theory of $\cal C$,
i.e., the set of all sentences that hold in every member of $\cal C$.  Let
$\nu(x)$ (resp. $\mu(x, y, z)$) be a formula of $L$ with one free variable
$x$ of sort $\cal R$ (resp. free variables $x$, $y$ and $z$ all of sort
$\cal R$).  If in some structure in the class $\cal C$,  $\nu(x)$ defines
the set of natural numbers with the intended interpretation of the
constant 1 and $\mu(x, y, z)$ defines the multiplication relation on the
set of real numbers, then there is a primitive recursive reduction of
second-order arithmetic to $\cal T$.
\end{Theorem}
\Proof
Given a sentence $\phi$ in second-order arithmetic,
let $\phi^{*}$ be the translation of $\phi$ into the language $L$ used in the proof
of Theorem~\ref{thm:mult-undec}.
There is a primitive recursive function that maps any formula in $L$ to a
logically equivalent one in which all instances of multiplication are unnested,
i.e., multiplication only appears in atomic predicates of the
form $xy = z$ where $x$, $y$ and $z$ are variables; see, e.g., \citet{hodges-model}.
Let $\phi^{+}$ be the result of applying this function to $\phi^{*}$ and then
replacing each atomic predicate of the form $xy = z$ by $\mu(x, y, z)$.
If we then set $\phi'' \IsDef \pPeano \And \pMult \Imp \phi^{+}$,
$\phi''$ holds in every member of $\cal C$ iff $\phi$ is true.
\Done

In fact, both Theorems~\ref{thm:mult-undec} and~\ref{thm:linear-undec}
can easily be strengthened to allow the formulas $\nu(x)$ and $\mu(x, y, z)$
to have additional free variables acting as parameters: if for some structure
and some choice of values for the parameters, $\nu(x)$ defines the natural numbers, then the conclusion of Theorem~\ref{thm:mult-undec} will obtain, while if
also $\mu(x, y, z)$ defines the graph of multiplication, then the conclusion
of Theorem~\ref{thm:linear-undec} will also obtain. The formulations without
parameters are all we need in the sequel.

\Section{Metric spaces}\label{sec:metric-spaces}

We begin the main work of this paper by considering metric spaces. The
generality of the metric space axioms gives us considerable freedom
to construct spaces in which various arithmetic sets and relations are
definable as needed to apply the methods of Section~\ref{sec:prelim}.
Many of the same ideas will appear later
for normed spaces but in a more intricate form.

The elementary theory of metric spaces is known to be undecidable.
This was first proved by \citet{bondi:dec-ms}.
\citet{kurz-metric} give a very simple proof by encoding an
arbitrary reflexive symmetric binary relation $R$ (i.e. an undirected graph) as
a metric via:
$$ d(x,y) = \left\{ \begin{array}{ll}
                     0 & \mbox{if $x = y$} \\
                     1 & \mbox{if $x \not= y \And R(x,y)$} \\
                     2 & \mbox{if $\Not R(x,y)$.}
                 \end{array} \right. $$

This allows the decision problem for the theory of a reflexive symmetric
binary relation, known to be hereditarily undecidable \citep{rabin-simple}, to be
reduced to the theory of metric spaces. In this proof, few special properties
of $\real$ are needed and almost any other set of valuations would work; the
set of points takes centre stage and the set of scalars plays only
a supporting role.

The theory of a reflexive symmetric binary relation is undecidable, but
is recursively (indeed finitely) axiomatizable.
The arguments of~\citet{bondi:dec-ms} and~\citet{kurz-metric} do not preclude the possibility
that the theory of metric spaces might be recursively axiomatizable.
By exploiting the methods of Section~\ref{sec:prelim}, we obtain
a much stronger result:

\begin{Theorem}\label{thm:ms-undecidable}
There is a primitive recursive reduction of second-order arithmetic to the
theory of metric spaces $\MS$.
\end{Theorem}
\Proof
Let $\num$ be the set of
integers with the usual metric $d(x,y) = |x - y|$.
Clearly in this metric space the formula:
$$ \pN{x} \IsDef \ex{a\;b} d(a,b) = x$$
\noindent defines the natural numbers as a subset of the real numbers.
Applying Theorem~\ref{thm:mult-undec} completes the proof \Done

The theory of metric spaces, is therefore not arithmetical, i.e., it is not
definable by any formula of first-order arithmetic, and hence it is not
recursively enumerable and it is not recursively axiomatizable.

If $K$ is an ordered field, define a {\em metric space over $K$} to be a
structure for the language $\LM$ of metric spaces in which the scalar sort
$\ScaSort$ and its operations are interpreted in $K$ and which satisfies the
metric space axioms.  Let $\cal C$ be the class of all structures for $\LM$
that are metric spaces over $K$, where $K$ ranges over all real closed fields.
Then $\cal C$ is clearly a recursively axiomatizable class of structures and so
the set of sentences of $\LM$ that are valid in $\cal C$ is recursively
enumerable.  Given Theorem~\ref{thm:ms-undecidable}, we must conclude that
there is a real closed field $K$ and a sentence of $\LM$ that holds in any
metric space over $\real$ but fails in some metric space over $K$.

\begin{figure}
\begin{center}
\includegraphics[angle=0,scale=0.9]{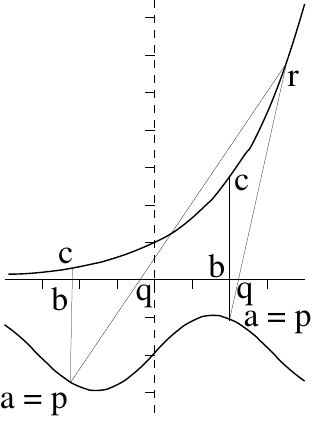}
\caption{Defining $\Func{exp}$ and $\Func{sin}$ in the  metric space $\sG$}
\label{fig:G}
\end{center}
\end{figure}

The situation is much the same even if we disallow multiplication:

\begin{Theorem}\label{thm:ms-additive-undecidable}
There is a primitive recursive reduction of second-order arithmetic to the
additive theory of metric spaces $\MSA$ .
\end{Theorem}
\Proof
We will exhibit a metric space $\sG$ such
that the set of natural numbers and the graph of the real multiplication
function are {\em additively definable} in $\sG$ i.e., definable using formulas that
do not involve multiplication. $\sG$ is the subspace of the euclidean
plane comprising the $x$-axis together with the graphs of two functions $e$ and
$s$ where $e$ is the exponential function, $e(x) = \Exp{x}$, and $s$ is defined
by $s(x) = \Sin{x} - 2$. Thus $\sG$ has three connected components: the graph of
$e$ lying strictly above the $x$-axis, the $x$-axis itself and the graph of $s$
lying strictly below the $x$-axis, as illustrated in Figure~\ref{fig:G} (which 
actually shows $\Exp{x/2}$ rather than $\Exp{x}$ for reasons of space).

Our first task is to show that the connected components of $\sG$ are additively
definable.  In the euclidean plane, a point $\Vq$ lies on the line segment
$[\Vp, \Vr]$ iff $d(\Vp, \Vr) = d(\Vp, \Vq) + d(\Vq, \Vr)$. A point $\Vp$ of
$\sG$ lies on the $x$-axis iff $\sG$ contains the entire line segment $[\Vp, \Vq]$
for some $\Vq \not= \Vp$. So the $x$-axis is additively definable in $\sG$. Now
if $f$ is a real-valued function of a real variable and $x$ is any real number,
then $(x, 0)$ is the point on the $x$-axis nearest to the point $(x, f(x))$ on
the graph of $f$.  Therefore, if $\Vp$ is a point of $\sG$ and $\Vq$ is the point
on the $x$-axis nearest to $\Vp$, then $d(\Vp, \Vq) > 3$ iff $\Vp = (x, e(x))$
for some $x$ with $e(x) > 3$, so the set of such $\Vp$ is additively definable.
But then the graph of $s$ comprises precisely those points $\Vp$ of $\sG$ for
which there are a point $\Vr = (x, e(x))$ with $e(x) > 3$ and a point $\Vq \not= \Vp$ on
the $x$-axis such that $\Vq$ lies on the line segment $[\Vp, \Vr]$ (see
Figure~\ref{fig:G}). Thus the graph of $s$ is additively definable and hence so
is the graph of $e$ (which comprises the points of $\sG$ that are neither on the
$x$-axis nor on the graph of $s$).

The point $\VO = (0, 0)$ is now additively definable in $\sG$ as the point on the
$x$-axis for which there is a point $\Vp$ on the graph of $e$ with $d(\Vp, \VO)
= 1$ and $d(\Vp, \Vq) > 1$ for every other point $\Vq$ on the $x$-axis. The
functions $\Func{exp}$ and $\Func{sin}$ are then additively definable: given a
real number $t$, there are collinear points $\Va$, $\Vb$ and $\Vc$ with $\Va$
on the graph of $s$, $\Vc$ on the graph of $e$ and $\Vb$ the point on the
$x$-axis closest to $\Vc$ with $d(\VO, \Vb) = |t|$ and with $d(\Vb, \Vc) \ge 1$
if $t \ge 0$ and $d(\Vb, \Vc) < 1$ if $t < 0$ (see Figure~\ref{fig:G}). With
this unique choice of $\Va$, $\Vb,$ and $\Vc$, $\Exp{t} = d(\Vb, \Vc)$ and
$\Sin{t} = 2 - d(\Vb, \Va)$.

For positive $x$, we may now define $\Log{x}$ by $\Exp{\Log{x}} = x$, then
define multiplication for positive real numbers using $xy = \Exp{\Log{x} + \Log{y}}$
and extend the definition to all real numbers using $0y = x0 = 0$, $(-x)y = x(-y) =
-xy$ and $(-x)(-y) = xy$.  The real number $\pi$ is additively definable as the
smallest $x > 0$ such that $\Sin{x} = 0$ and then the natural numbers are
additively definable as the set of $n \ge 0$ such that $\Sin{n\pi} = 0$. Thus
multiplication and the natural numbers are additively definable in the metric
space $\sG$ and we may conclude by Theorem~\ref{thm:linear-undec} that there is
a primitive recursive reduction of second-order arithmetic to the additive
theory of any class of metric spaces including $\sG$.
\Done

\Subsection{Decidability of the $\fAE$ fragment}

A sentence is said to be $\fAE$ if it is in
prenex normal form with no universal quantifier in
the scope of an existential one, i.e. it has the following form for some $n \geq 0$,
and $m \geq 0$ with $\phi$ quantifier-free:
$$ \all{x_1 \ldots x_n} \ex{y_1 \ldots y_m} \phi $$

\noindent the set of $\fEA$ sentences being defined analogously exchanging
`$\forall$' with `$\exists$'.

The set of valid first-order $\fAE$
sentences with no function symbols is decidable
\citep{bernays-schonfinkel}: in fact, such a sentence with $n$
initial universal quantifiers is valid
iff it holds in all interpretations with at most $\Max\{n, 1\}$ elements;
but then it is a finite problem to enumerate all such interpretations.
By working in many-sorted logic, this can be generalized to
some important cases where function symbols occur
\citep{fontaine-thesis}. We will prove the decidability of the set of valid
$\fAE$ sentences in the language of metric spaces using similar ideas
exploiting the fact that if $K \subseteq M$ and $d$ is a metric on $M$
then the restriction of $d$ to $K \times K$ is also a metric on $K$.
In fact our decision procedure will decide validity for a superset of the $\fAE$ sentences.
We say a sentence is:

\begin{itemize}

\item $\fAE_p$ if it is prenex and no universal quantifier over points is in the scope of an
existential quantifier (of any sort);

\item $\fEA_p$ if it is prenex and no existential quantifier over points is in the scope of a
universal quantifier (of any sort).

\end{itemize}

We have the following analogue of the theorem of Bernays and
Sch\"{o}nfinkel:

\begin{Theorem}\label{thm:ea-v-finite-models}
Let $\phi$ be an $\fEA_p$ sentence in the language of metric spaces,
and let $n$ be the number of existential quantifiers of
the point sort in $\phi$.
Then $\phi$ is satisfiable in a metric space iff it is satisfiable in a
finite metric space with no more than $\Max\{n, 1\}$ points.
\end{Theorem}
\Proof
The right-to-left direction of the theorem is immediate.
For the left-to-right direction, assume that the $\fEA_p$
sentence $\phi$ is satisfiable in some metric space $M$. As
existential quantifiers commute up to logical equivalence, we can assume
without loss of generality that $\phi$ consists of a block of $n \geq 0$ {\em existential} quantifiers over points followed by a block
comprising {\em universal} quantifiers over points and scalar quantifiers of
either kind.
We write this as follows:
$$ \phi \equiv \ex{\Vx_1 \ldots \Vx_n} \all{\overline{\Vy} / \Func{Q}\overline{z}} \psi $$
where $\psi$ is a quantifier-free formula whose free variables are contained in
$$
\{\Vx_1, \ldots, \Vx_n, \Vy_1, \ldots, \Vy_k, z_1, \ldots, z_l\}.
$$
If $n = 0$, we may replace $\phi$ by the logically equivalent formula $\ex{\Vx}\phi$
(hence replacing $n$ by $1 = \Max\{n, 1\}$), and so we may assume that
$n \ge 1$.
We have that $  \rho \IsDef \all{\overline{\Vy} / \Func{Q} \overline{z}}\psi$ holds
for some points $\Vx_1, \ldots, \Vx_n \in M$.  But then {\it a fortiori}, $\rho$
and hence $\phi$ hold in the subspace $K = \{\Vx_1, \ldots, \Vx_n\}$ of $M$. But
$K$ has at most $n$ points and we are done.
\Done

\begin{Corollary}
An $\fAE_p$ sentence in the language of metric spaces with $n$ universally
quantified point variables, which we can write as $$
\all{{\Vx}_1 \ldots {\Vx}_n} {\exists{\overline{\Vy}} / \Func{Q}
\overline{z}}.\; \phi $$

\noindent holds in all metric spaces iff it holds in all finite metric
spaces with at most $\Max\{n, 1\}$ points.
\end{Corollary}
\Proof
Apply the theorem to the negation of the sentence.
\Done

These ideas lead to a decision procedure for valid $\fAE_p$ sentences:

\begin{Theorem}\label{thm:ea-v-decidable}
The set of valid $\fAE_p$ sentences in the language of metric
spaces is decidable.
\end{Theorem}
\Proof
Since $\phi$ is valid iff $\Not \phi$ is not satisfiable, it suffices
to describe a decision procedure for satisfiable $\fEA_p$ sentences.
If $\phi$ is an $\fEA_p$ sentence, then
as in the proof of Theorem~\ref{thm:ea-v-finite-models}, we may assume $\phi$
has the form
$\ex{\Vx_1 \ldots \Vx_n} \all{\overline{\Vy} / \Func{Q}\overline{z}} \psi $
where $\psi$ is quantifier-free and $n \ge 1$,
and then $\phi$ is satisfiable iff it is satisfiable in a metric space comprising
just the interpretations of $\Vx_1, \ldots, \Vx_n$ under a satisfying
assignment for $\all{\overline{\Vy} / \Func{Q}\overline{z}} \psi $.
So if we replace each subformula of $\phi$ of the form $\all{\Vy}\rho$,
by the conjunction $\rho[\Vx_1/\Vy] \And \ldots \And \rho[\Vx_n/\Vy]$ we obtain
a sentence that is equisatisfiable with $\phi$ and has no point universal
quantifiers. So we may assume $\phi$ has the form
$\ex{\Vx_1 \ldots \Vx_n} \psi$ where $\psi$ contains only scalar quantifiers.
Now if $M$ is a finite metric space with $n$ points $\Vp_1, \ldots, \Vp_n$,
say, define a function $f_M : M \mapsto \real^n$ by $f_M(\Vp) = (d(\Vp, \Vp_1),
\ldots, d(\Vp, \Vp_n))$.  If we equip $\real^n$ with the metric $d_{\infty}$
induced from the $\infty$-norm, $d_{\infty}(\Vv, \Vw) = \Max\{\Abs{\Vv_i - \Vw_i}
\ST 1 \le i \le n\}$, then it is easy to check that $f_M$ is an isometric
embedding of $M$ in $(\real^n, d_{\infty})$.  It follows that $\ex{\Vx_1
\ldots \Vx_n}\psi$ is satisfiable in general iff it is satisfiable in
$(\real^n, d_{\infty})$.  Thus if we choose fresh variables $x_{ij}$, $1 \le i,
j \le n$, and let $\psi'$ be the result of replacing each subterm $\Vx_s =
\Vx_t$ in $\psi$ by $x_{s1} = x_{t1} \And \ldots \And x_{sn} = x_{tn}$ and each
subterm $d(\Vx_s, \Vx_t)$ by $\Max\{\Abs{x_{s1} - x_{t1}}, \ldots, \Abs{x_{sn}
- x_{tn}}\}$, then $\ex{\Vx_1 \ldots \Vx_n}\psi$ is satisfiable iff $\phi'
  \IsDef \ex{x_{11}\;x_{12} \ldots x_{nn}}\psi'$ is satisfiable.  But $\phi'$
contains no point variables so we may apply a decision procedure for real
closed fields to complete the proof.
\Done

\Subsection{Undecidability of the $\fEA$ fragment}
The following result shows that Theorem~\ref{thm:ea-v-decidable} is
the best possible decidability result of its type:

\begin{Theorem}\label{thm:ms-ea-valid-undec}
If $\cal C$ is any class of metric spaces that includes the metric space
$\num$, then the set of $\fEA$ sentences that are valid in $\cal C$ is
undecidable.
\end{Theorem}
\Proof
We will prove the equivalent claim that the set of $\fAE$ sentences that are
satisfiable in $\cal C$ is undecidable.
Note that the formula $\pN{x}$ used in the proof of
Theorem~\ref{thm:ms-undecidable} is purely existential, and so the corresponding
sentence $\pPeano$ of Section~\ref{sec:prelim}
is logically equivalent to an $\fAE$ sentence.
Let $\phi(x_1, \ldots, x_k)$ be a quantifier-free formula in the language of arithmetic
and consider the following sentence in the language of metric spaces:
\begin{eqnarray*}
\phi_1 &\IsDef& \pPeano \And \ex{x_1 \ldots x_k} \pN{x_1} \And \ldots \And \pN{x_k} \And \phi(x_1, \ldots, x_k).
\end{eqnarray*}
$\phi_1$ is logically equivalent to an $\fAE$ sentence and $\phi_1$ is satisfiable in $\num$ and
hence in $\cal C$ iff $\phi(x_1, \ldots, x_k)$ is satisfiable over the
natural numbers. Thus a decision procedure for $\fAE$ sentences that are
satisfiable in $\cal C$ would lead to a decision procedure for satisfiability
of quantifier-free formulas in arithmetic and, in particular, for systems of
Diophantine equations, contradicting the famous resolution of Hilbert's $10^{\mbox{\scriptsize th}}$
problem by \citet{matiyasevich-paper}.
\Done

\Section{Undecidability of theories of normed spaces}\label{sec:undec}

The theory $\NS{1}$ of 1-dimensional normed spaces reduces easily to the theory
of the real numbers, since every such space is isomorphic to $\real$ with absolute
value as the norm.  Thus $\NS{1}$ is decidable. We will show in this section that this is the strongest possible
positive decidability result of its type: even the additive theory $\NSA{2}$ of
2-dimensional normed spaces is undecidable.  In fact, $\NSA{2}$  is not even arithmetical.

The main argument giving undecidability is in Section~\ref{subsec:undec-reducing}.
We exhibit a 2-dimensional normed space, $\sX$, and describe
geometric constructions in that space of the set of natural numbers and of the
the graph of the multiplication function. Formalising these constructions
in the additive language of normed spaces and
applying the methods of Section~\ref{sec:prelim} immediately gives
a reduction of second-order arithmetic to the (additive) theory of any class
of normed spaces including $\sX$. Taking a product with a Hilbert space of
appropriate dimension, the construction lifts into any desired dimension
${}\ge2$.

In Section~\ref{subsec:many-one-degrees}, we obtain tighter estimates of the
degrees of unsolvability of the normed space theories.  We prove a kind of
Skolem-L\"{o}wenheim theorem for normed spaces and use it to give reductions of
the normed space theories to fragments of third-order arithmetic.
We find that for any integer $d \ge 2$, the theory $\NS{d}$ is many-one
equivalent to second-order arithmetic, as is the theory $\NS{\fin}$ of all
finite-dimensional normed spaces.  We then strengthen the results of
Section~\ref{subsec:undec-reducing}: using a variant of the approach
of Section~\ref{sec:prelim}, we show that the theory $\NS{\infty}$ of
infinite-dimensional normed spaces and the theory $\NS{\relax}$ of all normed
spaces are both many-one equivalent to the set of true $\Pi^2_1$ sentences in
third-order arithmetic. All of this goes through for the purely additive
theories with little extra work. The results also hold equally well for
Banach spaces: even though, by Theorem~\ref{thm:incompleteness}, the theory
$\NS{\relax}$ of all normed spaces is a proper subset of the theory
$\BS{\relax}$ of all Banach spaces, the two theories turn out to be
many-one equivalent.

\Subsection{Reducing second-order arithmetic to the theory of a normed space}\label{subsec:undec-reducing}

\begin{figure}
\begin{center}
\includegraphics[angle=0,scale=0.9]{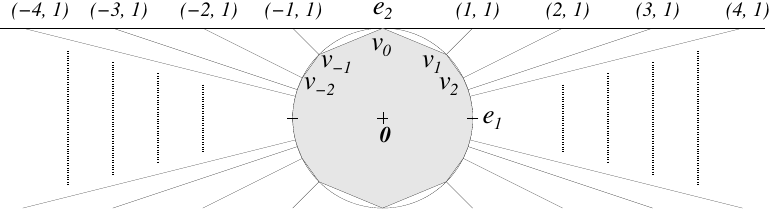}
\caption{The unit disc $D$ in the space $\sX$}
\label{fig:infinigon}
\end{center}
\end{figure}

To apply the results of Section~\ref{sec:prelim}, we will exhibit a particular
2-dimensional normed space, $\sX$, and give additive predicates that, in $\sX$,
define the natural numbers as a subset of the scalars and the graph of the
scalar multiplication function.  We define the norm by describing its unit
disc.  Let $C$ be the unit circle in $\real^2$ with respect to the standard
euclidean norm. For each $i \in \zahl$, let $l_i$ be the line passing through
$\VO$ and the point  $(i, 1)$. Then $l_i$ meets $C$ in two points $\Vv_i$, say,
in the upper half-plane and $-\Vv_i$ in the lower (see
Figure~\ref{fig:infinigon}). The set $E$ comprising the $\pm\Vv_i$ together with
the two points $\Ve_1 = (1, 0)$ and $-\Ve_1$ is a closed and bounded subset of
$\real^2$ and is symmetric about the origin.  If we write $D$ for the convex
hull of $E$, $D$ satisfies the requirements for a unit disc. Let us define $\sX$
to be $\real^2$ with the norm $\Norm{\_}$ that has $D$ as its unit disc. Note
that as $D$ is symmetric with respect to the $x$-axis and $y$-axis, $\Norm{\_}$
is invariant under reflection in these axes.

If we let $S$ be the boundary of $D$, i.e., $S$ is the set of unit vectors
under $\Norm{\_}$, then clearly $S$ consists of an infinite family of line
segments, $\pm[\Vv_{i}, \Vv_{i+1}]$, together with the points $\pm\Ve_1$.
The extreme points of $D$ comprise the set $E$, i.e., the $\pm\Vv_i$ and
$\pm\Ve_1$. Any neighbourhood of $\Ve_1$ or $-\Ve_1$ contains infinitely many
extreme points of $D$; moreover, no other point of $S$, or indeed of $\sX$, has
this property.

We now define formulas in the additive language $\LA$ that express
various topological and geometric properties that will let us
define a set of vectors in $\sX$ whose norms comprise the natural
numbers.
\begin{eqnarray*}
\pEP{\Vv} &\IsDef& \all{\Vu\;\Vw} \Norm{\Vu} = \Norm{\Vv} = \Norm{\Vw} \And \Vv = \frac{1}{2}(\Vu + \Vw) \Imp \Vu = \Vv = \Vw\\
\pO{\Vv, \Vw} &\IsDef& \Norm{\Vv - \Vw} = \Norm{\Vv + \Vw}\\
\pACC{\Vv} &\IsDef& \pEP{\Vv} \And  (\all{\epsilon} \epsilon > 0 \Imp \ex{\Vu} \\
    && \quad \Norm{\Vu} = \Norm{\Vv} \And \pEP{\Vu} \And \Vu \not= \Vv \And \Norm{\Vu - \Vv} < \epsilon)\\
\pB{\Vp, \Vq} &\IsDef& \Norm{\Vp} = \Norm{\Vq} = 1 \And \pACC{\Vp} \And \pEP{\Vq} \And \pO{\Vq, \Vp}
\end{eqnarray*}
\noindent
So $\pEP{\Vv}$ holds iff $\Vv$ is an extreme point of the disc $D_{\Norm{\Vv}}$
centred on the origin and of radius $\Norm{\Vv}$ (this is true in $\sX$ iff $\Vv$
lies on the $x$-axis or on one of the lines $l_i$); $\pO{\Vv, \Vw}$ holds iff
$\Vv$ is equidistant from the points $\pm\Vw$; $\pACC{\Vv}$ holds iff $\Vv$ is
a point of accumulation in the set of extreme points of the disc
$D_{\Norm{\Vv}}$ (by the remarks above this is true in $\sX$ iff $\Vv$ lies on
the $x$-axis); and $\pB{\Vp, \Vq}$ holds iff $\Vp$ is an accumulation point in
the set of extreme points of the unit disc $D$ and $\Vq$ is an extreme point of the unit disc equidistant from the points $\pm\Vp$.

If $\Vp = \pm\Ve_1$ and $\Vq = \pm\Ve_2$, we refer to $\Vp$ and $\Vq$ as a
{\em standard basis pair}.  Since the norm on $\sX$ is invariant under reflection in
the $y$-axis, if $\Vv$ lies on the $y$-axis, then $\pO{\Vv, \Ve_1}$ holds in
$\sX$.  The following lemma gives the converse, which means that the predicate
$\pB{\Vp, \Vq}$ characterises the standard basis pairs in $\sX$.

\begin{Lemma}\label{lma:standard-basis}
{\em(i)} $\pO{\Vv, \Ve_1}$ holds in $\sX$ iff $\Vv$ lies on the $y$-axis, whence
{\em(ii)} $\pB{\Vp, \Vq}$ holds  in $\sX$ iff $\Vp = \pm\Ve_1$ and $\Vq = \pm
\Ve_2$.
\end{Lemma}
\Proof
We have already observed that the points $\pm\Ve_1$ are the only accumulation
points in the set of extreme points of the unit disc.  Thus {\em(ii)} follows
from {\em(i)} since {\em(i)} implies that the vectors $\pm\Ve_2$ are the only
unit vectors that are equidistant from $\pm\Ve_1$.  By the
remarks above, we have only to prove that if $\Vv$ is equidistant from
$\pm\Ve_1$, then $\Vv$ lies on the $y$-axis.  Replacing $\Vv$ by $-\Vv$
if necessary, we may assume that $\Vv$ lies in the upper half plane.
So, writing $\Vv = (a, b)$, we may assume $b \ge 0$ and what we have to prove is
that if $\Vv$ is equidistant from $\pm\Ve_1$ then $a = 0$.

So assume that $\Vv$ is equidistant from the points $\pm\Ve_1$, which
means that $\Vv$ lies in the intersection of the sets $F = \Ve_1 + \lambda S$
and $G = -\Ve_1 + \lambda S$, where $\lambda = \Norm{\Vv - \Ve_1} = \Norm{\Vv +
\Ve_1}$.  By the triangle inequality, $2 = \Norm{\Ve_1 + \Ve_1} \le \Norm{\Ve_1
- \Vv} + \Norm{\Ve_1 + \Vv} = 2\lambda$, so $\lambda \ge 1$.  The upper half of
the set $F$ comprises the graph of a function $f : [1 - \lambda, 1 + \lambda]
\rightarrow \real$ and the upper half of $G$ comprises the graph of a function
$g : [-1 - \lambda, -1 + \lambda] \rightarrow \real$.  Since $\Vv = (a, b)$ is
in the upper half-plane by assumption, $a$ must lie in the intersection $[1 -
\lambda, -1 + \lambda]$ of the domains of $f$ and $g$ and we have $b = f(a) =
g(a)$.  As the norm on $\sX$ is invariant under reflection in the $y$-axis, we
have $f(x) = g(-x)$ for $x \in [1 - \lambda, - 1 + \lambda]$, thus $f(0) =
g(0)$ and the point $(0, f(0))$ lies in the intersection of the two graphs.
Now $f$ is strictly increasing on $[1 - \lambda, 1]$ and strictly decreasing
on $[1, 1 + \lambda]$ and $g(x) = f(x+2)$.
So in the (possibly empty) closed interval where $f$ and $g$ are both
defined and $g$ is increasing, we have $g(x) > g(x-2) = f(x)$, while where
$f$ and $g$ are both defined and $f$ is decreasing we have $f(x) >
f(x+2) = g(x)$. Thus $f(a) = g(a)$ implies that $a$ is in the interval
where $f$ is increasing and $g$ is decreasing and there can be
at most one such $a$.
Hence we must have $(a, b) = (0, f(0))$ so that $a = 0$ as required.
\Done

With a few more definitions, we can give a formula of $\LA$ that in $\sX$
characterizes the natural numbers.
\begin{eqnarray*}
\pXAX{\Vv, \Vp, \Vq} &\IsDef& \Vv = \VO \Or (\pACC{\Vv}  \And \Norm{\Vv + \Vp} = \Norm{\Vv} + \Norm{\Vp})\\
\pYAX{\Vv, \Vp, \Vq} &\IsDef& \pO{\Vv, \Vp}  \And \Norm{\Vv + \Vq} = \Norm{\Vv}
+\Norm{\Vq} \\
\pZ{\Vv, \Vp, \Vq} &\IsDef& \pXAX{\Vv, \Vp, \Vq} \And \pEP{\Vv + \Vq} \\
\pNat{x} &\IsDef& \ex{\Vv\;\Vp\;\Vq} x = \Norm{\Vv} \And \pB{\Vp, \Vq} \And \pZ{\Vv, \Vp, \Vq}
\end{eqnarray*}
Thus in $\sX$, if $\Vp$ and $\Vq$ are a standard basis pair: $\pXAX{\Vv, \Vp,
\Vq}$ holds iff $\Vv$ lies on the $x$-axis on the same side as $\Vp$;
$\pYAX{\Vv, \Vp, \Vq}$ holds iff $\Vv$ lies on the $y$-axis on the same side as
$\Vq$; and for $\pZ{\Vv, \Vp, \Vq}$ and $\pNat{x}$ we have:

\begin{Lemma}\label{lma:nat}
{\em (i)} If $\Vp$ and $\Vq$ are a standard basis pair in $\sX$, $\pZ{\Vv, \Vp,
\Vq}$ holds iff $\Vv = x \Vp$ for some $x \in \nat$, whence
{\em (ii)} $\pNat{x}$ holds in $\sX$ iff $x \in \nat$.
\end{Lemma}
\Proof
The right-to-left direction of the claim about $\pZ{\Vv, \Vp, \Vq}$ is easy to
check. So assume $\pZ{\Vv, \Vp, \Vq}$ holds. By Lemma~\ref{lma:standard-basis}, $\Vp =
\pm\Ve_1$ and $\Vq = \pm\Ve_2$. Also $\Vv$ lies on the $x$-axis on the same
side as $\Vp$.  Thus as $\pEP{\Vv + \Vq}$ holds, $\Vv + \Vq = \Vv \pm \Ve_2$ is
the point of intersection of the line $y = \pm 1$ and one of the lines $l_i$
(since it cannot lie on the $x$-axis).  Thus $\Vv$ is indeed a natural number
multiple of $\Vp = \pm\Ve_1$.  The claim about $\pNat{x}$ follows, since
$\pB{\Vp, \Vq}$ implies that $\Norm{\Vp} = 1$.
\Done

The above lemma will give us the undecidability of the theory of any class of
normed spaces that includes the 2-dimensional normed space $\sX$.  The next lemma
lets us transfer information about definability in $\sX$ to definability in
normed spaces and Banach spaces of higher dimensions.

\begin{Lemma}\label{lma:arbcard}
For any $d \in \{2, 3, 4, \ldots\} \cup \{\infty\}$, there is a Banach space
$\sX^d$ with $\Dim{\sX^d} = d$ such that for any formula $\rho(x_1, \ldots, x_k)$ of
$\LN$ with the indicated free variables (all scalar), there is a
formula $\rho^{*}(x_1, \ldots, x_k)$ of $\LN$ with the same free variables such that
under any assignment of real numbers to the $x_i$, $\rho^{*}(x_1, \ldots, x_k)$
holds in $\sX^d$ iff $\rho(x_1, \ldots, x_k)$ holds in $\sX$.  Moreover, $\rho^{*}$ is
additive if $\rho$ is.
\end{Lemma}
\Proof
If $V$ and $W$ are normed spaces, their {\em 1-sum}, $V + W$, is the product
vector space $V \times W$ equipped with the norm defined by $\Norm{(\Vp, \Vq)}
= \Norm{\Vp}_V + \Norm{\Vq}_W$.  $V+W$ has dimension $\Dim{V} + \Dim{W}$ and is
a Banach space iff $V$ and $W$ are both Banach spaces.
The subspaces $V \times 0$ and $0 \times W$ are isomorphic to $V$ and $W$
respectively, and the extreme points of the unit disc in $V + W$ comprise
the points $(\Vv, \VO)$ and $(\VO, \Vw)$ where $\Vv$ and $\Vw$ are extreme
points of the unit discs in $V$ and $W$ respectively.

Let $W$ be the euclidean space $\real^{d-2}$ if $d \not = \infty$ or any
infinite-dimensional Hilbert space, e.g., $l_2$, if $d = \infty$, and let $\sX^d
= \sX + W$.
Now every unit vector in
the Hilbert space $W$ is an extreme point of the unit disc (a counter-example
would give rise to a counter-example in a 2-dimensional subspace and hence a
counter-example in $\real^2$). On the other
hand, the unit disc in $\sX$ has only countably many extreme points.  Moreover a
point of $\sX^d$ lies in $\sX \times 0$ iff it is equidistant from $\pm\Vu$ for
every unit vector $\Vu \in 0 \times W$ (as may be seen by noting that for any
unit vectors $\Vx \in \sX$ and $\Vw \in W$,
there is an isomorphism from $\real^2$ under the 1-norm
to the subspace of $\sX + W$ spanned by $(\Vx, \VO)$ and $(\VO, \Vw)$
that maps $\Ve_1$ to $(\Vx, \VO)$ and $\Ve_2$ to $(\VO, \Vw)$).
It follows that if we define:
\begin{eqnarray*}
\pU{\Vu} &\IsDef& \all{\delta}
              1 > \delta \ge 0
         \Imp \ex{\Vb}
                    \Norm{\Vb - \Vu} = \delta
              \And \Norm{\Vb} = 1
              \And \pEP{\Vb}\\
\pX{\Vv} &\IsDef& \all{\Vu} \pU{\Vu} \Imp \pO{\Vv, \Vu}
\end{eqnarray*}
then $\pU{\Vu}$ holds iff $\Vu$ is a unit vector in $0 \times W$ and $\pX{\Vv}$
holds iff $\Vv$ is in $\sX \times 0$.

Let $\rho^{*}$ be the relativization of $\rho$ to $\pX{\Vv}$, i.e., let $\rho^{*}$ be
obtained from $\rho$ by replacing every subformula of the form  $\ex{\Vv}\phi$ by
$\ex{\Vv}\pX{\Vv} \And \phi$ and every subformula of the form $\all{\Vv}\phi$ by
$\all{\Vv}\pX{\Vv} \Imp \phi$.  Clearly $\rho^{*}$ is in $\LN$ and, as $\pX{\Vv}$ is
additive, $\rho^{*}$ is additive if $\rho$ is.  Since $\pX{\Vv}$ holds iff $\Vv$
belongs to  $\sX \times 0$, under any assignment of real numbers to
the $x_i$, $\rho^{*}(x_1, \ldots, x_k)$ holds in $\sX^{d}$ iff $\rho(x_1, \ldots,
x_k)$ holds in $\sX$.
\Done

We write  $\NS{}$, $\NS{n}$, $\NS{\fin}$ and $\NS{\infty}$
for the theories of normed spaces where the dimension is respectively unconstrained,
constrained to be $n$, constrained to be finite and constrained to
be infinite. We write $\BS{}$, $\BS{n}$ etc. for the theories of Banach
spaces with the corresponding constraints on the dimension.  As
finite-dimensional normed spaces are Banach spaces, $\BS{n} = \NS{n}$ and
$\BS{\fin} = \NS{\fin}$.

\begin{Theorem}\label{thm:L-reduction}
There is a primitive recursive reduction of second-order arithmetic to any of
the theories $\BS{}$, $\BS{\infty}$, $\NS{}$, $\NS{n}$, $\NS{\fin}$, and
$\NS{\infty}$ ($n \ge 2$).
\end{Theorem}
\Proof
What we need to apply Theorem~\ref{thm:mult-undec} is provided by part
{\em(ii)} of Lemma~\ref{lma:nat} using Lemma~\ref{lma:arbcard} in the cases of
$\NS{n}$ for $n > 2$, $\BS{\infty}$ and $\NS{\infty}$.
\Done

Theorem~\ref{thm:L-reduction} is already a satisfyingly sharp result, since as we observed at the
beginning of this section, the
theory of 1-dimensional normed spaces reduces to the theory of the real numbers. But with
a little more work, we can show that scalar multiplication can be defined in
our space $\sX$ in the additive language $\LA$ and so get a reduction of
second-order arithmetic to purely additive normed space theory.
To this end we define some more geometric predicates.
``$\Func{ESD}$'' stands for ``extreme points, same direction''.
\begin{eqnarray*}
\pESD{\Vv, \Vw} &\IsDef& \pEP{\Vv + \Vw} \And \Norm{\Vv + \Vw} = \Norm{\Vv} + \Norm{\Vw}
\end{eqnarray*}
\noindent
I.e., $\pESD{\Vv, \Vw}$ holds iff $\Vv + \Vw$
is an extreme point of the disc $D_{\Norm{\Vv + \Vw}}$ and equality holds
in the triangle inequality for $\Vv$ and $\Vw$. I claim that $\pESD{\Vv, \Vw}$
holds in any normed space iff either $\Vv = \Vw = \VO$ or there is an extreme point $\Vu$
of the unit disc such that $\Vv = x\Vu$ and $\Vw = y\Vu$ for some
non-negative $x$ and $y$.
Thus $\pESD{\Vv, \Vw}$
holds in $\sX$ iff $\Vv$ and $\Vw$ lie on the same side of
the origin on the $x$-axis or on one of the lines $l_i$.
My claim follows easily from the following lemma:

\begin{figure}
\begin{center}
\includegraphics[angle=0,scale=0.9]{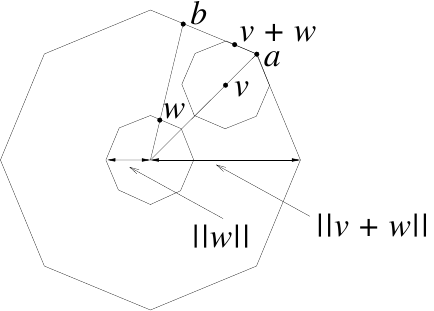}
\caption{%
If $\Norm{\Vv + \Vw} = \Norm{\Vv} + \Norm{\Vw}$ and $\Va \not= \Vb$,
then $\Vv + \Vw \in (\Vb, \Va)$}
\label{fig:SD}
\end{center}
\end{figure}

\begin{Lemma}\label{lma:esd-lemma}
Let $\Vv$ and $\Vw$ be non-zero vectors in a normed space.
If $\Vv + \Vw$ is an extreme point of the disc $D_{\Norm{\Vv+\Vw}}$ of
radius $\Norm{\Vv+\Vw}$ and if $\Norm{\Vv + \Vw} = \Norm{\Vv} + \Norm{\Vw}$,
then $\Vv = \frac{\Norm{\Vv}}{\Norm{\Vw}}\Vw = \frac{\Norm{\Vv}}{\Norm{\Vv+\Vw}}(\Vv+\Vw)$.
\end{Lemma}
\Proof
Under the given hypotheses on $\Vv$ and $\Vw$,
let $\Va = \frac{\Norm{\Vv+\Vw}}{\Norm{\Vv}}\Vv$
and
$\Vb = \frac{\Norm{\Vv+\Vw}}{\Norm{\Vw}}\Vw$
(see Figure~\ref{fig:SD}). As $\Norm{\Vv+\Vw} = \Norm{\Vv}+\Norm{\Vw}$,
we have:
\begin{eqnarray*}
\Vv + \Vw
    &=& \frac{\Norm{\Vv}}{\Norm{\Vv+\Vw}}\Va + \frac{\Norm{\Vw}}{\Norm{\Vv+\Vw}}\Vb\\
    &=& \frac{\Norm{\Vv}}{\Norm{\Vv}+\Norm{\Vw}}\Va + \left(1-\frac{\Norm{\Vv}}{\Norm{\Vv}+\Norm{\Vw}}\right)\Vb .
\end{eqnarray*}
Thus $\Vv + \Vw$ is a proper convex combination of $\Va$ and $\Vb$.
As $\Norm{\Va} = \Norm{\Vb} = \Norm{\Vv+\Vw}$
and $\Vv + \Vw$ is an
extreme point of the disc $D_{\Norm{\Vv+\Vw}}$, we must have $\Va = \Vb$,
i.e., $\frac{\Norm{\Vv+\Vw}}{\Norm{\Vv}}\Vv =
\frac{\Norm{\Vv+\Vw}}{\Norm{\Vw}}\Vw$ implying $\Vv =
\frac{\Norm{\Vv}}{\Norm{\Vw}}\Vw$ and so also $\Vv =
\frac{\Norm{\Vv}}{\Norm{\Vv+\Vw}}(\Vv+\Vw)$.
\Done

\begin{figure}
\begin{center}
\includegraphics[angle=0,scale=0.9]{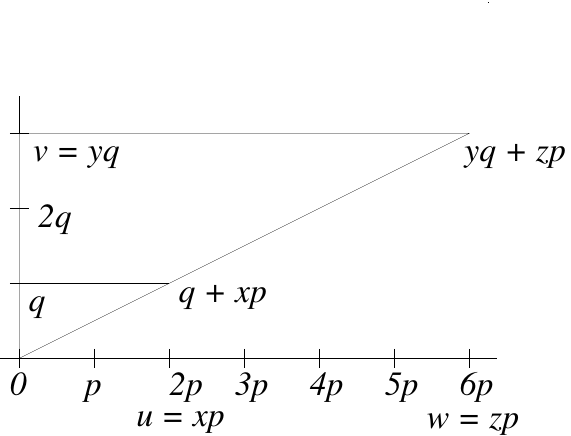}
\caption{$z = xy$}
\label{fig:TIMES}
\end{center}
\end{figure}

We now give the geometric predicate that will allow us to define
multiplication (see Figure~\ref{fig:TIMES}).
\begin{eqnarray*}
\pNTIMES{x, y, z} &\IsDef& \ex{\Vp\;\Vq\;\Vu\;\Vv\;\Vw} x = \Norm{\Vu} \And y = \Norm{\Vv} \And z = \Norm{\Vw} \And {}\\
        && \pB{\Vp, \Vq} \And \pZ{\Vu, \Vp, \Vq} \And \pYAX{\Vv, \Vp, \Vq} \And \pXAX{\Vw, \Vp, \Vq} \And {} \\
        && \pESD{\Vq + \Vu, \Vv + \Vw}
\end{eqnarray*}
\noindent

\begin{Lemma}\label{lma:ntimes}
In $\sX$, $\pNTIMES{x, y, z}$ holds iff $x \in \nat$, $y, z \in \real_{{\ge}0}$ and $z = xy$.
\end{Lemma}
\Proof
By reference to Figure~\ref{fig:TIMES}, it is easy to see that the
right-to-left direction of the lemma holds (put $\Vp = \Ve_1$, $\Vq = \Ve_2$,
$\Vu = x\Vp$, $\Vv = y\Vq$ and
$\Vw = z\Vp$).  Conversely, let $\Vp$, $\Vq$, $\Vu$, $\Vv$ and $\Vw$ be
witnesses to the truth of the existential formula $\pNTIMES{x, y, z}$, so that
$\Norm{\Vu} = x $, $\Norm{\Vv} = y$, $\Norm{\Vw} = z$. Since $\Vp$ and $\Vq$
are a standard basis pair and $\pZ{\Vu, \Vp, \Vq}$, by Lemma~\ref{lma:nat} we
have that $x \in \nat$ and $\Vu = x\Vp$.  Also, since $\pYAX{\Vv, \Vp, \Vq}$
and $\pXAX{\Vw, \Vp, \Vq}$ hold, we have that $\Vv = y \Vq$ and $\Vw = z\Vp$.
Now $\Vq + \Vu = \Vq + x\Vp$ lies on the line $l_x$ passing through the point
$(x, 1)$. Moreover, since $\pESD{\Vq + \Vu, \Vv + \Vw}$ holds, $\Vv + \Vw$ also
lies on $l_x$. But this means that the right-angled triangle $A$ with vertices
$\VO$, $\Vq$ and $\Vq + \Vu = \Vq + x \Vp$ is similar to and parallel to the
triangle $B$ with vertices $\VO$, $\Vv = y\Vq$ and $\Vv + \Vw = y\Vq + z \Vp$.
Hence $z = xy$ completing the proof.
\Done

\begin{Lemma}\label{lma:rtimes}
There is a formula $\pRTIMES{x, y, z} $ in the additive language $\LNA$
which holds in $\sX$ iff $z = xy$.
\end{Lemma}
\Proof
Consider the following formulas of $\LNA$:
\begin{eqnarray*}
\pZTIMES{x, y, z} &\IsDef& \pNTIMES{x, y, z} \Or \pNTIMES{-x, y, -z} \Or {} \\
    &&  \pNTIMES{x, -y, -z} \Or \pNTIMES{-x, -y, z} \\
\pQTIMES{x, y, z} &\IsDef& \ex{m\;n\;t} n \not= 0 \And {} \\
   &&   \pZTIMES{n, x, m}  \And \pZTIMES{m, y, t} \And \pZTIMES{n, z, t} \\
\pRTIMES{x, y, z} &\IsDef& \all{\epsilon} \epsilon > 0 \Imp ( \ex{\delta} \delta > 0 \And {} \\
   &&  \quad ( \all{r\;t} |x - r| < \delta \And \pQTIMES{r, y, t} \Imp |z - t| < \epsilon) ).
\end{eqnarray*}

By Lemma~\ref{lma:ntimes}, in $\sX$, $\pNTIMES{x, y , z}$ defines the graph of
the multiplication function restricted to $\nat \times \real_{{\ge}0}$.
The predicate $\pZTIMES{x, y, z}$ therefore defines the graph of multiplication
restricted to $\zahl \times \real$.
In the formula $\pQTIMES{x, y, z}$, the matrix of the right-hand side of the definition
asserts that $nx = m$ and that $my = t = nz$, so that, when $n \not= 0$, $z =
(m/n)y = xy$, so $\pQTIMES{x,y, z}$ defines the graph of multiplication
restricted to $\rat \times \real$.
By continuity, we have that $\pRTIMES{x, y, z}$ defines the graph of
multiplication without restriction completing the proof of the lemma.
\Done

We write $\NSA{}$, $\NSA{n}$, $\BSA{}$ etc. for the additive subtheories of $\NS{}$, $\NS{n}$, $\BS{}$ etc.
\begin{Theorem}\label{thm:LA-reduction}
There is a primitive recursive reduction of second-order arithmetic to any of
the theories $\BSA{}$, $\BSA{\infty}$, $\NSA{}$, $\NSA{\fin}$, $\NSA{n}$, and
$\NSA{\infty}$ ($n \ge 2$).
\end{Theorem}
\Proof
What we need to apply Theorem~\ref{thm:linear-undec} is provided by part
{\em(ii)} of Lemma~\ref{lma:nat} and Lemma~\ref{lma:rtimes}
using Lemma~\ref{lma:arbcard} in the cases of
$\NSA{n}$ for $n > 2$, $\BSA{\infty}$ and $\NSA{\infty}$.
\Done

\Subsection{The many-one degrees of theories of normed spaces}\label{subsec:many-one-degrees}
Theorems~\ref{thm:L-reduction} and~\ref{thm:LA-reduction} show that the
decision problems for our theories of normed spaces and Banach spaces
are at least as hard as that for the theory of second-order arithmetic.
We now consider the converse problem of reducing the normed space
and Banach space theories to theories of higher-order arithmetic.

As usual, writing $|A|$ for the cardinality of a set $A$, let $\aleph_0 =
\Card{\nat}$ be the first infinite cardinal and let $\C = 2^{\aleph_0} =
\Card{\real}$ be the cardinality of the continuum. If $A$ is any non-empty
finite or countably infinite set, the set $\real^A$ of real-valued functions on
$A$ has cardinality $\C$.  In particular, the set $\Seq{\real}$ of countably
infinite sequences of real numbers has cardinality $\C$.  If $V$ is a vector
space we write $\Card{V}$ for the cardinality of its set of vectors.  Note that
$\Card{V}$ is either 1 or at least $\C$.  We write $\Dimc$ for some fixed vector
space with a basis $B$ of cardinality $\C$, say $B = \{\Vb_x \mid x \in
\real\}$.  Clearly $\Card{\Dimc} \ge \Card{B} = \C$ and, conversely, as any
element of $\Dimc$ is a finite sum $\Sigma_{m=0}^{k}c_m\Vb_{x_m}$ for some
$c_m, x_m \in \real$, $\Card{\Dimc}$ is at most $\Card{\Seq{(\real \times
\real)}} = \C$.  Thus a vector space has cardinality at most $\C$ iff it is
isomorphic to a subspace of $\Dimc$.  The following Skolem-L\"{o}wenheim
theorem thus implies that any satisfiable first-order property of normed spaces
or Banach spaces is satisfiable in a space given by equipping some subspace of
$\Dimc$ with a norm.

\begin{Theorem}\label{thm:skolem}
Let $V$ be a real vector space.
Then $V$ has a subspace $W$ with $|W| \le \C$ that is an elementary substructure of $V$, i.e., a sentence $\phi$ in the language $\LN$ of normed spaces
holds in $V$ iff it holds in $W$.
Moreover, $W$ may be taken to be a Banach space if $V$ is a Banach space.
\end{Theorem}
\Proof
We will construct $W$ using a certain function
$F : \nat \times \Seq{\real} \times \Seq{V} \rightarrow V$.
Let us first show that for any such function there is a subset $W$
of $V$ of cardinality at most $\C$
that is $F$-closed in the sense that $F[\nat \times \Seq{\real} \times \Seq{W}] \subseteq W$.
To see this, define a transfinite sequence of subsets
$W_{\alpha}$ of $V$ as follows, where $\alpha$ is any ordinal and
$\lambda$ is any limit ordinal:
\begin{eqnarray*}
W_0 &=& \{\VO\} \\
W_{\alpha+1} &=& W_{\alpha} \cup F[\nat \times \Seq{\real} \times \Seq{(W_{\alpha})}] \\
W_{\lambda} &=& \bigcup_{\alpha < \lambda} W_{\alpha}.
\end{eqnarray*}
Let $\aleph_1$ be the smallest uncountable cardinal and let $W = W_{\aleph_1}$.
By transfinite induction, one may show that $\Card{W_{\alpha}} \le \C$ for
$\alpha \le \aleph_1$, and so in particular $\Card{W} \le \C$.  Now if $(k, s,
\Vx) \in \nat \times \Seq{\real} \times \Seq{W}$ then I claim $F(k, s, \Vx) \in
W$. For if $\alpha$ is the least ordinal such that $\Vx_m \in W_{\alpha}$ for
all $m \in \nat$, then $\alpha < \aleph_1$ (since $\alpha$ can be written as a
countable union of countable ordinals and hence is countable).  Thus $F(k, s, \Vx)
\in W_{\alpha+1} \subset W$ and $W$
is indeed an $F$-closed subset of $V$ of cardinality at most $\C$.

To define the function $F$, let the formulas of $\LN$ be enumerated as $\psi_1, \psi_2,
\ldots$.  We fix a total ordering on the variables of $\LN$ and choose a vector
variable $\Vv$, and then given $(k, s, \Vx) \in \nat \times \Seq{\real} \times
\Seq{V}$, we define $F(k, s, \Vx)$ as follows:
\begin{enumerate}
\item if $k = 0$ and the $\Vx_m$ converge in $V$ to a limit $\Vp$,
we set $F(k, s, \Vx) = \Vp$;
\item if $k > 0$, consider the formula $\psi \IsDef \ex{\Vv}\psi_k$ and let
$x_0, \ldots, x_m$ and $\Vv_0, \ldots, \Vv_n$ list its free scalar and
vector variables in order. We interpret $x_i$ as $s_i$
and $\Vv_j$ as $\Vx_j$. If $\psi$ is true in $V$ under this interpretation,
then there is a $\Vq$ in $V$ such that $\psi_k$ becomes true if we
extend the interpretation by interpreting $\Vv$ as $\Vq$, we
choose such a $\Vq$ and set $F(k, s, \Vx) = \Vq$;
\item in all other cases, we set $F(k, s, \Vx) = \VO$.
\end{enumerate}

Now let $W \subseteq V$ be an $F$-closed subset of cardinality at most
$\C$ as constructed above. Clause~2 of the definition of $F$ ensures that
the Tarski-Vaught criterion applies so that $W$ is an elementary
substructure of $V$; see, e.g., \citet{hodges-model}.  In particular,
$W$ is a vector space over some subfield of $\real$.  Clause~1 implies
that the 1-dimensional subspaces of this vector space are metrically
complete, so the field of scalars of $W$ may be taken to be $\real$ so that $W$
is a subspace of $V$.  Finally, if $V$ is a Banach space, clause~1 implies
that $W$ is also a Banach space.
\Done

It will simplify our syntactic constructions to extend the language $\A^2$ of
second-order arithmetic as follows: first let $\AR^2$ be the result of adding
to $\A^2$ a sort $\ScaSort$ for the real numbers together with function and
predicate symbols for the operations of the ordered field $\real$ and for the
injection $\iota : \NatSort \rightarrow \ScaSort$ of $\nat$ into $\real$; then
let $\AV^2$ be obtained from $\AR^2$ by adding a sort $\VecSort$ of vectors,
together with function symbols for the vector space operations on $\VecSort$
with scalars in $\ScaSort$ and for a function symbol $\gamma : \VecSort \times
\ScaSort \rightarrow \ScaSort$. The intended interpretation of $\VecSort$ in
$\AV^2$ is the vector space $\Dimc$ with $\gamma$ the operation that maps a
pair $(\Vv, x)$ to the coefficient $c_x$ of the basis element $\Vb_x$ in the
expression of $\Vv$ as a linear combination of elements of the basis $B$. We
choose the symbols so that the language $\LV$ of vector spaces is a sublanguage
of $\AV^2$.

A {\em standard model} of one of the languages $\A^2$, $\AR^2$ or $\AV^2$ is
one in which (up to isomorphism) all the sorts and symbols of the language have
their intended interpretations. In particular, in a standard model, the sort
$\SetSort$ is interpreted as the full powerset $\Pow{\nat}$ of the set of natural
numbers.  Let $T^2_{A}$, resp. $T^2_{AV}$, resp. $T^2_{AR}$,  be the set of all sentences
of $\A^2$, resp. $\AV^2$, resp.  $\AR^2$, that are true in a standard model
(and hence in all standard models).  In the light of the following lemma, to
reduce a decision problem to $T^2_{A}$, i.e., second-order arithmetic, it is
sufficient to reduce it to $T^2_{AV}$.

\begin{Lemma}
There are primitive recursive reductions of $T^2_{AV}$ and $T^2_{AR}$ to the theory
$T^2_{A}$ of true sentences of second-order arithmetic.
\end{Lemma}
\Proof
It is well-known that using suitable encodings, the real numbers may be
constructed, e.g., via Dedekind cuts, as a definitional extension of
second-order arithmetic; see \citet{simpson-subsystems}.  Unwinding the
definitions provides a primitive recursive reduction of $T^2_{AR}$ to $T^2_{A}$.  (The
unwinding process requires occurrences of function symbols first to be unnested
so that they can be replaced by predicates as in the proofs of
Theorem~\ref{thm:mult-undec} and
Lemma~\ref{lma:NS-BS-reduce-to-A2N}).  So it suffices to give a primitive
recursive reduction of $T^2_{AV}$ to $T^2_{AR}$.

Now in $\AR^2$ we can encode the elements of sets such as $\real \times \real$,
$\Seq{\real}$, $\Seq{(\real \times \real)}$ etc. as real numbers.  Given a
vector $\Vv = \Sigma_{m=0}^k c_m\Vb_{x_m} \in \Dimc$, we can arrange for the
$c_m$ to be non-zero and for
the $x_m$ to be listed in strictly increasing order, and then encode $\Vv$ as
the real number that encodes the sequence $s$, with $s_m = (c_m, x_m), 0 \le m
\le k$ and $s_m = (0, 0), m > k$.  Using this encoding we can define the
vector space operations on $\Dimc$ together with the function $\gamma$.
Unwinding these definitions gives the required primitive recursive reduction of
$T^2_{AV}$ to $T^2_{AR}$.
\Done

Now let $\AN^2$ be $\AV^2$ extended with a predicate symbol $\NRM$ of
type $\VecSort \times \ScaSort$.  A {\em standard model} of a sentence
of $\AN^2$ is to be one which extends a standard model of $\AV^2$,
i.e., one in which all the sorts and symbols of $\AV^2$ have their
intended interpretations while the interpretation of $\NRM$ is
arbitrary.

\begin{Lemma}\label{lma:NS-BS-reduce-to-A2N}
There are primitive recursive functions, $\phi \mapsto \phi_N$ and $\phi \mapsto \phi_B$,
which map sentences of the language $\LN$ of normed spaces to sentences of
$\AN^2$, such that the standard models of $\phi_N$ (resp. $\phi_B$) comprise
precisely those standard models in which $\pNRM{\Vv, x}$ defines a norm on a
subspace of $\Dimc$ that provides a model (resp. Banach space model) of $\phi$.
Moreover $\phi$ has a model (resp. Banach space model) iff
$\phi_N$ (resp. $\phi_B$) has a standard model.
\end{Lemma}
\Proof
There is a primitive recursive function mapping $\phi$ to a logically equivalent
sentence $\phi_1$ in which all occurrences of the norm operator are unnested,
i.e., in which the norm operator only appears in atomic formulas of the form
$\Norm{\Vv} = x$ where $\Vv$ and $x$ are variables.
Let $\phi_2$ be the result of replacing each subformula
$\Norm{\Vv} = x$ in $\phi_1$ by $\pNRM{\Vv, x}$. Then $\phi_2$ is a sentence of
$\AN^2$.  Let $\phi_3$ be the relativization of $\phi_2$ to the domain of the
relation defined by $\pNRM{\Vv, x}$, i.e., obtain $\phi_3$ from $\phi_2$ by
replacing subformulas of the form  $\ex{\Vv}\psi$ by
$\ex{\Vv}(\ex{x}\pNRM{\Vv, x}) \And \psi$ and subformulas of the form
$\all{\Vv}\psi$ by $\all{\Vv}(\ex{x}\pNRM{\Vv, x}) \Imp \psi$.

There is a sentence $\pQ_N$ of $\AN^2$ asserting that $\pNRM{\Vv, x}$
defines a relation that {\em(i)} is a partial function, {\em(ii)} has a domain
that is closed under the vector space operations and {\em(iii)} satisfies the
conditions for a norm on the vectors in its domain.  As completeness may
be defined using quantification over countably infinite sequences of vectors,
which is available in $\AN^2$, there is a sentence $\pQ_B$ of $\AN^2$ asserting
that the metric given by $\pNRM{\Vv, x}$ is complete.  We take $\phi_N \IsDef \pQ_N
\And \phi_3$ and $\phi_B \IsDef \pQ_N \And \pQ_B \And \phi_3$.

Now if $\phi$ has a normed space model (resp. Banach space model), then by
Theorem~\ref{thm:skolem}, it has a model that is isomorphic to a subspace $W$
of $\Dimc$ under some norm $\Norm{\_}$.  Extending the standard
interpretation of $\AV^2$ to interpret $\pNRM{\Vv, x}$ as $\Vv \in W \And
\Norm{\Vv} = x$ gives a standard model of $\phi_N$ (resp. $\phi_B$).  Conversely, a
standard model of $\phi_N$ (resp. $\phi_B$) gives a normed space model (resp. Banach
space model) isomorphic to a subspace of $\Dimc$ under the norm defined by the
interpretation of $\pNRM{\Vv, x}$.
\Done

\begin{Theorem}\label{thm:NS-fin-reduces-to-T2}
There are primitive recursive reductions of the theories
$\NS{\fin}$ and $\NS{n}$, $n \in \nat$, to
$T^2_{AV}$ and hence to second-order arithmetic.
\end{Theorem}
\Proof
Let a natural number $n$ and a real number $x$ be given.  In $T^2_{AV}$, we can
define the subspace $\real^n$ of $\Dimc$ spanned by the $\Vb_m, 1 \le m \le n,
m \in \nat$; we can define the subset $\rat^n$ of $\real^n$ comprising the
points with rational coordinates; since $\rat^n$ is countable, we can view $x$ as an encoding of an
arbitrary subset $\VQ^n_x$  of $\rat^n$; using the
coefficient function $\gamma$, we can define the euclidean norm on $\Dimc$ with
respect to the basis $B$.  Thus there is a formula $\Delta(n, x, \Vp, t)$ of
$\AV^2$ that holds in a standard model iff every open disc in $\real^n$ centred on $\Vp$
meets both $t\VQ^n_x$ and its complement $\rat^n \Diff t\VQ^n_x$.
But then if $\VQ^n_x$ is the set $\rat^n \cap D$ of rational points in the unit
disc $D$ of a norm $\Norm{\_}$ on $\real^n$, $\Delta(n, x, \Vp, t)$ holds iff
$\Norm{\Vp} = t$.

We complete the proof for $\NS{\fin}$, the proof for $\NS{n}$ being very
similar.
As a sentence is valid iff its negation is unsatisfiable, it is sufficient to
give a primitive recursive function $\phi \mapsto \phi_F$ from $\LN$ to $\AV^2$ such
that $\phi$ is satisfiable in a finite-dimensional normed space iff $\phi_F$  is
true.
Applying Lemma~\ref{lma:NS-BS-reduce-to-A2N}, we have a sentence $\phi_N$ of
$\AN^2$ that has a standard model iff $\phi$ is satisfiable.  Choose variables
$n$ of sort $\NatSort$ and $x$ of sort $\ScaSort$ that do not appear in $\phi_N$
and let $\psi$ be the result of replacing each occurrence of $\pNRM{\Vp, t}$ in
$\phi_N$ by $\Delta(n, x, \Vp, t)$.  Setting $\phi_F \IsDef \ex{n\; x}\psi$, $\phi_F$ holds in a
standard model of $\AV^2$ iff there are $n \in \nat$ and $x \in \real$ such
that $\VQ^n_x$ is a set of rational points whose closure is the unit disc of a
norm on $\real^n$ and $\phi$ holds under this norm on $\real^n$.  Since any
$n$-dimensional normed space is isomorphic to one given by defining a norm
on $\real^n$, $\phi_F$ is true iff $\phi$ is satisfiable.
\Done

Using the terminology of recursion theory we have the following corollary
concerning degrees of unsolvability; see, e.g., \citet{rogers-recursive-functions} for definitions.

\begin{Corollary}\label{thm:NS-fin--many-one}
The theories $\NS{\fin} = \BS{\fin}$, $\NSA{\fin} = \BSA{\fin}$ and $\NS{n} =
\BS{n}$, $\NSA{n} = \BSA{n}$, $n \ge 2$, all have the same many-one degree as the theory $T^2_{A}$ of
second-order arithmetic.
\end{Corollary}
\Proof
This is immediate from Theorems~\ref{thm:LA-reduction} and~\ref{thm:NS-fin-reduces-to-T2}.
\Done

Now let $\A^3$ be the language of third-order arithmetic.  This is $\A^2$
extended with an additional sort $\SetSort_2$ called ``type 2'' whose intended
interpretation is $\Pow{\Pow{\nat}}$.  $\A^3$ has a predicate symbol $\in$
of type $\SetSort \times \SetSort_2$ to denote the membership relation and a
supply of type 2 variables $u = u_1, u_2, \ldots$, but we shall only need the
first of these.  A sentence of $\A^3$ is said to be $\Sigma^2_1$ (resp.
$\Pi^2_1$) if it has the form $\ex{u}\psi(u)$ (resp. $\all{u}\psi(u)$)
where $\psi(u)$ contains no quantifiers over type 2 variables.

\begin{Theorem}\label{thm:NS-BS-reduce-to-Pi21}
There are primitive recursive reductions of each of the theories
$\NS{}$, $\BS{}$, $\NS{\infty}$ and $\BS{\infty}$
to the set of true $\Pi^2_1$ sentences.
\end{Theorem}
\Proof
As with $\A^2$ we are free to work in a definitional extension $\AV^3$ of $\A^3$
that includes the language $\LV$ of vector spaces (with $\Dimc$ as the
intended interpretation of the vector sort).
Let $a$, $\Vv$ and $x$ be variables
of sort $\SetSort$, $\VecSort$ and $\ScaSort$ respectively.
There is a formula $\pU{a, \Vv, x}$ of $\AV^2 \subseteq \AV^3$ with the indicated
free variables that in a standard model of $\A^2$ defines the graph of a
bijection mapping $a \in \Pow{\nat}$ to $(\Vv, x) \in \Dimc \times \real$.
This gives an encoding of all relations between $\Dimc$ and $\real$, i.e.,
all subsets of $\Dimc \times \real$, as type 2 sets.

To complete the proof, let us first consider $\NS{}$.
As a sentence is valid iff its negation is unsatisfiable, it suffices to give a
primitive recursive function $\phi \mapsto \phi_{1}$ from
the language $\LN$ of normed spaces to the set of $\Sigma^2_1$ sentences such that
$\phi$ is satisfiable iff $\phi_{1}$ is true.
Given a sentence $\phi$ in the language of normed spaces, apply
Lemma~\ref{lma:NS-BS-reduce-to-A2N}, to give a sentence $\phi_N$ of $\AN^2$ that has a
standard model iff $\phi$ is satisfiable.  Let $\psi(u)$ be obtained from $\phi_N$ by
replacing all instances of $\pNRM{\Vv, x}$ by $\ex{a} a \in u \And \pU{a,
\Vv, x}$ and let $\phi_{1}$ be $\ex{u}\psi(u)$. Then $\phi_{1}$ is a $\Sigma^2_1$
formula that is true iff $\phi_N$ has a standard model.
So $\phi_{1}$ is true iff $\phi$ is satisfiable.

For $\BS{}$, we use a primitive recursive function $\phi \mapsto \phi_{2}$ from
$\LN$ to the set of $\Sigma^2_1$ sentences such that $\phi$ is satisfiable in a
Banach space iff $\phi_{2}$ is true.  The construction of $\phi_{2}$ is identical
to that of $\phi_{1}$ except that we use the sentence $\phi_B$ from
Lemma~\ref{lma:NS-BS-reduce-to-A2N} rather than $\phi_N$.

Finally, for $\NS{\infty}$ and $\BS{\infty}$, there is a formula
$\pI{u}$ of $A^3$ with no type 2 quantifiers which holds iff $u$ encodes
a relation between $\Dimc$ and $\real$ whose domain is an infinite-dimensional
subspace of $\Dimc$. Relativization of $\phi_{1}$ and $\phi_{2}$ to $\pI{u}$
gives the reductions required to complete the proof.
\Done

\begin{figure}
\begin{center}
\includegraphics[angle=0,scale=1.6]{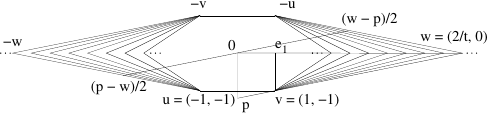}
\caption{The unit discs of the spaces $\sS_t$}
\label{fig:steeples}
\end{center}
\end{figure}

We complete our study of the degrees of unsolvability of the normed space
and Banach space theories by exhibiting a primitive recursive reduction
of the set of true $\Pi^2_1$ sentences to the theories $\NSA{}$ and $\BSA{}$.
To do this we need Banach spaces in which an arbitrary
subset of the open interval $(0, 1)$ can be defined in a uniform way.
We begin by considering the special case of a singleton set.
So let $t \in (0, 1)$ be given
and define points of $\real^2$ by  $\Vu = (-1, -1)$, $\Vv = (1, -1)$
and $\Vw = (\frac{2}{t}, 0)$. Let $\sS_t$ be $\real^2$ equipped with the norm
$\Norm{\_}_t$ whose unit circle comprises the hexagon with vertices
$\pm\Vu$, $\pm\Vv$ and $\pm\Vw$ (see Figure~\ref{fig:steeples}).
One finds using the ordinary euclidean norm $\Norm{\_}_e$ that $\Norm{\Vv -
\Vu}_t = \frac{\Norm{\Vv - \Vu}_e}{\Norm{\Vw}_e} = \frac{2}{2/t} = t$.
Let the line through $\Vv$ and $\Vw$ meet the $y$-axis at the point $\Vp$.
Then the line segment $[\Vp, \Vw]$ is a translate of the line segment
$[(\Vp - \Vw)/2, (\Vw - \Vp)/2]$ which is a diameter
of the unit
disc in $\sS_t$. So $\Norm{\Vw-\Vp}_t = 2$ and $\Norm{\Vw-\Vv}_t = 2 - \Norm{\Vv-\Vp}_t$. But the triangles $\VO\Vp\Vw$ and $\Ve_1\Vv\Vw$ are similar
and so $\Norm{\Vv-\Vp}_t = \Norm{\Vw-\Vp}_t(\frac{\Norm{\Vv-\Vp}_e}{\Norm{\Vw-\Vp}_e}) = 2(\frac{\Norm{\Ve_1}_e}{\Norm{\Vw}_e}) = t$ whence
$\Norm{\Vw - \Vv}_t = 2 - \Norm{\Vv - \Vp}_t = 2 - t$.
By symmetry, each edge of the hexagon that comprises the
unit circle in $\sS_t$ has length $t$ or $2 - t$ in the $\sS_t$ norm.

Let us say that two vectors $\Vp$ and $\Vq$ in a normed space $V$ are {\em
adjacent} if $\Vp$ and $\Vq$ are distinct extreme points of the set
$S_{\Norm{\Vp}}$ of vectors of length $\Norm{\Vp}$ and
$\Norm{\frac{1}{2}(\Vp +\Vq)} = \Norm{\Vp}$.  This implies that the line segment
$[\Vp, \Vq]$ is the intersection of some affine line with the set
$S_{\Norm{\Vp}}$.  If $\Vp$ and $\Vq$ are adjacent unit vectors then $\Norm{\Vp
- \Vq} \le 2$ with equality iff $\Vp$ and $-\Vq$ are also adjacent unit
vectors, in which case the linear transformation that maps $\Ve_1$ to $\Vp$
and $\Ve_2$ to $\Vq$ defines an isomorphism between $\real^2$ under
the 1-norm and the subspace of $V$ spanned by $\Vp$ and $\Vq$.
Now consider the following formulas in $\LNA$, the first of which
formalises the notion of adjacency.
\begin{eqnarray*}
\pADJ{\Vp, \Vq} &\IsDef& \pEP{\Vp} \And \pEP{\Vq} \And \Vp \not= \Vq \And \Norm{\Vp} = \Norm{\Vq} = \Norm{(\Vp + \Vq)/2}\\
\pH{\Vu, \Vv, \Vw} &\IsDef&
\begin{array}[t]{@{}l}
\pADJ{\Vu, \Vv} \And \pADJ{\Vv, \Vw} \And \pADJ{\Vw, -\Vu} \And {} \\
\Norm{\Vv - \Vu} < 2\Norm{\Vv} \And \Norm{\Vw - \Vv} <
2\Norm{\Vv} \And \Norm{\Vw + \Vu} < 2\Norm{\Vv}\\
\end{array} \\
\pT{x} &\IsDef& \ex{\Vu\; \Vv\; \Vw} \Norm{\Vu} = 1 \And \pH{\Vu, \Vv, \Vw} \And x = \Norm{\Vv - \Vu} < 1
\end{eqnarray*}

Clearly $\pADJ{\Vu, \Vv}$, $\pADJ{\Vv, \Vu}$ and
$\pADJ{-\Vu, -\Vv}$ are all equivalent and so
in any normed space, $\pH{\Vu, \Vv, \Vw}$ implies that the vectors $\Vu,
\Vv, \Vw, -\Vw, -\Vu, -\Vw$ are the vertices of a hexagon
inscribed in the set $S_{\Norm{\Vu}}$ of vectors of length $\Norm{\Vu}$.
$\pH{\Vu, \Vv, \Vw}$ also includes a condition on the length of the edges of
this hexagon that will presently help us pick out elements of $\sS_t$ when it is
embedded in a larger space.  Now in $\sS_t$, if $\Vu$ and $\Vv$ are unit vectors
and $\Norm{\Vv - \Vu} < 1$, $\pH{\Vu, \Vv, \Vw}$ can only hold if $\Vu$ and
$\Vv$ are the end-points of one of the edges of $S_1$ whose length in the $\sS_t$
norm is $t$, so in $\sS_t$, $\pT{x}$ defines the singleton set $\{t\}$.

We now need a generalisation of the 1-sum construction that we used in the
proof of Lemma~\ref{lma:arbcard}.  Let $V_i$, $i \in I$, be an arbitrary family
of normed spaces and write $\Norm{\_}_i$ for the norm on $V_i$.  If $\Vf$ is a
member of $\Pi_{i \in I} V_i$ and if $J$ is a finite subset of $I$, let $n(\Vf,
J) = \Sigma_{j \in J}\Norm{\Vf_j}_j$. The {\em 1-sum} $\Sigma_{i\in I}V_i$
comprises those $\Vf$ for which $n(\Vf, J)$ is bounded as $J$ ranges over all
finite subsets of $I$. We define $\Norm{\Vf}$ to be the supremum of the $n(\Vf,
J)$.  As is easily verified, $\Sigma_{i\in I}V_i$ is a normed space and is a
Banach space iff the $V_i$ are all Banach spaces.  There is a natural
isomorphism between the summand $V_i$ and the subspace of $\Sigma_{i\in I}V_i$
comprising those $\Vf$ such that $\Vf_j = \VO$ whenever $j \not= i$ and we may
identify $V_i$ with that subspace.  Under this identification, the extreme
points of the unit disc in $\Sigma_{i\in I}V_i$ comprise the union of the
extreme points of the unit discs of the $V_i$.  If $\Vp \in V_i$ and $\Vq \in
V_j$ are unit vectors and $i \not= j$, then in the 1-sum, $\Norm{\Vp - \Vq} =
2$.

If $T$ is any subset of the interval $(0, 1)$, let $\sS_T = \Sigma_{t \in
T}\sS_t$. Then $\sS_T$ is the 1-sum of Banach spaces and hence is itself a
Banach space.  I claim that the formula $\pT{x}$ that defines $t$ in the space
$\sS_t$ defines $T$ in the 1-sum $\sS_T + V$ where $V$ is any normed space
whose unit circle contains no hexagons.  For, assume that $\pT{t}$ holds for
some $t$.  Then there are extreme points $\Vu$, $\Vv$ and $\Vw$ of the unit
disc in $\sS_T + V$ such that $t = \Norm{\Vv - \Vu} < 1$ and $\Norm{\Vw - \Vv}
< 2\Norm{\Vv} = 2$.  Now as $\Norm{\Vv - \Vu}, \Norm{\Vw - \Vv} < 2$, $\Vu$,
$\Vv$ and $\Vw$ are either all in $\sS_T$ or all in $V$ (viewed as subspaces of
$\sS_T + V$), and as they lie on a hexagon contained in the unit circle they
must all lie in $\sS_T$.  But then $\Vu$, $\Vv$ and $\Vw$ must belong to the
same summand of $\sS_T$ and that summand must be $\sS_t$, so $t \in T$.
Conversely, if $t \in T$, then $\pT{t}$ holds in $\sS_t$, and then, as the
extreme points in the unit disc of $\sS_t$ are a subset of those of $\sS_T +
V$, $\pT{t}$ must hold in $\sS_T + V$.

\begin{Theorem}\label{thm:x-t-defines-t}
There is a formula $\pT{x}$ in $\LN$ with the indicated scalar free variable
such that {\em(i)} in any normed space $\pT{x}$ defines a subset of the interval
$(0, 1)$ and {\em(ii)} for any set $T \subseteq (0, 1)$,
and any normed space $V$ whose unit circle contains no hexagons,
$\pT{x}$ defines $T$ in the 1-sum $\sS_T + V$.
\end{Theorem}
\Proof
Taking $\pT{x}$ as defined above, we have already proved {\em(ii)},
while {\em(i)} is immediate from the definition of $\pT{x}$.
\Done

\begin{Theorem}\label{thm:Pi21-reduces-to-NS-BS}
There are primitive recursive reductions of the set of all true $\Pi^2_1$
sentences to each of the theories $\NSA{}$, $\BSA{}$, $\NSA{\infty}$ and $\BSA{\infty}$.
\end{Theorem}
\Proof
It suffices to produce a primitive recursive function $\phi \mapsto \phi_A$ from the
set of $\Sigma^2_1$ sentences to $\LNA$ such that {\em (i)} $\phi$ is
true iff $\phi_A$ has a normed space model and {\em(ii)} whenever $\phi_A$ has
a normed space model it also has an infinite-dimensional Banach space model.
So let $\phi$ be a $\Sigma^2_1$ sentence $\ex{u}\psi(u)$.

We work in the 1-sum $\sX + \sS_T$ where
$\sX$ is the 2-dimensional normed space defined at the beginning of
Section~\ref{subsec:undec-reducing}
and illustrated in Figure~\ref{fig:infinigon}
and $\sS_T$ is as above for some $T \subseteq (0, 1)$.
Consider the following formulas of $\LNA$:
\begin{eqnarray*}
\pEPX{\Vv} &\IsDef& \pEP{\Vv} \And \Not\ex{\Vu\; \Vw} \pH{\Vu, \Vv, \Vw}\\
\pX{\Vv} &\IsDef& \ex{\Vu\; \Vw} \pEPX{\Vu} \And \pEPX{\Vw} \And \Vv = \Vu + \Vw.
\end{eqnarray*}

In $\sX + \sS_T$, $\pEPX{\Vv}$ holds iff $\Vv$ is an extreme point of the
disc of radius $\Norm{\Vv}$ in the summand $\sX$ and so $\pX{\Vv}$
holds iff $\Vv$ is a sum of such extreme points, which is true iff
$\Vv \in \sX$.
If, as in the proof of Lemma~\ref{lma:arbcard}, we relativize the earlier
definitions of the formulas $\pNat{x}$ and $\pRTIMES{x, y , z}$ to $\pX{\Vv}$
then the resulting formulas will define the set of natural numbers and the
graph of the multiplication function in $\sX + \sS_T$ just as in
Theorem~\ref{thm:LA-reduction}. As in Section~\ref{sec:prelim}, there are
sentences $\pPeano$ and $\pMult$ of $\LNA$ asserting that the relativized
versions of $\pNat{x}$ and $\pRTIMES{x, y, z}$ do indeed define the natural numbers and real multiplication respectively.

Let $\psi_1$ be obtained from $\psi(u)$ as follows:
first, replace each subformula of the form $x^{\SetSort} \in u$ by
$\pT{\frac{1}{3}(x^{\ScaSort} + 1)}$ and translate all other formulas as in the
reduction of second-order arithmetic of Theorem~\ref{thm:mult-undec},
using $\pD(n, x)$ to represent sets of natural numbers as real numbers,
using $\pS(x)$ to single out canonical representatives and
using the relativized $\pNat{x}$ as the predicate for the natural numbers;
then, as in the proof of Theorem~\ref{thm:linear-undec}, eliminate
multiplication using the relativized $\pRTIMES{x, y, z}$.
Now let $\phi_A \IsDef \psi_1 \And \pPeano \And \pMult$.
By construction $\phi_A$ contains no terms of the form $a\Vv$, so $\phi_A$
is indeed in $\LNA$.

We may now check conditions {\em(i)} and {\em(ii)}. First, assume $\phi_A$ has a
model, and in that model let $U = \{S \mid \pT{\frac{1}{3}(\sharp S + 1)}\}$
where  $\sharp$ is the injection of $\Pow{\nat}$ into the interval $[0, 3/2]$
defined in Section~\ref{sec:prelim}.  Then as $\pPeano$ and $\pMult$ hold,
$\psi(u)$ must hold in the standard model when $u$ is interpreted as $U$, so $\phi$,
i.e., $\ex{u}\psi(u)$ is true.  Conversely, if $\phi$ is true, so that $\psi(u)$ holds
when $u$ is interpreted as $U$ say, then if we put $T = \{\frac{1}{3}(\sharp S
+ 1) \mid S \in U\} \cup (0, 1/3)$, $\phi_A$ is satisfied in the normed space $\sX
+ \sS_T$ (since if $\frac{1}{3}(x + 1) \in (0, 1/3)$ then $x < 0$ and $\pS(x)$
is false).  Now $\sX + \sS_T$ is a Banach space and is infinite-dimensional so
if $\phi_A$ has a model it has an infinite-dimensional Banach space model.
\Done

\begin{Corollary}\label{cor:NS-BS-many-one}
The theories $\NSA{}$, $\BSA{}$, $\NSA{\infty}$, $\BSA{\infty}$,
$\NS{}$, $\BS{}$, $\NS{\infty}$ and $\BS{\infty}$
all have the same many-one degree as the set of all true $\Pi^2_1$ sentences,
\end{Corollary}
\Proof
This follows immediately from Theorems~\ref{thm:Pi21-reduces-to-NS-BS}
and~\ref{thm:NS-BS-reduce-to-Pi21}.
\Done

As a final remark on degrees of unsolvability, close analogues of
the above results on normed spaces and Banach spaces hold for metric spaces:
there is a Skolem-L\"{o}wenheim theorem stating that any (complete) metric
space has an elementarily equivalent (complete) subspace of cardinality at most
$\C$; the theory of countable metric spaces is many-one equivalent to
second-order arithmetic; and the theory of arbitrary metric spaces is many-one
equivalent to the set of true $\Pi^2_1$ sentences.
(For the analogue of the space $\sX + \sS_T$ in the proof
of Theorem~\ref{thm:Pi21-reduces-to-NS-BS}, choose $\Vv \in \real^2$
such that $d(\Vv, \sG) \ge 2$ where $\sG$ is the
space of Theorem~\ref{thm:ms-additive-undecidable}
and, for $T \subseteq (0, 1)$, let
$\sH_T \IsDef \{\Vv\} \cup \{ \Vu \in \real^2 \ST d(\Vu, \Vv) - 1 \in T\}$.
Then, in place of $\sX + \sS_T$, use $\sG \cup \sH_T$ and design the
various formulas needed using the fact that $\Vv$ is the only
isolated point.)

\Section{Quantifier elimination for theories of inner product spaces}\label{sec:qelim}

The main idea of this section is that in the first-order theory of inner
product spaces over $\real$ it should take at most $k$ degrees of freedom to
decide the validity of a formula with $k$ vector variables. The key result
implies that if a formula $\phi$ has free vector variables $\Vv_1, \ldots \Vv_m$
and has $k$ vector variables in all, then in all dimensions ${} \ge k$,
$\phi$ is equivalent to a system of constraints on the inner products
$\inner{\Vv_i, \Vv_j}$. The proof is via a process that eliminates
vector quantifiers in favour of blocks of scalar quantifiers.
It follows that to decide a sentence with $k$ vector variables we need only
decide it in $\real^n$ for $n = 0, 1, \ldots, k$ and that is easy after a simple
syntactic transformation given a decision procedure for formulas that do not
involve vectors, i.e., for the language of a real closed field.

In the paper that our title echoes, \citet{tarski-decision} gave the first
quantifier elimination procedure for a real closed field and hence a decision
procedure of the kind that we need.  Apparently the first actual computer
implementation of an algorithm for this problem was by \citet{collins}. A
relatively simple procedure due to Cohen and H\"ormander
\citep{hormander-pdo2,garding-history,bochnak-realag} has been implemented by
several people including one of the present authors.  Collins's method of
cylindrical algebraic decomposition has complexity exponential in the number of
bound variables. The best known algorithms are exponential in the number of
quantifier alternations (see~\citet{basu-algorithms}), but work on implementation
of these algorithms is in its early stages. Since our syntactic transformations
replace vector quantifiers by blocks of scalar quantifiers, these recent
improvements are significant for the complexity of our decision procedure.

We write $\IP{}$, resp., $\IP{\fin}$, resp., $\IP{\infty}$ for the theories of
real inner product spaces where the dimension is unconstrained, resp.,
constrained to be finite, resp., constrained to be infinite, and $\HS{}$,
$\HS{\fin}$ and $\HS{\infty}$ for the theories of Hilbert spaces with the
corresponding constraints on the dimension.  By the well-known fact that finite
dimensional inner product spaces are complete, $\HS{\fin} = \IP{\fin}$.  We
will show that all of these theories are decidable and that $\IP{} = \IP{\fin}
= \HS{} = \HS{\fin}$ and that $\IP{\infty} = \HS{\infty}$.

Let us agree on some terminology and notation.
Given a formula $\phi$ of $\LI$, let $v(\phi)$ and $s(\phi)$
denote the sets of free vector variables and free scalar variables of $\phi$
respectively.
If $\Tuple{\Vv} = (\Vv_1, \ldots, \Vv_m)$ is a sequence of vector variables
and $\Tuple{x}= (x_1, \ldots, x_n)$ is a sequence of scalar variables, let us
write $\phi(\Tuple{\Vv}, \Tuple{x})$ to indicate that
$v(\phi) \subseteq \{\Vv_1, \ldots, \Vv_m\}$ and
$s(\phi) \subseteq \{x_1, \ldots, x_n\}$.
Let $V$ be an inner product space.
If $\phi$ is a sentence of $\LI$, we
write $V \models \phi$ to indicate that $\phi$ holds in $V$.
More generally, if $\phi(\Tuple{\Vv}, \Tuple{x})$ is any formula
in $\LI$, and if $\Tuple{\Vp} \in V^m$
and $\Tuple{c} \in \real^n$, we write $V \models \phi(\Tuple{\Vp}, \Tuple{c})$
to indicate that $\phi$ holds in $V$ if each $\Vv_i$ is interpreted as $\Vp_i$
and each $x_j$ is interpreted as $c_j$.
Note that if the formula $\phi$ contains no constants or variables of vector
sort, then $\phi$ is a formula in the first-order language of an ordered field
and, for any $V$, $V \models \phi(\emptyset, \Tuple{c})$ iff
$\phi(\Tuple{c})$ holds in the ordered field $\real$.

For $k \in \nat$, let us say that formulas $\phi_1(\Tuple{\Vv}, \Tuple{x})$ and $\phi_2(\Tuple{\Vv}, \Tuple{x})$ with
the same free variables are {\em $k$-equivalent} iff for every inner product space
$V$ of dimension at least $k$, and every $\Tuple{\Vp} \in V^{\Size{\Tuple{\Vv}}}$ and
every $\Tuple{c} \in \real^{\Size{\Tuple{x}}}$,
$V \models \phi_1(\Tuple{\Vp}, \Tuple{c})$ iff $V \models \phi_2(\Tuple{\Vp}, \Tuple{c})$,
i.e., $\phi_1$ and $\phi_2$ are equivalent in the theory of all spaces of
dimension at least $k$. So, for example, the sentences $\ex{\Vv\;\Vw} \all{x}
\Vv \not= x\Vw \And \Vw \not= x\Vv$ and $\VO = \VO$ are 2-equivalent, but not
1-equivalent.
Providing they have the same free variables, logically equivalent formulas are
$k$-equivalent for any $k$.

If $V$ is an inner product space and $\Tuple{\Vp} \in V^m$, recall that the {\em Gram matrix}
of $\Tuple{\Vp}$ is the positive semidefinite symmetric $m \times m$ matrix $G = G(\Tuple{\Vp})$
with $G_{ij} = \inner{\Vp_i, \Vp_j}$.
If $\Tuple{\Vv} = (\Vv_1, \ldots, \Vv_m)$ is a sequence of vector variables, let us write $\pUG(\Tuple{\Vv})$ for
the sequence of terms of $\LI$ defined inductively by:
$$
\begin{array}{rcl}
\pUG(\Vv_1) &\IsDef& (\inner{\Vv_1, \Vv_1}) \\
\pUG(\Vv_1, \ldots, \Vv_m) &\IsDef& \pUG(\Vv_1, \ldots, \Vv_{m-1}) \frown 
 (\inner{\Vv_1, \Vv_m}, \ldots,  \inner{\Vv_m, \Vv_m})
\end{array}
$$
where $\frown$ denotes concatenation.
Thus $\pUG(\Tuple{\Vv})$ enumerates the upper triangle
of the formal Gram matrix of $\Tuple{\Vv}$ by column.  Let us say
that a formula $\phi(\Tuple{\Vv}, \Tuple{x})$ is {\em special} iff it
has the form $\psi(\pUG(\Tuple{\Vv}), \Tuple{x})$ where
$\psi(w_1, \ldots, w_{m(m+1)/2}, \Tuple{x})$ contains no variables or
constants of vector sort, i.e., $\psi$ is a formula in the language of
an ordered field.  Note that if $\phi$ and $\psi$ are special, then so
are $\Not \phi$, $\phi \circ \psi$ and $\ex{x}\phi$, where $\circ$ is
any binary propositional connective, e.g., ${\And}$, ${\Or}$, ${\Imp}$
or ${\Iff}$.

Our main theorem will inductively transform a formula of $\LI$ containing
$k$ vector variables into a
$k$-equivalent special formula. The following two lemmas give the two main
ingredients of the proof.

\begin{Lemma}\label{lma:standard-form}
There is a primitive recursive function, $\phi \mapsto \PreSpecial{\phi}$, such
that, for any formula $\phi$ of $\LI$, $\PreSpecial{\phi}$ is equivalent
to $\phi$ in the theory of real inner product spaces and the only terms of
vector sort in $\PreSpecial{\phi}$ are variables occurring as operands of the
inner product operator $\inner{\_,\_}$. Moreover $\PreSpecial{\phi}$ is
quantifier-free if $\phi$ is.
\end{Lemma}
\Proof
There is a primitive recursive $p$ such that $p(\phi)$ results from $\phi$ by
replacing each vector equation $\Va = \Vb$ by the equivalent scalar equation $
\inner{\Va - \Vb, \Va - \Vb} = 0$.  I claim that there exists a primitive recursive
$q$ such that $q(\psi)$ results from $\psi$ by repeatedly applying the
following equations as left-to-right rewrite rules until no redexes remain.
\[
\begin{array}{rcl@{\quad}rcl}
    \inner{\Va, \VO} &=& 0
 &  \inner{-\Va, \Vb} &=& -\inner{\Va, \Vb}
\\
    \inner{\VO, \Va} &=& 0
 &  \inner{\Va, -\Vb} &=& -\inner{\Va, \Vb}
\\
    \inner{t\Va, \Vb} &=& t \inner{\Va, \Vb}
 &  \inner{\Va, \Vb + \Vc} &=& \inner{\Va, \Vb} + \inner{\Va, \Vc}
\\
    \inner{\Va, t\Vb} &=& t \inner{\Va, \Vb}
 &  \inner{\Va + \Vb, \Vc} &=& \inner{\Va, \Vc} + \inner{\Vb, \Vc}.
\end{array}
\]
Thus, if $\psi$ contains no vector equations, $q(\psi)$ will contain no terms
of vector sort other than variables occurring as operands of $\inner{\_, \_}$.
So given $q$, we may take $\PreSpecial{\phi} = q(p(\phi))$ to complete the
proof.  For the existence of $q$ one can either apply a general result of
\citet{hofbauer-mpo} or use the following construction.  Let the {\em weight}
of a redex be the total number of constant and function symbols it contains.
There is a primitive recursive $f$ such that, if $\psi$ contains a redex, then
$f(\psi)$ results from $\psi$ by applying one rule to a redex of maximal
weight.  Let ${n(\psi)}$ be the number of redexes of maximal weight in $\psi$
and let $g(\psi) = f^{n(\psi)}(\psi)$.  Now let ${k(\psi)}$ be $0$ if $\psi$
has no redexes and be the maximal weight of a redex in $\psi$ otherwise and let
$q(\psi) = g^{k(\psi)}(\psi)$.  Then $q$ is primitive recursive and $q(\psi)$
results from $\psi$ by applying rewrite rules until no redexes remain.
\Done

\begin{Lemma}\label{lma:gram-matrix}
Let $M$ be a symmetric $m \times m$ matrix with real coefficients, let $V$
be an inner product space of dimension at least $m$ and let
$\Vp_1, \ldots, \Vp_{m-1} \in V$ be such that
$$
G(\Vp_1, \ldots, \Vp_{m-1}) = (M_{ij})_{1 \le i, j < m}
$$
where $G(\Vp_1, \ldots, \Vp_{m-1})$ is the Gram matrix of the $\Vp_i$.
The following are equivalent:
\begin{itemize}
\item[{\em(i)}] there exists $\Vp_m \in V$ such that
$M = G(\Vp_1, \ldots, \Vp_{m-1}, \Vp_m)$;
\item[{\em(ii)}] there exist $b_1, \ldots, b_m \in \real$
such that
$$
\begin{array}{rrcll}
\mbox{\Lab{a}} & M_{im} &=& \sum_{j=1}^{m-1} b_jM_{ij} \quad & \mbox{for $1 \le i < m$}\\
\mbox{\Lab{b}} & M_{mm} &=& \sum_{i=1}^{m-1}\sum_{j=1}^{m-1} b_ib_jM_{ij} + b_m^2. &
\end{array}
$$
\end{itemize}
\end{Lemma}
\Proof
For both parts of the proof, let $W$ be the subspace spanned by the $\Vp_i$
with $1 \le i < m$ and note that $W$ is a proper subspace of $V$ since $\Dim{V}
\ge m$.

$\mbox{\em{(i)}} \Imp \mbox{\em{(ii)}}$:
Given $\Vp_m \in V$ such that 
$G(\Vp_1, \ldots, \Vp_{m-1}, \Vp_m) = M$, there is a unit vector $\Vc$
orthogonal to $W$ such that $\Vp_m$ lies in the subspace spanned by $W$ and $\Vc$.
Then we can write $\Vp_m = \sum_{i=1}^{m-1}b_i\Vp_i + b_m\Vc$ for some $b_i \in \real$.
{\Lab{a}} and {\Lab{b}} follow for this choice of the
$b_i$ using the expression for $\Vp_m$
to expand the inner products $M_{im} = \inner{\Vp_i, \Vp_j}$, $1 \le i \le m$.

$\mbox{\em{(i)}} \Pmi \mbox{\em{(ii)}}$:
Given $b_1, \ldots, b_m \in \real$ satisfying
{\em(a)} and {\em(b)}, choose a unit vector $\Vc$ orthogonal to $W$ and
let $\Vp_m = \sum_{i=1}^{m-1}b_i\Vp_i + b_m\Vc$.
We must show that the equation $M_{ij} = \inner{\Vp_i, \Vp_j}$ holds
for $1 \le i, j \le m$.
But this is so by assumption when $1 \le i, j < m$,
by {\Lab{a}} when $1 \le i < j = m$,
by symmetry and {\Lab{a}} when $1 \le j < i = m$
and by {\Lab{b}} when $ i = j = m$.
\Done

We now give the main theorem of this section. In this theorem, we need to count the number of vector variables in a formula. This is to be done in a very frugal way, by ignoring
variable binding and simply counting the number of distinct variable names
that appear labelled with the vector sort: so that, for example, $(\all{\Vv\;\Vw} \Vv + \Vw = \VO) \Imp
(\all{\Vv} \Vv = \Vw)$ contains just two variables, $\Vv$ and $\Vw$.

\begin{Theorem}\label{thm:special-formula}
There is a primitive recursive function, $\phi \mapsto \Special{\phi}$,
such that, for any formula
$\phi \in \LI$ containing $k$ vector variables counted in the sense described above, $\Special{\phi} \in \LI$ is a special formula that is
$k$-equivalent to $\phi$.
\end{Theorem}

\Proof
We will show by induction that every formula $\phi$ with $k$ vector
variables is $k$-equivalent to a special formula and it will be clear from
the proof that a suitable special formula can be calculated as a
primitive recursive function of $\phi$.

We may replace $\all{\ldots} \ldots$ by $\Not \ex{\ldots} \lnot \ldots$
throughout, so the cases we have to consider are:
\Lab{i} quantifier-free formulas (and hence in particular atomic formulas),
\Lab{ii} logical negation,
\Lab{iii} scalar existential quantification,
\Lab{iv} vector existential quantification
and \Lab{v} the binary propositional connectives.

{\em(i):} If $\phi$ is quantifier-free with free variables
$\Vv_1, \ldots, \Vv_k$,
Lemma~\ref{lma:standard-form} provides a quantifier-free formula $\PreSpecial{\phi}$
that is 0-equivalent to $\phi$ and in which vector terms only occur in
in terms of the form $\inner{\Vv_i, \Vv_j}$.
Replacing each occurrence of
$\inner{\Vv_i, \Vv_j}$ in $\PreSpecial{\phi}$ where $i > j$ by $\inner{\Vv_j, \Vv_i}$
then gives a special formula that is 0-equivalent and hence
$k$-equivalent to $\phi$ for any $k$.

For steps \Lab{ii} to \Lab{iv},
assume that $\phi$ is $k$-equivalent to a special
formula, $\sigma$ with the same free variables.

\Lab{ii} Like $\phi$, $\Not \phi$ contains $k$ vector variables and 
is $k$-equivalent to the special formula $\Not \sigma$.

\Lab{iii} Again like  $\phi$,  $\ex{x} \phi$ contains $k$ vector variables and is
$k$-equivalent to the special formula $\ex{x} \sigma$.

\Lab{iv} If $\Vv$ does not appear free in $\phi$, then $\ex{\Vv} \phi$ contains
either $k$ or $k+1$ vector variables and is logically equivalent to $\phi$, and
hence $k$-equivalent and so also $(k+1)$-equivalent to the special formula $\sigma$.
If $\Vv$ does appear free in $\phi$, then $\phi$ and $\ex{\Vv}\phi$ both contain $k$
vector variables.
Let $m = |v(\phi)|$ and $n = |s(\phi)|$.
Let $\Tuple{\Vv} = (\Vv_1, \ldots, \Vv_m)$ enumerate $v(\phi)  = v(\sigma)$
so that $\Vv \equiv \Vv_m$
and let $\Tuple{z} = (z_1, \ldots, z_n)$ enumerate $s(\phi) = s(\sigma)$.
Since $\sigma$ is special it has the form $\psi(\pUG(\Tuple{\Vv}), \Tuple{z})$, where $\psi(w_1, w_2, \ldots, w_{m(m+1)/2}, \Tuple{z})$ is a formula
in the language of an ordered field.
Let $\chi \IsDef \psi[\inner{\Vv_1, \Vv_1}/w_1, \ldots, \inner{\Vv_{m-1}, \Vv_{m-1}}/w_{(m-1)m/2})]$ be the result of substituting the terms of the
sequence
$\pUG(\Vv_1, \ldots, \Vv_{m-1})$ for $w_1, \ldots, w_{(m-1)m/2}$ in $\psi$.
I claim that $\ex{\Vv} \phi \equiv \ex{\Vv_m} \phi$ is $k$-equivalent to the
special formula $\sigma_1$ defined as follows
where the $x_i$ and the $y_i$ are fresh variables:
\[
\begin{array}{rccl}
\sigma_1 &\IsDef& 
  \multicolumn{2}{l}{\ex{x_1\; \cdots x_m \; y_1\; \cdots\; y_m}} \\
&&&
    {\displaystyle \bigwedge_{i=1}^{m-1} x_i} = \sum_{j=1}^{m-1}y_j\inner{\Vv_i, \Vv_j}\\
&&\And&
    x_m = \sum_{i=1}^{m-1}\sum_{j=1}^{m-1}y_iy_j\inner{\Vv_i, \Vv_j} + y_m^2\\
&&\And&
    \chi[x_1/w_{(m-1)m/2+1}, \ldots, x_m/w_{m(m+1)/2}]
\end{array}
\]
To see that $\ex{\Vv_m}\phi$ is indeed $k$-equivalent to $\sigma_1$,
let $V$ be an inner product space of dimension at least $k$,
let $\Vp_1, \ldots, \Vp_{m-1} \in V$ and let $\Tuple{c} \in \real^{n}$.
We have to show that $V \models (\ex{\Vv_m}\phi)(\Vp_1, \ldots, \Vp_{m-1}, \Tuple{c})$
iff $V \models \sigma_1(\Vp_1, \ldots, \Vp_{m-1}, \Tuple{c})$:

$\Imp$: assume $V \models (\ex{\Vv_m}\phi)(\Vp_1, \ldots, \Vp_{m-1}, \Tuple{c})$,
so there is $\Vp_m \in V$, such that $V \models \phi(\Vp_1, \ldots, \Vp_m, \Tuple{c})$.
Since $\phi$ and $\sigma$ are $k$-equivalent, 
$V \models \sigma(\Vp_1, \ldots, \Vp_m, \Tuple{c})$,
i.e.,
$V \models \psi(\pUG(\Vp_1, \ldots, \Vp_m), \Tuple{c})$.
Applying Lemma~\ref{lma:gram-matrix} to the $\Vp_i$ and the Gram matrix $M =
G(\Vp_1, \ldots, \Vp_m)$, we obtain $b_1, \ldots, b_{m} \in \real$ satisfying
equations {\Lab{a}} and {\Lab{b}} of the lemma, so that if we
interpret $x_i$ as $M_{im}$ and $y_i$ as $b_i$, $1 \le i \le m$,
the matrix of $\sigma_1$ holds, so that $V \models
\sigma_1(\Vp_1, \ldots, \Vp_{m-1}, \Tuple{c})$ as required.

$\Pmi$: assume $V \models \sigma_1(\Vp_1, \ldots, \Vp_{m-1}, \Tuple{c})$, so that
there are $a_i, b_i \in \real$, $1 \le i \le m$,
such that the matrix of $\sigma_1$ holds if we interpret the $x_i$ as the
$a_i$, the $y_i$ as the $b_i$ and $\Vv_1, \ldots, \Vv_{m-1}$ as $\Vp_1,
\ldots, \Vp_{m-1}$.
Let $M$ be the $m \times m$ matrix with  $M_{ij} = \inner{\Vp_i, \Vp_j}$, $1 \le i, j < m$ and $M_{im} = M_{mi} = b_i$, $1 \le i \le m$.
Then the assumptions of Lemma~\ref{lma:gram-matrix} hold for $M$
as do equations {\Lab{a}} and {\Lab{b}} of the lemma, which thus gives us
$\Vp_m \in V$ such that $M_{ij} = \inner{\Vp_i, \Vp_j}$,  $1 \le i, j \le m$.
The final conjunct in the matrix of $\sigma_1$ then implies that $V \models
\psi(\pUG(\Vp_1, \ldots, \Vp_m), \Tuple{c})$, i.e. $V \models
\sigma(\Vp_1, \ldots, \Vp_m, \Tuple{c})$.  As $\phi$ and $\sigma$ are $k$-equivalent, we
must have $V \models \phi(\Vp_1, \ldots, \Vp_m, \Tuple{c})$, so that $V \models
(\ex{\Vv_m}\phi)(\Vp_1, \ldots, \Vp_{m-1}, \Tuple{c})$ as required.

\Lab{v} At this point, the proof for formulas in prenex normal form would be
complete. However, putting a formula into prenex normal form can cause an
exponential explosion in the number of variables that it contains when counted
in our frugal sense and this would make algorithms based on our results less
efficient.  So for the final step, assume that $\phi$ contains $k$ vector
variables and is $k$-equivalent to the special formula $\sigma$, while $\phi'$
contains $k'$ vector variables and is $k'$-equivalent to the special formula
$\sigma'$.
Let ${\circ}$ be any binary propositional connective.
It is easy to see that, $\phi \circ \phi'$ is $\Max\{k, k'\}$-equivalent
to $\sigma \circ \sigma'$.
But if $k''$ is the number of vector variables in $\phi \circ \phi'$,
we must have $k'' \ge \Max\{k, k'\}$, so $\phi \circ \phi'$ is also
$k''$-equivalent to the special formula  $\sigma \circ \sigma'$.
\Done

The only special feature of the field of real numbers used in the above
proof is that it is {\em euclidean}, i.e., all
positive elements have square roots (which is needed to ensure the existence of
a unit vector in any given direction). Over a non-euclidean field,
both the proof and the statement of the theorem break down: over the
field of rational numbers, there is a countable infinity of distinct
isomorphism classes of 1-dimensional inner product spaces indexed
by square-free positive integers, the class corresponding to
$m$ being characterized by the sentence $\ex{\Vv} \inner{\Vv, \Vv} = m$.

The theorem immediately gives us an effective quantifier elimination procedure in the infinite-dimensional case:

\begin{Corollary}\label{:inf-dim-qelim}
There is a primitive recursive function, $\phi \mapsto \QElim{\phi}$,
such that if $\phi \in \LI$, $\QElim{\phi} \in \LI$ is a quantifier-free
formula that is equivalent to $\phi$ modulo either of the theories
$\IP{\infty}$ and $\HS{\infty}$.
\end{Corollary}
\Proof
First calculate the special formula $\Special{\phi}$ given by the
theorem; $\Special{\phi}$ will be equivalent to $\phi$ in any
infinite-dimensional inner product space and has the form 
$\psi(\pUG(\Vv_1, \ldots, \Vv_m), \Tuple{x})$ where
the $\Vv_i$ are the free vector variables of $\phi$ and $\psi$
is a formula in the language of an ordered field.
Now apply quantifier elimination for real closed fields to $\psi$,
giving an equivalent quantifier-free formula, $\chi$, say,
and put $\QElim{\phi} = \chi(\pUG(\Vv_1, \ldots, \Vv_m), \Tuple{x})$.
\Done

It follows that $\IP{\infty}$ are $\HS{\infty}$ are both decidable
and actually coincide.
However, as we are also interested in the decision problem for $\IP{}$,
$\IP{\fin}$, $\HS{}$ and $\HS{\fin}$, we will take a different line,
using the following theorem to justify an alternative approach.

\begin{Theorem}\label{thm:stability}
Let $\phi$ be a sentence of $\LI$ containing $k$ vector variables and let $V$ be
any inner product space of (possibly infinite) dimension $d \ge k$. Then $\phi$
holds in $V$ iff it holds in $\real^{k}$.
\end{Theorem}
\Proof
By Theorem~\ref{thm:special-formula}, $\phi$ is $k$-equivalent to a special
formula, but a special formula with no free variables is just a sentence in the
language of an ordered field  and its truth is independent of the choice of
vector space, so any space of dimension at least $k$, e.g., $\real^{k}$, will
serve to test the truth of $\phi$.
\Done

\begin{Lemma}\label{lma:IPn-reduces-to-Rn}
There is a primitive recursive function that maps a sentence $\phi$ of $\LI$ and a
natural number $n$ to a sentence $\Res{\phi}{n}$ in the language
of an ordered field
such that $\real^n \models \phi \Iff \Res{\phi}{n}$, i.e.,
$\phi$ holds in $\real^n$ iff $\Res{\phi}{n}$ holds in the
ordered field $\real$.
\end{Lemma}
\Proof
We describe a primitive recursive algorithm that constructs the sentence
$\Res{\phi}{n}$ and show that it holds in $\real$
iff $\phi$ holds in the standard
$n$-dimensional inner product space $\real^n$, which proves the lemma.

If $n = 0$,  $\Res{\phi}{0}$ is obtained from $\phi$ by deleting all vector
quantifiers, replacing all inner products by scalar 0 and replacing all vector
equations by the scalar equation $0 = 0$. Evidently $\real^0
\models \phi$ iff $\real^0 \models \Res{\phi}{0}$.

If $n \ge 1$, pick $n$ fresh vector variables $\Vb_1, \ldots, \Vb_n$
and, for each vector variable $\Vv$ occurring in $\phi$, pick
$n$ fresh scalar variables $x^\Vv_1, \ldots, x^\Vv_n$. Replace each vector
quantifier $\all{\Vv}$ (resp. $\ex{\Vv}$) in $\phi$ by the string of scalar
quantifiers $\all{x^\Vv_1\;\cdots\;x^\Vv_n}$ (resp.
$\ex{x^\Vv_1\;\cdots\;x^\Vv_n}$) and replace all other occurrences of $\Vv$ by
$x^\Vv_1\Vb_1 + \ldots + x^\Vv_n\Vb_n$. Let the resulting formula be $\phi_1(\Vb_1, \ldots, \Vb_n)$.
Clearly, $\real^n \models \phi$ iff $\real^n \models_{\real^n} \phi_1(\Tuple{\Ve})$
where $\Tuple{\Ve} = (\Ve_1, \ldots, \Ve_n)$ 
is the standard basis for $\real^n$.
By Lemma~\ref{lma:standard-form}, $\phi_1$ is
equivalent to a special formula $\phi_2(\pUG(\Vb_1, \ldots, \Vb_n))$
where $\phi_2(w_1, \ldots, w_{n(n+1)/2})$ is a formula in the language
of an ordered field.
Writing the Kronecker symbol $\delta_{ij}$
to stand for the constant $1$ when $i = j$ and the constant $0$ otherwise,
define $\Res{\phi}{n} \IsDef
        \phi_2[\delta_{ij}/w_{n(i-1)+j} \ST 1 \le i \le j \le n]$.
We then have $\real^n \models \phi$ iff
 $\real^n \models \phi_1(\Tuple{\Ve})$ iff
 $\real^n \models \phi_2(\pUG(\Tuple{\Ve}))$ iff
$\Res{\phi}{n}$ holds in $\real$.
\Done

In the construction of $\Res{\phi}{n}$ in the above proof, an alternative way of
eliminating vector variables from the formula $\phi_1$ is to rearrange vector
equations into the form  $t_1\Vb_1 + \ldots + t_n\Vb_n = \VO$
which may then be replaced by $t_1 = t_2 = \ldots = t_n = 0$ before
applying the method of Lemma~\ref{lma:standard-form} to eliminate
inner products. This is more
efficient and also avoids introducing multiplication, which might be
practically beneficial when working in the additive fragment of an extended
language including a richer supply of vector constants.

In Section~\ref{subsubsec:poss-qe} we defined sentences $\pD_{{\le}n}$ for $n \in \nat$
that hold in a vector space iff the space has finite dimension less than or
equal to $n$. Let us define $\pD_0 := \pD_{{\le}0}$ and $\pD_{n+1} := \pD_{{\le}n+1}
\And \Not \pD_{{\le}n}$ so that the sentence $\pD_n$ holds iff the dimension is
exactly $n$. We use these sentences to reduce the theories of interest to the
theory $\IP{{\le}k}$ of inner product spaces of dimension at most $k$.

\begin{Theorem}\label{thm:ip-quant-elim}
Let $\phi$ be a sentence of $\LI$ containing $k$ vector variables; if $k = 0$, let
$\Star{\phi} \IsDef \Res{\phi}{0}$, otherwise define $\Star{\phi}$ by,
$$
\Star{\phi} \IsDef (\pD_0 \And \Res{\phi}{0}) \Or (\pD_1 \And \Res{\phi}{1}) \Or
\ldots \Or (\pD_{k-1} \And \Res{\phi}{k-1}) \Or (\Not \pD_{{\le}(k-1)} \And \Res{\phi}{k}).
$$
Then $\Star{\phi}$ is equivalent to $\phi$ in any of the theories $\IP{}$,
$\IP{\fin}$, $\IP{\infty}$, $\HS{}$, $\HS{\fin}$, and $\HS{\infty}$.
\end{Theorem}
\Proof
Let $V$ be any inner product space. If $V$ has infinite dimension or finite
dimension $d \ge k$, then $\pD_n$ is false in $V$ for $n \le k - 1$ and $\Not
\pD_{{\le}k-1}$ is true, so $\Star{\phi}$ is equivalent in $V$ to $\Res{\phi}{k}$. But,
by Lemma~\ref{lma:IPn-reduces-to-Rn}, $\Res{\phi}{k}$ is true iff $\phi$ is true in $\real^k$, and
by Theorem~\ref{thm:stability}, $\phi$ is true in $V$ iff it is true in $\real^k$. If $V$ has
finite dimension $d < k$, then $\Star{\phi}$ is equivalent to $\Res{\phi}{d}$ which
is valid iff $\phi$ holds in $\real^d$ iff $\phi$ holds in $V$, since $V$ and
$\real^d$ are isomorphic. So irrespective of the dimension of $V$, $\phi$ holds
iff $\Star{\phi}$ holds.
Noting that our methods of proof make no assumptions
about completeness this completes the proof of the theorem.
\Done

\begin{Corollary}\label{:ip-valid-dimensions}
For every sentence $\phi$ of $\LI$ there is a subset $D_\phi$ of $\nat \cup
\{\infty\}$ such that $\phi$ holds in an inner product space $V$ iff
$\Dim{V} \in D_\phi$.
Moreover $D_\phi$ is either a finite subset of $\nat$ or the complement
of a finite subset of $\nat$ and can be effectively computed from $\phi$.
\end{Corollary}
\Proof
First, calculate $\phi^{*}$ as in the theorem and then apply the quantifier elimination
algorithm for the first-order theory of real arithmetic to determine the truth
values of the sentences  $\Res{\phi}{i}$ that appear in $\phi^{*}$.
Now simplify to
give either {\em(i)} a (possibly empty) disjunction of the form $\pD_{i_1} \Or
\ldots \Or \pD_{i_m}$ or {\em(ii)} a disjunction of the form $\pD_{i_1} \Or \ldots
\Or \pD_{i_m} \Or \Not \pD_{{\le}(k-1)}$ (where $k > i_m$ is the number of vector
variables in $\phi$). In both cases, the truth of the result is determined
by a set $D_\phi$ of dimensions: in case {\em(i)}, we have $D_\phi = \{i_1, \ldots,
i_m\}$ which is a finite subset of {\nat}, while in case {\em(ii)} $D_\phi$ is the complement in $\nat \cup
\{\infty\}$ of the finite subset $\{0, \ldots, k - 1\} \Diff \{i_1, \ldots,
i_m\}$.
Let us represent $D_\phi$ as a pair $(t, X)$, where $t \in \{0,
1\}$ and $X$ is a finite set of natural numbers, $D_\phi$ being given by $X$ when
$t = 0$ and its complement when $t = 1$.
Since the construction of $\phi^{*}$ and the $\Res{\phi}{i}$ is primitive
recursive, we have an effective procedure for computing the representation of
$D_\phi$.
\Done

\begin{Corollary}\label{:ip-axiomatizability}
A class $\cal C$ of structures for the language $\LI$ is axiomatizable (resp.
recursively axiomatizable) iff it comprises all inner product spaces $V$ such
that $\Dim{V} \in D$ for some $D \subseteq \nat \cup \{\infty\}$ that is
either finite or contains $\infty$  (resp.  either finite or the complement of
a recursively enumerable subset of $\nat$).
\end{Corollary}
\Proof
Recall that a class of structures for a language is said to be (recursively)
axiomatizable iff it comprises all models of some (recursive) set of axioms.
If $A$ is any set of sentences of $\LI$, then, by the previous corollary, $V$ is
a model of $A$ iff $\Dim{V} \in \bigcap_{\phi\in A}D_\phi$ where each $D_\phi$ is either
a finite set of natural numbers or the complement in $\nat \cup \{\infty\}$ of
a finite set of natural numbers. A subset $D$ of $\nat \cup \{\infty\}$ can be written
as such an intersection iff it is either a finite set of natural numbers or
contains $\infty$.

The assertion about recursive axiomatizability is an easy exercise in
recursion theory: in one direction, test for non-membership of $D$ using an
algorithm that on input $d$, enumerates the sentences of $A$ checking for each
sentence in turn whether it excludes models of dimension $d$; in the other
direction, observe that a non-empty {r.e.} set of finite dimensions may be
excluded by an {r.e.} set of axioms and then use the well-known trick of
replacing the {r.e.} set $\phi_1, \phi_2, \phi_3, \ldots$ by the recursive set $\phi_1, \phi_1
\And \phi_2, \phi_1 \And \phi_2 \And \phi_3, \ldots$ to get a recursive axiomatization.
\Done

\begin{Theorem}\label{thm:ip-decidable}
The theories $\IP{}$, $\IP{\fin}$, $\IP{\infty}$, $\HS{}$, $\HS{\fin}$  and
$\HS{\infty}$ are all decidable. Moreover
$\IP{} = \IP{\fin} = \HS{} = \HS{\fin}$ and
$\IP{\infty} = \HS{\infty}$.
\end{Theorem}
\Proof
By Corollary~\ref{:ip-valid-dimensions}, given a sentence $\phi$ of $\LI$, we can
effectively calculate the set $D_\phi \subseteq \nat \cup \{\infty\}$ of
dimensions in which $\phi$ holds and for some finite $X \subseteq \nat$, either
$D_\phi = X$ or $D_\phi = (\nat \cup \{\infty\}) \Diff X$.  If $D_\phi = X$, then $\phi$
does not belong to any of the theories listed.  If $D_\phi = (\nat \cup
\{\infty\}) \Diff X$, then $\phi$ certainly belongs to both $\IP{\infty}$ and
$\HS{\infty}$, while $\phi$ belongs to $\IP{}$, $\IP{\fin}$, $\HS{}$ or
$\HS{\fin}$ iff $X$ is empty. Thus we have an effective procedure for deciding
membership for each of the theories.  Since the theories $\IP{\infty}$ and
$\HS{\infty}$ have a common decision procedure they are equal and similarly
$\IP{}$, $\IP{\fin}$, $\HS{}$ and $\HS{\fin}$ are all equal.
\Done

For $d \in \nat$ there is exactly one inner product space of dimension
$d$ up to isomorphism.  Corollary~\ref{:ip-valid-dimensions} implies
that there is exactly one infinite-dimensional inner product space up to
elementary equivalence.  By contrast, it can be shown that, up to
elementary equivalence, there are $\C = \Card{\real}$ distinct
$d$-dimensional normed spaces for each $d$, $2 \le d \in \nat \cup
\{\infty\}$.

\Section{Decidable fragments of the theory of normed spaces}\label{sec:decidable-fragments}

Although we have shown that the general theory of normed spaces is
undecidable, there are some significant decidable fragments.
In this section, we will find that the purely universal and purely existential
fragments are both decidable via reductions to the first-order theory of the
real numbers. The reduction for purely existential sentences is
very simple, but for purely universal sentences, the reduction involves an
interesting geometrical construction.
In Section~\ref{sec:ea-ae-fragments} we will find
that the $\fEA$ and $\fAE$ fragments are undecidable,
so these results are the best possible of their type.

Consider a sentence in the language of normed spaces that is in prenex normal
form and contains no universal quantified vector variables: clearly such a
sentence $\phi$ holds in all normed
spaces iff it holds in the trivial normed space 0. We therefore obtain
a decision procedure for valid sentences of this form by striking out all vector quantifiers,
replacing all norm expressions by $0$ and all vector equations by $0 = 0$ and
then applying the decision procedure for the first-order theory of the real numbers.
In particular, the set of valid purely existential sentences is decidable.

As we shall now see the set of true purely universal sentences in the language of
normed spaces is also decidable, but the decision procedure and its
verification are much less trivial:  the crux of the argument lies in deciding
satisfiability of a set of bounds on the norms of a finite set of vectors, so
we start by considering how to define a norm satisfying a system of
constraints.

A subset of $X$ of a vector space $V$ is said to be {\em symmetric} if $X = -X$
where $-X =\{-\Vv \mid \Vv \in X\}$.  Given a subset $Y$ of $V$ we define the
{\em symmetric convex hull} of $Y$, written $\Sconv{Y}$, to be the intersection
of the set of all symmetric convex sets containing $Y$.  $\Sconv{Y}$ is itself
symmetric and convex and it is easy to verify that $\Sconv{Y}$ is the convex
hull of $Y \Union -Y$.  If $\Vv \in \Sconv{Y}$, then, by symmetry, $-\Vv \in
\Sconv{Y}$ and then, by convexity, the line segment $[-\Vv, \Vv]$ is contained
in $\Sconv{Y}$, i.e., $c\Vv \in \Sconv{Y}$ for any $c$ with $|c| \leq 1$.

\begin{Lemma}\label{lma:sconv}
Let $X = \{\Vx_1, \ldots, \Vx_n\}$ be a non-empty finite subset of  a
vector space. Then the symmetric convex hull of $X$ is given by:
$$\Sconv{X} =
  \left\{ \sum_{i=1}^n c_i \Vx_i \;\left|\; \sum_{i=1}^n |c_i| \leq 1 \right.\right\}.
$$
\end{Lemma}
\Proof
Write $D = \{ \sum_{i=1}^n c_i \Vx_i \mid \sum_{i=1}^n |c_i| \leq 1 \}$.  It is
easy to check that $D$ is convex, symmetric and contains $X$, so $\Sconv{X} \subseteq D$.  Conversely, let
$\Vv \in D$, so $\Vv = \sum_{i=1}^n c_i \Vx_i$ for some $c_i$ where $c =
\sum_{i=1}^n |c_i| \leq 1$.  If $c = 1$, then $\Vv$ is a convex combination of
the points $\pm\Vx_i$ and $\Vv \in \Sconv{X}$ by the
remarks above.  If $c = 0$, then trivially $\Vv  = \VO  \in \Sconv{X}$.
So assume $0 < c < 1$, so that $\Vv = c\sum_{i=1}^n (c_i/c)
\Vx_i$ and we have $\sum_{i=1}^n |c_i/c| = 1$. Hence $\Vv$ can be written as
$c\Vw$ where $|c| \leq 1$ and $\Vw \in \Sconv{X}$ (by
the case $c = 1$ just considered) so by the remarks above $\Vv \in
\Sconv{X}$.
\Done

\begin{Lemma}\label{lma:risconv}
Let $Y = \{\Vx_1,\ldots,\Vx_n\}$ be a non-empty finite subset of a vector space
$V$ and let $D = \Sconv{Y}$ be its symmetric convex hull. Then {\em(i)} $D$ is
the unit disc of a norm on the subspace $W$ of $V$ spanned by $Y$ and
{\em(ii)} if $S$ is the unit circle for this norm, then
$$
D \Diff S = \left\{ \sum_{i=1}^n c_i \Vx_i\;\left|\;\sum_{i=1}^n |c_i| < 1 \right.\right\}.
$$
\end{Lemma}
\Proof
For {\em(i)}, as $D$ is certainly convex, it will satisfy the criteria for  the
unit disc of a norm on $W$ if it meets every line through the origin in $W$ in
a line segment $[-\Vv, \Vv]$ with $\Vv \not= 0$.
By the Minkowski-Weyl theorem, $D$, which is the convex hull of a finite set of
points, can be written as the intersection of a finite set of closed
halfspaces.  Hence if $l$ is any line through the origin in $W$, $l \cap D$ is
the intersection of $l$ and a finite set of closed half-lines, and hence, as it
is non-empty, bounded and symmetric about the origin, it must be 
the line segment, $[-\Vv, \Vv]$ for some $\Vv$. We have only to show that $\Vv
\not= \VO$.  To see this let $\Vw$ be any point of $l \Diff \{\VO\}$.  Since $\Vw
\in W$, there are $c_i$ such that $\Vw = \sum_{i=1}^n c_i \Vx_i$.  If we let $c
= \sum_{i=1}^n |c_i|$, then $c \not= 0$, and by Lemma~\ref{lma:sconv}, $\Vw/c
\in D$, but then $\Vw/c \in D \cap l = [-\Vv, \Vv]$ and as $\Vw/c \not= 0$ we
must have $\Vv \not= 0$.

For {\em(ii)}, note that $\Vv \in D \Diff S$ iff there is a $d > 1$ such that
$d\Vv \in D$. By Lemma~\ref{lma:sconv}, $d\Vv \in D$ iff $d\Vv$ can be written
as $\sum_{i=1}^n d_i \Vx_i$ with $\sum_{i=1}^n |d_i| \leq 1$ and this holds
for $d > 1$ iff $\Vv$ can be written as $\sum_{i=1}^n c_i \Vx_i$ with
$\sum_{i=1}^n |c_i| < 1$.
\Done

\begin{Lemma}\label{lem:normexistsunit}
Let $\Vx_1,\ldots,\Vx_n$ and $\Vy_1,\ldots,\Vy_m$ be vectors in a vector space
$V$. Then there exists a norm $\|\_\|$ on $V$ such that $\|\Vx_i\| \leq 1$ for
all $i, 1 \leq i \leq n$ and $\|\Vy_j\| \geq 1$ for all $j, 1 \leq j \leq m$ iff no
$\Vy_k$ is expressible as $\Vy_k = \sum_{i=1}^n c_i \Vx_i$ with $\sum_{i=1}^n
|c_i| < 1$.
\end{Lemma}
\Proof
If a norm satisfies the stated bounds, then it is indeed impossible that any
$\Vy_k = \sum_{i=1}^n c_i \Vx_i$ with $\sum_{i=1}^n |c_i| < 1$, for then by the
triangle inequality $\|\Vy_k\| \leq  \sum_{i=1}^n \|c_i \Vx_i\| = \sum_{i=1}^n |c_i|
\|\Vx_i\| \leq \sum_{i=1}^n |c_i| < 1$, contradicting $\|\Vy_k\| \geq 1$.

Conversely, suppose no $\Vy_k$ is expressible as $\Vy_k = \sum_{i=1}^n c_i
\Vx_i$ with $\sum_{i=1}^n |c_i| < 1$. By Lemma~\ref{lma:risconv}, we can define
a norm $\|\_\|_0$ on the span $V_0$ of $\Vx_1,\ldots,\Vx_n$ with $D =
\mbox{sconv}(\{\Vx_1,\ldots,\Vx_n\})$ as its unit disc. Let $V_1$ be a
complementary subspace of $V_0$, so that every
$\Vv \in V$ is uniquely expressible as $\Vv = \Vv_0 + \Vv_1$ for $\Vv_0 \in
V_0$ and $\Vv_1 \in V_1$. Let $\|\_\|_1$ be an arbitrary norm on $V_1$, e.g.
defined using an inner product w.r.t. some basis. For any $B > 0$, the norm
$\|\Vv_1\|_B = B \|\Vv_1\|_1$ is also a norm on $V_1$, and $\|\Vv\| =
\|\Vv_0\|_0 + \|\Vv_1\|_B$ is a norm on $V$. I claim that for sufficiently
large $B$, this satisfies the constraints $\|\Vy_j\| \geq 1$. First, if $\Vy_j
\in V_0$, then this follows immediately since the assumption implies, by Lemma
\ref{lma:risconv}, that $\Vy_j$ is not in $\{\Vw \mid \|\Vw\|_0 < 1\}$.  On the
other hand, all the $\Vy_j \not\in V_0$ can be written $\Vy_j = \Vw_j + \Vz_j$
for $\Vw_j \in V_0$, $\Vz_j \in V_1$ with $\Vz_j$ nonzero. To ensure
$\|\Vy_j\| = \|\Vw_j\|_0 + B \|\Vz_j\|_1 \geq 1$, it suffices to choose
$B > \Max\{1/\Norm{\Vz_1}_1, \ldots, 1/\Norm{\Vz_m}_1\}$. \Done

\begin{Theorem}\label{thm:normexists}
Let $\Vx_1,\ldots,\Vx_n$ and $\Vy_1,\ldots,\Vy_m$ be vectors in a vector space
$V$, and let $b_1,\ldots,b_n$ and $d_1,\ldots,d_m$ be real numbers such that $b_i \not= 0$ for some $i, 1 \le i \le n$. Then
there exists a norm $\|\_\|$ on $V$ such that $\|\Vx_i\| \leq b_i$ for
all $i, 1 \leq i \leq n$, and $\|\Vy_j\| \geq d_j$ for all $j, 1 \leq j \leq m$, iff
the following conditions hold:

\begin{itemize}

\item For all $1 \leq i \leq n$, $b_i \geq 0$;

\item For all $1 \leq i \leq n$, if $b_i = 0$ then $\Vx_i = 0$;

\item No $\Vy_j$ is expressible as $\Vy_j = \sum_{i=1}^n c_i \Vx_i$ with
$\sum_{i=1}^n |c_i| b_i < d_j$.

\end{itemize}
\end{Theorem}
\Proof
If a norm satisfying the claimed inequalities exists, then all three properties
follow immediately from the norm properties, the last one using the triangle
inequality just as in the proof of Lemma~\ref{lem:normexistsunit}.

Conversely, suppose the three properties hold. In order to construct a norm
satisfying the inequalities, we can assume without loss of generality that all
$b_i > 0$, because by the second property, if $b_i = 0$ then $\Vx_i = 0$ and so
any norm at all satisfies $\|\Vx_i\| \leq b_i$. Similarly, we can assume that
each $d_j > 0$ because if $d_j \leq 0$ the constraint $\|\Vy_j\| \geq d_j$ is
automatically satisfied.

Define $\Vu_i = \Vx_i / b_i$ and $\Vv_j = \Vy_j / d_j$. Note that no $\Vv_j$ is
expressible as $\Vv_j = \sum_{i=1}^n c_i \Vu_i$ with $\sum_{i=1}^n
|c_i| < 1$, because then $\Vy_j = d_j \sum_{i=1}^n c_i \Vu_i =
\sum_{i=1}^n (d_j c_i / b_i) \Vx_i$, and $\sum_{i=1}^n |d_j c_i / b_i| b_i =
d_j \sum_{i=1}^n |c_i| < d_j$, contrary to the third condition. Therefore by
Lemma~\ref{lem:normexistsunit}, there is a norm on $V$ satisfying
$\|\Vu_i\| \leq 1$ for $1 \leq i \leq n$ and $\|\Vv_j\| \geq 1$ for $1 \leq
i \leq m$. I.e.,  $\|\Vx_i\| \leq b_i$ for
all $1 \leq i \leq n$ and $\|\Vy_j\| \geq d_j$ for all $1 \leq j \leq m$.
\Done

We can immediately obtain a simpler result if we seek conditions allowing us to
set the specific values of the norms of a finite set of vectors:

\begin{Corollary}\label{cor:normexisteq}
Let $\Vx_1,\ldots,\Vx_n$ be vectors in a real vector space $V$ and let
$b_1,\ldots,b_n$ be real numbers. Then there exists a norm $\|\_\|$ on $V$ such
that $\|\Vx_i\| = b_i$ for all $i, 1 \leq i \leq n$ iff:

\begin{itemize}

\item For all $i, 1 \leq i \leq n$, $b_i \geq 0$;

\item For all $i, 1 \leq i \leq n$, if $b_i = 0$ then $\Vx_i = 0$;

\item For each $k, 1 \leq k \leq n$ there are no real numbers $c_1,\ldots,c_n$
such that some $\Vx_k = \sum_{i=1}^n c_i \Vx_i$ with $\sum_{i=1}^n |c_i| b_i <
b_k$.

\end{itemize}
\end{Corollary}
\Proof
The case when each $b_i = 0$ is evident.  If some $b_i \not= 0$, then apply Theorem~\ref{thm:normexists} with $m = n$, $\Vx_i = \Vy_i$ and $b_i
= d_i$.
\Done

\begin{Corollary}\label{cor:ms-a-dec}
The set of valid purely universal sentences in the language of normed spaces
is decidable.
\end{Corollary}
\Proof
If $\sigma$ is a purely universal sentence in prenex
normal form $\all{\ldots}\psi$, $\sigma$ is true iff $\Not \psi$ is unsatisfiable.
So it suffices to give a decision procedure for satisfiable quantifier-free
formulas. So let $\phi$ be quantifier-free say
with free variables given by
$v(\phi) = \{\Vx_1,\ldots,\Vx_n\}$ and $s(\phi) = \{u_1,\ldots,u_m\}$. Introduce
additional scalar variables $b_1,\ldots,b_k$, one $b_i$ for each norm expression
$\|\Vy_i\|$ appearing in $\phi$. (Each such vector $\Vy_i$ can be written as
$p_1 \Vx_1 + \cdots + p_n \Vx_n$ for polynomials $p_i$, though the $p_i$ may
themselves involve other norm expressions.) Satisfiability of $\phi$ in a
normed space is equivalent to satisfiability of $\phi' \And \Ands_{i=1}^k
\|\Vy_i\| = b_i$, where $\phi'$ is $\phi$ with each $\|\Vy_i\|$ replaced by its
corresponding $b_i$, in a bottom-up fashion so that $\phi'$ does not contain
the norm operator. But by the corollary, this is equivalent to the
satisfiability in a vector space of the following formula:
$$
\phi'' \IsDef
\begin{array}{l}
    \phi' \And {}\\
   \Ands_{i=1}^k b_i \geq 0 \And {}\\
   \Ands_{i=1}^k (b_i = 0 \Imp \Vy_i = 0) \And {}\\
   \Ands_{i=1}^k (\all{c_1 \ldots c_k}
        |c_1| b_1 + \cdots + |c_k| b_k < b_i
        \Imp \Vy_i \not= c_1 \Vy_1 + \cdots + c_k \Vy_k).
\end{array}
$$

\noindent
The decision procedure of Theorem~\ref{thm:ip-decidable} applied to the
existential closure of $\phi''$ will then decide satisfiability of $\phi''$
and hence of $\phi$.
\Done

Note that, if the formula $\phi$ is satisfiable, then our methods give a
norm on $\real^n$, whose unit disc may be taken to be a polyhedron, together with a
satisfying assignment for $\phi$ in $\real^n$ under that norm. Thus, at
least in principle, the above decision procedure can be extended to give
a counter-example if the input purely universal sentence is false.  It
is also noteworthy that the only instances of multiplication introduced
in the passage from $\phi$ to $\phi''$ are in the last conjunct of $\phi''$.  For
the case where the input sentence is purely additive, one can
develop a more efficient algorithm using a parametrised linear
programming technique.

If $K$ is an ordered field, define a {\em normed space over $K$} to be a
structure for the language $\LN$ of normed spaces in which the scalar sort and
its operations are interpreted in $K$ and which satisfies the usual axioms for
a norm.  The proofs above go through over any real closed field $K$ (for a proof
of the Minkowski-Weyl theorem that does not appeal to separation properties
that are only valid over $\real$, see, for example, \citet{Weyl-etkp}).
We therefore have a decision procedure for the
purely universal fragment of the theory of
normed spaces over any real closed field $K$. As this decision procedure is
independent of $K$, we may conclude that a universal sentence in $\LN$ holds
for all normed spaces over a real closed field $K$ iff it holds for all real
normed spaces.

\Section{The $\fEA$ and $\fAE$ fragments of the theory of normed spaces}\label{sec:ea-ae-fragments}

In this section we shall see that the $\fEA$ and $\fAE$ fragments of the
theory of normed spaces are
both undecidable.  Thus the results of Section~\ref{sec:decidable-fragments}
for the purely existential and
purely universal fragments are the best of their type.
The proofs given here do make use of multiplication, but the constructions
they use have since been adapted to give undecidability results for the additive
$\fEA$ and $\fAE$ fragments (over $\real$) and for theories of normed spaces over an
arbitrary ordered field ~\citep{arthan:additive-aia, arthan:any-field}.

The plan of this
section is as follows: we first prove undecidability for the $\fEA$
fragment by giving a purely existential characterization of the natural
numbers in a certain normed space $\sK$ (cf.  the proof of
Theorem~\ref{thm:ms-ea-valid-undec}); then we prove undecidability of
the $\fAE$ fragment using a normed space $\sL$ whose unit circle
includes an encoding of a periodic function; finally, we show that a
small adjustment to $\sL$ allows us to prove undecidability of the set
of all valid sentences of the form $\phi \Imp \psi$ where $\phi$ and $\psi$ are
purely universal, which, up to a logical equivalence, covers
undecidability for both the $\fAE$ and $\fEA$ fragments.

The first proof for the $\fEA$ fragment is based on an extremely simple
method for encoding the natural numbers in the unit disc of a
2-dimensional normed space.

\begin{figure}
\begin{center}
\includegraphics[angle=0,scale=0.4]{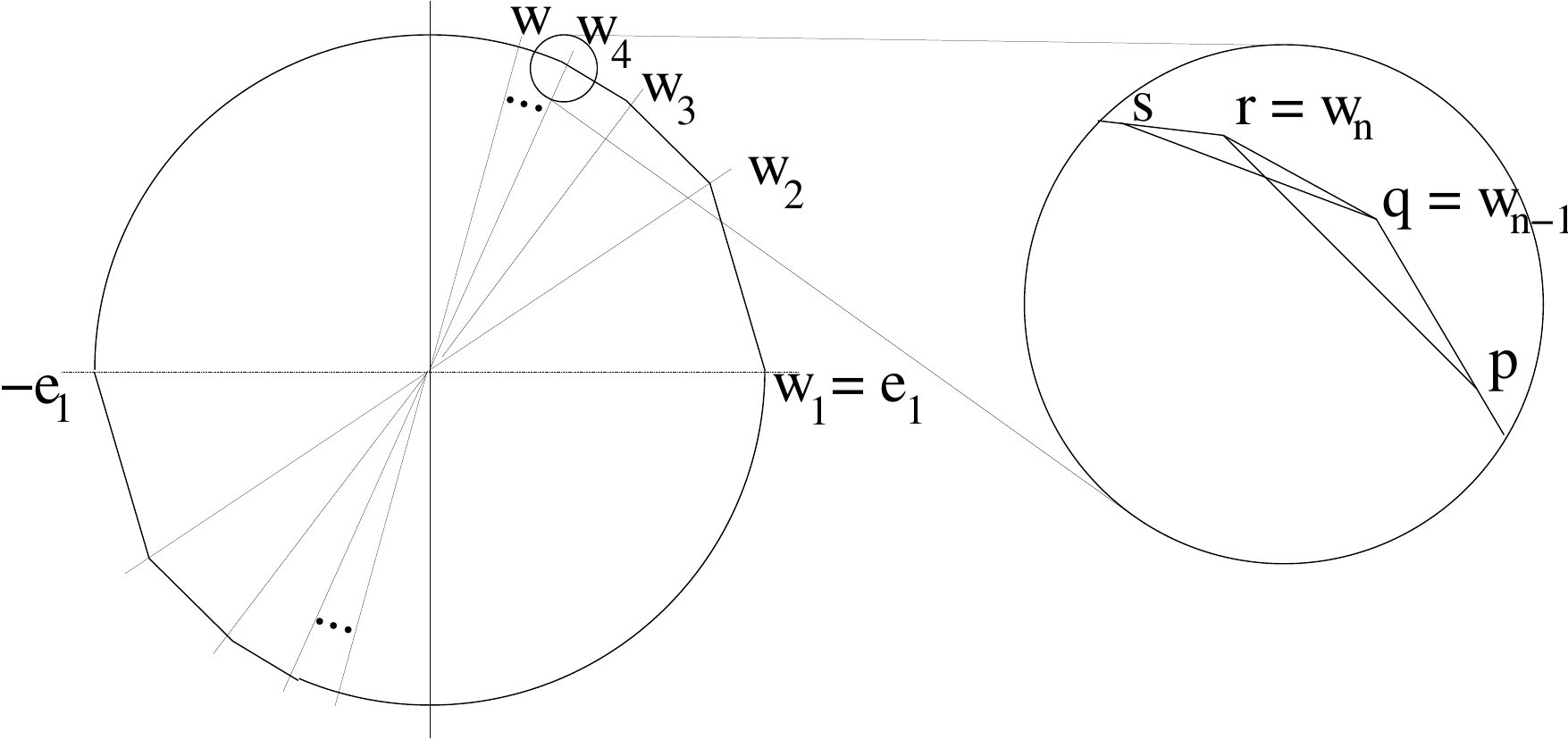}
\caption{The unit circle in the space  $\sK$ with a detail illustrating the predicate $\pA{\Vq, \Vr}$}
\label{fig:k}
\end{center}
\end{figure}

\begin{Theorem}\label{thm:ns-ae-def-nat}
There is a purely existential formula $\pN{x}$ in the language of normed spaces
such that for any $d \in \{2, 3, 4, \ldots\} \cup \{\infty\}$, there is a Banach
space $\sK^d$ of dimension $d$ in which $\pN{x}$ defines the natural numbers.
\end{Theorem}
\Proof
We consider the case $d = 2$ first.
Using the usual euclidean norm in the plane $\real^2$, define $\Vw_1 = \Ve_1$
and then, working anticlockwise around the unit circle, take  $\Vw_{n}$ to be
the unit vector with $\Norm{\Vw_{n}-\Vw_{n-1}}_e = \frac{1}{n!}$ for $n = 2, 3,
\ldots$ as illustrated in Figure~\ref{fig:k} (but not to scale).  Then
$\Norm{\Vw_{n} - \Vw_{1}}_e < \Sigma_{n=2}^{\infty}\frac{1}{n!} = e - 2 < 1 <
\sqrt{2} = \Norm{\Ve_1 - \Ve_2}_e$, and so
the $\Vw_n$ are all in the north-east quadrant and tend to a
limit $\Vw$.
Evidently $\Norm{\Vw}_e = 1$ and we may define $\sK = \sK^2$ to be $\real^2$  with
the norm $\Norm{\_}$ whose unit disc is the symmetric convex hull of
$A \cup \{\Vw_1, \Vw_2, \ldots\}$ where $A$ is the arc running anticlockwise from $\Vw$ to
$-\Ve_1$.
Note that this norm agrees with the euclidean norm on vectors
in the north-west and south-east quadrants, so, in particular, if $\Vp$
and $\Vq$ are unit vectors in the north-east quadrant, $\Norm{\Vp - \Vq} = \Norm{\Vp - \Vq}_e$.
Define predicates $\pA{\Vq, \Vr}$ and $\pN{x}$ as follows:
\begin{eqnarray*}
\pA{\Vq, \Vr} &\IsDef&
\begin{array}[t]{l@{}l}
\ex{\Vp\;\Vs} &
        \Norm{\Vp} = \Norm{\Vq} = \Norm{\Vr} = \Norm{\Vs} = 1 \And {} \\ &
        \Norm{(\Vp + \Vq)/2} =
        \Norm{(\Vq + \Vr)/2} =
        \Norm{(\Vr + \Vs)/2} = 1 \And {} \\ &
        \Norm{(\Vp + \Vr)/2} < 1 \And
        \Norm{(\Vq + \Vs)/2} < 1
\end{array} \\
\pN{x} &\IsDef&
\begin{array}[t]{l@{}l}
\ex{\Vp\;\Vq\;\Vr} &
        \pA{\Vp, \Vq} \And
        \pA{\Vq, \Vr} \And
        \Norm{\Vp - \Vq} > \Norm{\Vq - \Vr} \And {} \\ &
        \Norm{\Vp - \Vq} = (x + 4) \Norm{\Vq - \Vr}.
\end{array}
\end{eqnarray*}
For any $n \ge 3$, $\pA{\Vw_{n-1}, \Vw_{n}}$ holds in $\sK$ (choose
$\Vp \in [\Vw_{n-2}, \Vw_{n-1})$ and $\Vs \in (\Vw_{n}, \Vw_{n+1}]$ as in Figure~\ref{fig:k}).
Conversely, assume the matrix of $\pA{\Vq, \Vr}$ holds for some $\Vp$, $\Vq$,
$\Vr$ and $\Vs$: the conditions imposed imply that $\Vp$, $\Vq$, $\Vr$ and
$\Vs$ are pairwise distinct and that the line segments $[\Vp, \Vq]$, $[\Vq,
\Vr]$ and $[\Vr, \Vs]$ lie in the unit circle of $\sK$; moreover, $\Vp$,
$\Vq$ and $\Vr$ cannot be collinear, otherwise we would have
$\Norm{(\Vp + \Vr)/2} = 1$ and similarly $\Vq$, $\Vr$ and $\Vs$ cannot
be collinear; thus $\Vq$ and $\Vr$ must be adjacent
isolated extreme points of the unit
disc in $\sK$ and,  for some $n \ge 3$, we have $\pm\{\Vq, \Vr\} =
\{\Vw_{n-1}, \Vw_{n}\}$.  Since $\Norm{\Vw_{n-1} - \Vw_{n}} = \frac{1}{n!}$ for
all $n \ge 2$, it follows that, $\pN{x}$ holds iff for some $n \ge 3$,
$\frac{1}{n!} = \frac{x + 4}{(n+1)!}$, which holds (with $n = x + 3$) iff $x
\in \nat$.  Clearly $\pN{x}$ is equivalent to a purely existential formula
and the theorem is proved
for $d=2$.

For $d > 2$, let $V$ be a Hilbert space of dimension $d - 2$, and define
$\sK^d$ to be the 2-sum of $\sK$ and $V$, i.e., the product
vector space $\sK \times V$ equipped with the norm defined by
$\Norm{(\Vp, \Vv)} = \sqrt{\Norm{\Vp}^2 + \Norm{\Vv}^2}$.
That this is a norm making $\sK^d$ into a Banach space is readily verified.
We identify $\sK$ with the subspace $\sK \times 0$ of $\sK^d$.
It can be shown that if $\Va$ and $\Vb$ are distinct unit vectors
and if the line segment $[\Va, \Vb]$ is contained in the unit sphere of
$\sK^d$, then $[\Va, \Vb]$ is parallel to $\sK$
(see \citep{arthan:additive-aia} for a proof).
Hence there are $\Vp,
\Vq \in \sK$, $\Vv \in V$ and $t \in (0, 1]$ such that $\Va = (t\Vp, \Vv)$,
$\Vb = (t\Vq, \Vv)$ and the line segment $[\Vp, \Vq]$ is contained in the unit
circle of $\sK$.
Because $\pN{x}$ only depends on the ratios of the distances
between adjacent extreme points of the unit sphere, $\pN{x}$ holds in $\sK^d$
iff it holds in $\sK$ and the proof is complete.
\Done

\begin{Corollary}\label{thm:ns-ae-valid-undec}
Let $d \in \{2, 3, 4, \ldots\} \cup \{\infty\}$ and
let $\cal C$ be any class of normed spaces that includes all Banach
spaces of dimension $d$.
The set of $\fEA$ sentences that are valid in $\cal C$ is undecidable.
\end{Corollary}
\Proof
Just as in the proof of Theorem~\ref{thm:ms-ea-valid-undec}, use the existence of
a structure in which a purely existential formula defines the subset $\nat$ of
$\real$ to reduce the satisfiability of systems of Diophantine equations to
satisfiability in $\cal C$.
\Done

Our next undecidability result concerns $\fAE$ sentences in
theories of normed spaces.
As in the proof of Theorem~\ref{thm:ms-ea-valid-undec}, given a
quantifier-free formula $\psi(x_1, \ldots, x_k)$ in the language of arithmetic, we
will exhibit an $\fEA$ sentence $\psi_1$ in the language of normed spaces
that is satisfiable iff $\psi(x_1, \ldots, x_k)$ is satisfiable in $\nat$.
However, the quantifier structure of the sentence $\pPeano$ no longer
suits our purposes. Instead, we will design $\psi_1$ so that
its models comprise spaces whose unit circle contains a representation
of a periodic function on the set $\real_{+}$ of positive real numbers, which
we then use to define $\nat$.
We begin by showing how a norm may be used to define a function
on $\real_{+}$.

Consider a 2-dimensional normed space with a basis $\Ve_1$ and $\Ve_2$
where $\Norm{\Ve_1} = 1$.  For $x \in (-1, 1)$, we have $\Norm{x\Ve_1} =
|x| < 1$ and so by the remarks of Section~\ref{subsubsec:normed-spaces}, the
vertical line $\{x\Ve_1 + y\Ve_2 \mid y \in \real\}$ must meet the unit circle
in exactly two points, one in the upper half-plane and one in the lower.  Thus
with respect to the given basis, the part of the unit circle lying above the
open line segment $(-\Ve_1, \Ve_1)$ forms the graph of a function.

\begin{figure}
\begin{center}
\includegraphics[angle=0,scale=0.4]{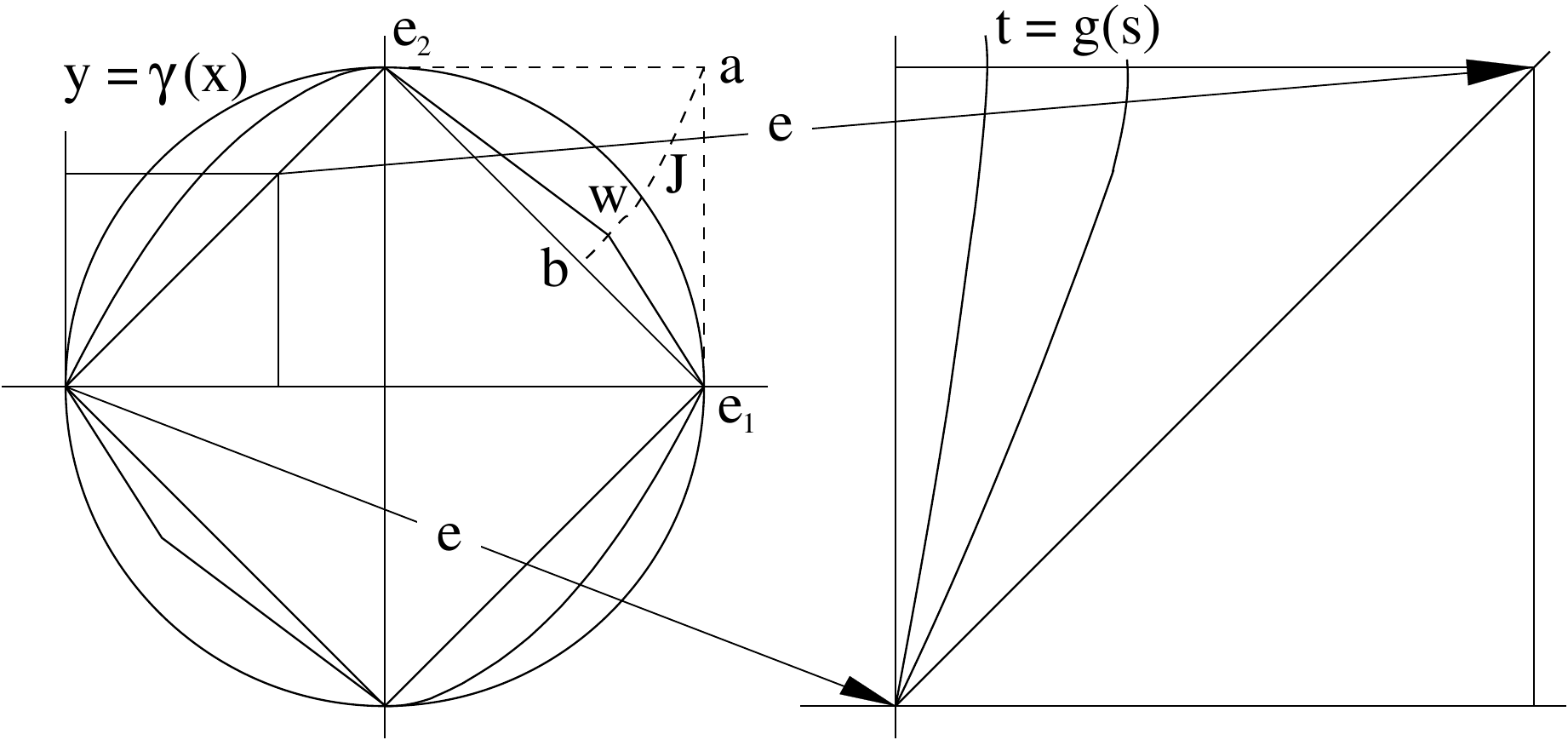}
\caption{The functions defined by $\pGamma{s, t}$ and $\pG{s, t}$ under various norms on $\real^2$}
\label{fig:l}
\end{center}
\end{figure}

Now consider the following formulas in which the vector variables $\Ve_1$ and
$\Ve_2$ occur free in addition to the scalar variables listed as parameters:
\begin{eqnarray*}
\pGamma{x, y} &\IsDef& -1 < x < 0 \And 0 < y < 1 \And \Norm{x\Ve_1 + y\Ve_2} = 1\\
\pG{s, t} &\IsDef& s > 0 \And t > 0 \And \Norm{{-}(1+t)\Ve_1 + (1+s)t\Ve_2} = (1+s)(1+t).
\end{eqnarray*}

Thus, for $s, t \not= -1$,  $\pG{s, t}$ is equivalent to $\pGamma{\frac{-1}{1+s}, \frac{t}{1+t}}$.
Assume that $\Ve_1$ and $\Ve_2$ are vectors
in some normed space $V$ and that the
following condition holds for all $x$ and $y$:
\begin{eqnarray*}
\pDef &\IsDef&
        \Norm{\Ve_1} = \Norm{\Ve_2} = 1 \And {}\\
    && (x\Ve_1 + y\Ve_2 = \VO \Imp x = y = 0) \And {}\\
    && (|x| > 0 \And |y| > 0 \And \Norm{x\Ve_1 + y\Ve_2} = 1 \Imp |x| < 1 \And |y| < 1).
\end{eqnarray*}
So $\Ve_1$ and $\Ve_2$ are unit vectors spanning a 2-dimensional subspace
of $V$ and, when we use them to define coordinates in that subspace,
the unit circle is contained in  the square with diagonal $[-\Ve_1-\Ve_2, \Ve_1+\Ve_2]$
and meets its boundary in the four points $\pm \Ve_1$ and $\pm \Ve_2$.
This means that $\pGamma{x, y}$ will hold iff $y = \gamma(x)$ where $\gamma$ is
the function whose graph comprises the part of the unit circle lying strictly
above the open line segment $(-\Ve_1, \VO)$,
while $\pG{s, t}$ will hold iff $t = g(s)$ where the graph
of $g$ is the image of the north-west quadrant of the unit circle under the
continuous bijection $e : (-1, 0) \times (0, 1) \rightarrow \real_{+} \times
\real_{+}$ defined by $e(s, t) = \left(\frac{1+s}{-s}, \frac{t}{1-t}\right)$.
The condition $\pDef$ ensures that $g(s)$ is well-defined
for all $s \in \real_{+}$.
Figure~\ref{fig:l} illustrates $\gamma$ and $g$ for various norms on $\real^2$.

\begin{Lemma}\label{lma:l-exists}
For some positive integer $M$, there is a 2-dimensional normed space $\sL$
containing vectors $\Ve_1$ and $\Ve_2$ for which $\pDef$ holds for
every $x$ and $y$, while
$\pG{s, t}$ holds iff $s > 0$ and $t = 2s + s^2 + \frac{1}{M}\Sin{s}$.
\end{Lemma}
\Proof
Define functions $g_r:\real_{+} \rightarrow \real_{+} $ for $2 > r > 0$
by $g_r(s) = 2s + s^2 + r\Sin{s}$.  Under the bijection $e$, $g_r$
corresponds to the function $\gamma_r: (-1, 0) \rightarrow (0, 1)$
where:
\[
\begin{array}{rcl}
\gamma_r(x) &=& \frac{g_r(\frac{x+1}{-x})}{1 + g_r(\frac{x+1}{-x})}.
\end{array}
\]
I claim that for all small enough $r>0$, $\gamma_r$ is a concave function,
i.e., the part of the plane lying below the graph of $\gamma_r$ is convex.
Assuming this, we can choose a positive integer $M$ such that $\gamma =
\gamma_{1/M}$ is concave. Noting that
$\gamma(x)$ tends to $0$ as $x$ tends to $-1$ and to $1$ from below as $x$ tends to $0$,
we can extend $\gamma$ to a concave function $\gamma^{*}$ on $[-1, 1]$ by taking $\gamma^{*}(-1) = 0$ and $\gamma^{*}(x) = 1 - x$ for $0 \le x \le 1$.
Let $\sL$ be $\real^2$ under the norm whose unit circle meets the upper
half-plane in the graph of $\gamma^{*}$.  Then in $\sL$, $\pDef$ holds for
every $x$ and $y$ and $\pGamma{x, y}$ defines $\gamma = \gamma_{1/M}$,  so that
$\pG{s, t}$ holds iff $s > 0$ and $t = g_{1/M}(s) = 2s + s^2 +
\frac{1}{M}\Sin{s}$.

It remains to prove that $\gamma_r$ is concave for all small enough $r>0$.
Certainly $\gamma_r$ is twice differentiable, and then, by standard results on
concave functions, it is sufficient to show that the second derivative
$\gamma_r''(x)$ is never positive for $x$ in $(-1, 0)$.
Differentiating the formula for $\gamma_r$ above twice gives:
\begin{eqnarray*}
\gamma_r'(x) &=& \frac{g_r'(\frac{x+1}{-x})}{x^2(1+g_r(\frac{x+1}{-x}))^2} \\
\gamma_r''(x) &=& \frac{g_r''(\frac{x+1}{-x})}{x^4(1 + g_r(\frac{x+1}{-x}))^2} +
        \frac{2g_r'(\frac{x+1}{-x})}{-x^3(1+g_r(\frac{x+1}{-x}))^2} -
        \frac{2(g_r'(\frac{x+1}{-x}))^2}{x^4(1+g_r(\frac{x+1}{-x}))^3}.
\end{eqnarray*}

Writing $s = \frac{x+1}{-x}$, so that $s > 0$ and $\frac{1}{-x} = 1 + s$,
and multiplying by the positive quantity $-x^3(1+g_r(s))^3$, we see that
$\gamma_r''(x)$ has the same sign as $h_r(s)$ where:
\begin{eqnarray*}
h_r(s) &=&
  (1+g_r(s))[(1+s)g_r''(s) + 2g_r'(s)] -
      2(1+s)(g_r'(s))^2  \\
  &=&
  (1+2s+s^2 + r\Sin{s})[(1+s)(2-r\Sin{s}) + 4 + 4s + 2r\Cos{s}] - \\
  && 2(1+s)(2+2s+r\Cos{s})^2.
\end{eqnarray*}
As $-1 \le \Sin{s}, \Cos{s} \le 1$, we have that $h_r(s) \le p(s, r)$, where:
\begin{eqnarray*}
p(s, r) &=&
  (1+2s+s^2 + r)[(1+s)(2+r) + 4 + 4s + 2r] - \\
  && 2(1+s)(2+2s-r)^2 \\
       &=&
   p_{0}(s) + p_{1}(s)r + p_{2}(s)r^2.
\end{eqnarray*}
each $p_i(s)$ being a polynomial of degree at most 3 in $s$ with constant
coefficients, say $p_i(s) = p_{i0} + p_{i1}s + p_{i2}s^2 + p_{i3}s^3$, $i =
0, 1, 2$.
Since $p_{0}(s) = p(s, 0) = 6(1+s)^3 - 8(1+s)^3 = -2 -6s -6s^2 -2s^3$,
each $p_{0j}$ is negative.  Let $q_j$ be the
coefficient of $s^j$ in $p(s, r)$ so $q_j = p_{0j} + p_{1j}r +
p_{2j}r^2$.  Since $p_{0j} < 0$, we may choose $\epsilon > 0$ such
that whenever $0 < r < \epsilon$, $q_j < 0$, $j = 0, 1, 2, 3$.  But then
if $0 < r < \epsilon$, we find that $p(s, r) < 0$ for all $s > 0$ whence
$\gamma_r''(x)$ is negative for all $x$ in $(-1, 0)$, since it has the
same sign as the quantity $h_r(s) \le p(s, r) < 0$. Thus $\gamma_r$ is
concave for $0 < r < \epsilon$ and the proof is complete.
\Done

Let the space $\sL$ and the positive integer $M$ be as given by the
lemma. In $\sL$, the following formula then defines the graph of the
positive half of the sine function when $\Ve_1$ and $\Ve_2$ are given
their usual interpretations:
\begin{eqnarray*}
\pSIN{s, t} &\IsDef& \pG{s, 2s + s^2 + \frac{1}{M}t}.
\end{eqnarray*}

Now consider the following formulas:
\begin{eqnarray*}
\pPeriodic &\IsDef&
        a > 0 \And {} \\
    &&  (0 < s < 2a \Imp (\pSIN{s, 0} \Iff s = a)) \And {} \\
    &&  (\pSIN{s, t} \Imp \pSIN{s + a, -t}) \\
\pN{x} &\IsDef& \pSIN{(x+1)a, 0}.
\end{eqnarray*}

In $\sL$ with the usual interpretation of $\Ve_1$ and $\Ve_2$, $\pPeriodic$
holds for all $s$ and $t$ if we interpret $a$ as $\pi$, in which case $\pN{x}$
holds iff $x$ is interpreted as a natural number.  On the other hand, if $V$ is
any normed space and there are $\Ve_1, \Ve_2 \in V$ and $a \in \real$ such that
$\pDef$ and $\pPeriodic$ hold for all $x$, $y$, $s$ and $t$, then $\pSIN{s, t}$
must define the graph of a function on $\real_{+}$ whose zeroes comprise
precisely the positive integer multiples of $a$, so that $\pN{x}$ defines the
natural numbers.

\begin{Theorem}\label{thm:ns-ea-valid-undec}
Let $d \in \{2, 3, 4, \ldots\} \cup \{\infty\}$ and
let $\cal C$ be any class of normed spaces that includes all Banach
spaces of dimension $d$.
The set of $\fAE$ sentences that are valid in $\cal C$ is
undecidable.
\end{Theorem}
\Proof
We will prove the equivalent claim that the set of $\fEA$
sentences that are satisfiable in $\cal C$ is undecidable.
Given a quantifier-free formula $\phi(x_1, \ldots, x_k)$ in the language
of arithmetic, define:
\begin{eqnarray*}
\phi_1 &\IsDef& \ex{\Ve_1\;\Ve_2\;a\;x_1 \ldots x_k}\all{x\;y\;s\;t} \\
    && \quad \pDef \And \pPeriodic \And \pN{x_1} \And \ldots \And \pN{x_k} \And \phi(x_1, \ldots x_k).
\end{eqnarray*}
Take $V = \sL \times W$ where $W$ is any vector space of dimension $d - 2$
under any norm extending that of the factor $\sL$;
if $\phi(x_1, \ldots, x_k)$ is satisfiable in $\nat$, $\phi_1$ will be satisfiable
in $\sL$ and hence in $V$.
Conversely, if $\phi_1$ is satisfiable in some normed space,
then under a satisfying assignment, the conditions
$\pDef$ and $\pPeriodic$ mean that $\pN{x}$ must define $\nat$
in $V$, so $\phi(x_1, \ldots, x_k)$ is satisfiable in $\nat$.
So, just as in the proof of Theorem~\ref{thm:ms-ea-valid-undec},
the existence of a decision procedure for $\fEA$ sentences satisfiable in $\cal
C$ would contradict the undecidability of satisfiability for quantifier-free
formulas in arithmetic.
\Done

To state our final result on undecidability, let us say a sentence is
$\fAIA$ if it has the form $A \Imp B$ where $A$ and $B$ are purely
universal.  With a small adjustment to the construction used to prove
Theorem~\ref{thm:ns-ea-valid-undec}, we now show that validity for
$\fAIA$ sentences is undecidable.  As $\fAIA$ sentences have both $\fEA$
and $\fAE$ equivalents, this provides an alternative proof for both
Corollary~\ref{thm:ns-ae-valid-undec} and
Theorem~\ref{thm:ns-ea-valid-undec}.

\begin{Theorem}\label{thm:ns-aia-valid-undec}
Let $d \in \{2, 3, 4, \ldots\} \cup \{\infty\}$ and
let $\cal C$ be any class of normed spaces that includes all Banach
spaces of dimension $d$.
The set of $\fAIA$ sentences that are valid in $\cal C$ is
undecidable.
\end{Theorem}
\Proof
If $d = 2$, let $\sL$ be the normed space constructed in the proof of
Lemma~\ref{lma:l-exists}.
Using the $\sL$-norm, let $C = \{\Vw \mid
\Norm{\Ve_2-\Vw} = 1\}$ be the unit circle centred at $\Ve_2$ and
consider the intersection $J = C \cap T$, where $T$ is the triangle with
vertices $\Ve_1$, $\Ve_2$ and $\Va = \Ve_1 + \Ve_2$  (see Figure~\ref{fig:l}).
$J$ meets the perimeter
of $T$ at the vertex $\Va$ with $\Norm{\Ve_1 - \Va}
= 1$ and at a point $\Vb$ on the edge $[\Ve_1, \Ve_2]$ with $t =
\Norm{\Ve_1 - \Vb} < 1$ (since $1 < \Norm{\Ve_2 - \Ve_1} < 2$).  $J$
is a continuous curve and so $\Norm{\Ve_1 - \Vw}$ takes on all values in
$[t, 1]$ as $\Vw$ ranges over $J$.  Since $T$ meets the unit disc of
$\sL$ in the edge $[\Ve_1, \Ve_2]$, it follows that there are
$i/j \in \rat$ and $\Vw \in J \subseteq T$ such that
$\Norm{\Ve_2 - \Vw} = 1 < \Norm{\Vw} $ and $ \Norm{\Ve_1 -
\Vw} = i/j < 1$.  Let $\sLO$ be the normed space whose unit disc is the
symmetric convex hull of the unit disc of $\sL$ and such a $\Vw$.
Then the $\sL$-norm and the $\sLO$-norm agree in the north-west quadrant and
so $\sL$ and $\sLO$  define the same functions $\gamma$ and $g$ and assign the
same lengths to the line segments that make up the north-east
quadrant of the unit circle of $\sLO$.  In $\sLO$, the
following formula holds iff $\Vp = s\Ve_1$, $\Vq = s\Ve_2$ and $\Vr =
s\Vw$ where $s = \pm 1$.
\begin{eqnarray*}
\pW{\Vp, \Vq, \Vr} &\IsDef&
    \Norm{\Vp} = \Norm{\Vq} = \Norm{\Vr} = \Norm{(\Vp + \Vr)/2} =
\Norm{(\Vq + \Vr)/2} = 1 \And {}\\
 && \Norm{\Vp - \Vr} = i/j \And \Norm{\Vq - \Vr} = 1 \And \Norm{(\Vp +
\Vq)/2} < 1
\end{eqnarray*}
\noindent Also, the following formula is invariant under
$\Vv \mapsto -\Vv$ and, when the free variables $\Ve_1$ and
$\Ve_2$ are given their usual interpretation in $\sLO$, holds iff $x = \pi$.
\begin{eqnarray*}
\pPi{x} &\IsDef& x < 4 \And \pSIN{x, 0}
\end{eqnarray*}
Now, given a quantifier-free formula $\phi(x_1, \ldots, x_k)$ in the
language of arithmetic,
define sentences $\psi$ and $\rho$ as follows:
\begin{eqnarray*}
\psi &\IsDef& \all{\Ve_1\;\Ve_2\;\Vw\;a\;x\;y\;s\;t}
\pW{\Ve_1, \Ve_2, \Vw} \And \pPi{a} \Imp \pDef \And \pPeriodic \\
\rho &\IsDef& \all{\Ve_1\;\Ve_2\;\Vw\;a\;x_1 \ldots x_k} \\
 && \quad \pW{\Ve_1, \Ve_2, \Vw} \And \pPi{a} \And \pN{x_1} \And \ldots \pN{x_k} \Imp \Not \phi(x_1, \ldots, x_k).
\end{eqnarray*}
By the above remarks on $\pW{\Vp, \Vq, \Vr}$ and $\pPi{x}$,
$\psi$ holds and $\pW{\Ve_1, \Ve_2, \Vw} \And \pPi{a}$ is
satisfiable in $\sLO$.
Also, in any normed space in which $\psi$ holds, $\pN{x}$ is true under an
assignment that satisfies $\pW{\Ve_1, \Ve_2, \Vw} \And \pPi{a}$ iff $x
\in \nat$. Thus if $\psi$ holds and $\pW{\Ve_1, \Ve_2, \Vw} \And \pPi{a}$ is satisfiable, then $\rho$ holds iff $\phi(x_1, \ldots, x_k)$
is unsatisfiable in $\nat$.  Thus $\psi \Imp \rho$ is valid in
a class of spaces including $\sLO$ iff $\phi(x_1, \ldots, x_k)$
is unsatisfiable and so a decision procedure for $\fAIA$ sentences that are
valid in such a class would lead to a
decision procedure for satisfiable quantifier-free sentences of
arithmetic, which is impossible.

For $d > 2$,  let $V$ be a Hilbert space of dimension $d - 2$, let
$W$ be the 2-sum of $\sLO$ and $V$, and identify $\sLO$ with the subspace
$\sLO \times 0$ of $W$. As in the proof of
Theorem~\ref{thm:ns-ae-def-nat}, if a line segment $[\Vu, \Vv]$
lies in the unit sphere of $W$, then it is parallel to $\sLO$.
Moreover, if also $\Norm{\Vu-\Vv} = 1$ then we must have that $\{\Vu,
\Vv\} = \pm\{\Ve_2, \Vw\} \subseteq \sLO$. This means that the formula $\pW{\Vp,
\Vq, \Vr}$ defines the same set of triples in $W$ as it does in $\sLO$.
The argument for $d = 2$ then shows that validity of $\fAIA$ sentences
in any class of spaces including $W$ is undecidable.
\Done

\Section{Related work and concluding remarks}\label{sec:concluding-remarks}

The reduction of second-order arithmetic to the theory of the real numbers
augmented with a predicate symbol for the integers has been known since the
1960s if not before. In descriptive set theory, the main ideas of
Section~\ref{sec:prelim} are used to show that a subset of $\real^n$ is
projective iff it is definable in the theory of the real numbers augmented with
a predicate for the integers; see for example, \citet{Moschovakis-descriptive},
Theorem 8B.4 or \citet{Kechris-cdst}, ex.~37.6. However, we know of no
published account of these ideas applied to problems of decidability.

\citet{scott-dimension} considers {\em geometric relations}, i.e.,
relations such as ``equidistant'' that are defined on affine euclidean
space in all dimensions and that are invariant under isometric
embeddings. He works with single-sorted languages whose variables
range over points and shows that a first-order sentence with $k + 1$ distinct
variables holds for every interpretation of its relation symbols as
geometric relations iff it holds in dimension $k$. This is clearly
closely related to our Theorem~\ref{thm:stability}.  (Scott's $k+1$ is
our $k$ because the constant vector $\VO$ costs us one variable.) He applies his
result to a formulation of euclidean geometry as a single-sorted theory with
``between'' and ``equidistant'' as primitive predicates and obtains
decidability and related results for theories $\cal E$, ${\cal E}_m$ and
${\cal E}_{\infty}$ analogous to our $\IP{\relax}$, $\IP{m}$ and
$\IP{\infty}$. Scott's proofs are based on semantic considerations
that apply to all geometric relations, in contrast with our more
algorithmic approach via a syntactic quantifier elimination
procedure. In his single-sorted Tarski-style formalism, the emphasis
is on geometry and the real numbers only arise implicitly as
equivalence classes for the equidistance relation, while our
two-sorted approach is closer to typical mathematical and engineering practice.

Before both our work and that of \citet{kurz-metric},
\citet{bondi:dec-ms} had proved the
undecidability of the theory of metric spaces. Let $\LB$ be the language of a
single binary predicate, intended to be interpreted as a symmetric relation.
Translated into our two-sorted framework, Bondi's proof defines a mapping
$\sigma \mapsto \sigma_0$ from sentences in $\LB$ to sentences in the language
$\LM$ of metric spaces  such that that $\sigma_0$ is valid if $\sigma$ is valid
and $\Not\sigma_0$ is valid if $\sigma$ is finitely refutable.  Undecidability
follows from a theorem of Lavrov. The only metric spaces Bondi uses are finite
subspaces of euclidean space, so her methods show that the theory of such
metric spaces is undecidable.  The methods of the present
paper give more information about the many-one degree of the theory and show
that the additive and the $\fEA$ fragments are undecidable.  As far as we know,
our Theorem~\ref{thm:ea-v-decidable} giving the decidability of the $\fAE$
fragment is the strongest known positive result on decision procedures for
theories of metric spaces.

\citet{bondi:dec-ns-hs} considers normed spaces and inner product spaces (she
actually writes ``Hilbert spaces'', but metric completeness plays no r\^{o}le
in her proofs). She first proves the undecidability of the theory of normed
spaces by a proof similar in structure to her proof for metric spaces sketched above. The method of the proof actually gives the
undecidability of any class of normed spaces containing all finite-dimensional
spaces, whereas the methods of the present paper give undecidability of any
class of spaces containing all spaces of any given (possibly infinite)
dimension $d > 1$.  Bondi then turns to the decision problem for inner product
spaces: she first gives recursive axiomatizations of the theory $T$ of all
non-trivial inner product spaces, of the theory $T_0$ of all
infinite-dimensional inner product spaces and of the theories $T_n$ of
$n$-dimensional inner product spaces, $n = 1, 2 \ldots$
She shows that the theory $T_0$ is model complete and so complete, since
any two models of $T_0$ contain isomorphic submodels.
Being recursively axiomatizable and complete, $T_0$ is
therefore decidable.
This argument shows the correctness of a decision procedure that
enumerates proofs rather than the more efficient procedures of
section~\ref{sec:qelim} above.
Writing $\Tfin$ for the theory of non-trivial
finite-dimensional inner product spaces, Bondi goes on to argue that $T
= \Tfin$ and concludes using a lemma of Ershov that $T$ is decidable.
Unfortunately, there is a significant gap in her proof that $T = \Tfin$:
she claims that a certain sentence in $T_0$ must belong to $T_n$ for
sufficiently large $n$, but gives no proof of this.  Her claim is true,
but it is unclear how to prove it without appealing to
Theorem~\ref{thm:special-formula} from the present paper. A precisely
analogous situation in which the analogue of $T = \Tfin$ fails can be
reached by adding a predicate $\pD(x)$ on scalars with the intended
interpretation that $\pD(x)$ hold in $V$ iff $x \le \Dim{V}$. This gives
extensions $T'$, $T'_0$ etc. of the theories $T$, $T_0$ etc. $\pD(x)$
can be defined by a recursive set of axioms and $T'_0$ and the $T'_n$
can be seen to be complete using Bondi's arguments.  However, $T' \not=
\Tfin'$, since $\ex{x}\Not \pD(x)$ holds in an inner product space $V$
iff $V$ is finite-dimensional.

A vector space over the real field is a special case of a module
over a ring.  Theories of modules over rings have been widely studied, often
with a view to applications in algebra.  However, most of this work has
concentrated on single-sorted theories in which quantification over the ring of
scalars is not allowed, the action of the ring on the module being represented
by function symbols $f_x$ indexed by ring elements such
that $f_x(\Vv) = x\Vv$ in the intended interpretations.
With this formulation, the
procedure of Baur and Monk gives quantifier elimination
relative to a set of predicates that specialise to our dimension predicates
$\pD_n$ when the ring is a field.
This procedure provides a powerful theoretical tool; see for
example, \citet{Prest}.  For modules over any B\'{e}zout domain,
\citet{vandendries-holly} give a quantifier elimination procedure
for formulas in which free scalar variables are allowed.  Their method is via a
reduction to the single-sorted language over a ring of polynomials and it is
unclear how it could be generalised to deal with scalar quantification.

\citet{granger-thesis} considers the theory of vector spaces equipped with a
bilinear form and proves a form of quantifier elimination for the natural
two-sorted formulations using model-theoretic arguments.  He gives an
interesting discussion of two-sorted formulations that attempt to decouple the
model theory of the underlying field from the model theory of vector spaces
over it. These formulations lie somewhere between the single-sorted formulation
of the Baur-Monk theorem and the two-sorted formulation adopted in the present
work. Granger's conclusion is that such a decouplng is in some sense not possible.

As already mentioned in Section~\ref{sec:ea-ae-fragments}, our results on the
undecidability of the $\fAE$ and $\fEA$ fragments of the theory of normed
spaces can be strengthened to the additive case.
\citet{arthan:additive-aia} does this by adapting the constructions used here
to prove Theorem~\ref{thm:ns-aia-valid-undec} so that scalar multiplication
becomes definable via a purely existential formula.

In the present paper, we have focussed on the case when the field of scalars
comprises the real numbers. However, all the results on inner product spaces
and Hilbert spaces in Section~\ref{sec:qelim} go through with the proofs
unchanged for an arbitrary real closed field. As discussed at the end of
Section~\ref{sec:decidable-fragments}, our positive decidability result for the
universal fragment of the theory of real normed spaces can also be adapted to
cover normed spaces over any real closed field.

One cannot hope to reduce second order arithmetic to a recursively
axiomatizable theory like the theory of normed spaces over a real
closed field. However, \cite{arthan:any-field} gives a construction over an
arbitrary ordered field of a 2-dimensional normed space that encodes the graph
of natural number multiplication.  Via a reduction of Robinson's theory $Q$,
this gives the undecidability of the additive theory $\NSA{\relax}({\cal C})$,
and hence, {\it a fortiori}, the full theory $\NS{\relax}({\cal C})$ for normed
spaces over any non-empty class of ordered fields $\cal C$  and similarly for
the theories $\NSA{\infty}({\cal C})$, $\NSA{n}({\cal C})$, $1 < n \in \nat$,
and $\NSA{\fin}({\cal C})$ with the indicated constraints on dimension.

\citet{Kopperman-ms, Kopperman-hs} considers formalisations of Hilbert spaces
and metric spaces in a family of infinitary languages, $L^t_{\pi,\varepsilon}$,
where $t$ amounts to a many-sorted signature and $\pi$ and $\varepsilon$ are
cardinals. $L^t_{\pi,\varepsilon}$ has $\pi + \epsilon$ variables and admits
conjunction and disjunction of any set $X$ of formulas where $\Card{X} < \pi$
and quantification over any set $Y$ of variables where $\Card{Y} < \epsilon$
(so the usual finitary language over a signature $t$ is $L^t_{\omega,\omega}$).
For $t$ a signature appropriate for Hilbert spaces, he gives a result for
$L^t_{\omega_1,\omega_1}$, redolent of our
Corollary~\ref{:ip-axiomatizability}, stating that any formula is equivalent to
a boolean combination of sentences $D_n$ asserting the dimension is $n \in \nat
\cup \{\infty\}$.  However, the infinitary languages are much more expressive
than the languages we consider: in the case of separable metric spaces, one can
encode the full metric structure of a countable dense subset in a single
sentence. Thus, by contrast with our undecidability results, Kopperman proves
quantifier elimination for separable complete metric spaces relative to a set
of formulas that define the possible countable dense subsets.  Unsurprisingly,
Kopperman's methods of proof are quite different from those of the present
work.

Special languages and logics for Banach spaces and similar structures have been
widely studied, largely from a model-theoretic perspective, with
applications in functional analysis in mind; see, for example,
\citet{Henson-Iovino}.  This work has typically involved logics that are weaker
than full first-order logic, since metric completeness makes conventional model
theory for Banach spaces less useful in the intended applications.
\citet{shelah-stern} have demonstrated the problems with conventional
model theory in this context using a construction with a similar flavour to our
construction of a sentence that holds in all Banach spaces but is not valid in
all normed spaces.

The work reported in the present paper was motivated by an interest in applying
mechanized theorem-proving to problems in pure mathematics and engineering.
For the potential applications, vector spaces and inner product spaces over the
real field are important, and, as we have seen, they admit more powerful
decision procedures than modules over an arbitrary ring. However, the
complexity of these decision procedures and the undecidability of theories of
normed spaces present some interesting challenges.

\section*{Acknowledgments}

We thank the referee for a very constructive and thorough report,
for helpful advice on the notation and
presentation, particularly in Section~\ref{sec:qelim}, and for providing a
neat improvement to the proof of Theorem~\ref{thm:ea-v-decidable}.  We
are grateful to Angus MacIntyre, Dana Scott and the referee for valuable
pointers to the literature. Finally, we are deeply indebted to Nadya
Kuzmina for helping us to follow the referee's pointers by kindly
undertaking to translate the papers by Irina Bondi from the original
Russian.

Arthan's work was supported in part by UK EPSRC grant EP/F02309X/1.
\bibliographystyle{elsarticle-harv}

\bibliography{all}


\end{document}